\input amstex
\documentstyle{amsppt}
\input epsf.tex
\input label.def

\MakeToc{toc}

\def\conj{\operatorname{conj}}
\def\C{{\Bbb C}}
\def\R{{\Bbb R}}
\def\Z{{\Bbb Z}}
\def\Q{{\Bbb Q}}
\def\A{\Bbb A}
\def\DD{\Bbb D}
\def\E{\Bbb E}
\def\U{\Bbb U}
\def\Rp#1{\R\roman P^{#1}}

\def\Im{\mathop{\roman{Im}}\nolimits}
\def\Re{\mathop{\roman{Re}}\nolimits}

\def\la{\langle}
\def\ra{\rangle}
\def\id{\mathop{\roman{id}}\nolimits}
\def\Aut{\mathop{\roman{Aut}}\nolimits}
\def\emptyset{\varnothing}
\def\oo{\varnothing}
\def\discr{\operatorname{discr}}
\def\rank{\operatorname{rank}}
\def\D{\Delta}
\def\til{\widetilde}
\def\s{\sigma}
\def\e{\varepsilon}
\def\a{\alpha}
\def\b{\beta}

\def\p{\pi}
\def\om{\Omega}
\def\sm{\smallsetminus}

\def\Fix{\operatorname{Fix}}
\redefine\P{\Cal P}
\def\RR{\Cal A}
\def\ss{\frac{\s}3}
\def\ID{\roman{ID}}
\def\Ch{\roman{Ch}}
\def\di#1#2#3{\la#1\frac#2#3\ra}
\def\dis#1#2{\la\frac#1#2\ra}
\def\dip{\la\frac23\ra}
\def\din{\la-\frac23\ra}
\def\co{\roman{comp}}
\def\bb#1#2{b^{#1}_{#2}}
\def\uu{\Cal U_2}
\def\vv{\Cal V_2}

\def\n{\nu}
\def\Br{\roman{Br}}
\def\h{\varkappa}
\def\ord{\operatorname{ord}}

\def\q{\frak q}
\def\boxedR{$\boxdot$}
\def\L{\Bbb L}
\def\h{\roman{h}}
\def\LDh{(\L,\bold{\D},\h)}
\def\LDhc{(\L,\bold{\D},\h,c)}

\let\rk=\remark
\let\endrk=\endremark
\let\ge\geqslant
\let\le\leqslant
\let\+\sqcup
\pagewidth{5.8 true in}
\pageheight{8.6 true in}

\nologo
\NoBlackBoxes
\leftheadtext{S.~Finashin, V.~Kharlamov}
 \topmatter
\title
Apparent contours of nonsingular real cubic surfaces
\endtitle
\author Sergey Finashin, Viatcheslav Kharlamov
\endauthor
\address Middle East Technical University,
Department of Mathematics\endgraf Ankara 06800 Turkey
\endaddress
\address
Universit\'{e} de Strasbourg et IRMA (CNRS)\endgraf 7 rue Ren\'{e}-Descartes 67084 Strasbourg Cedex, France
\endaddress
\thanks{The second author acknowledges a financial support by
the grant ANR-09-BLAN-0039-01 of {\it Agence Nationale de la Recherche}.}\endthanks
\abstract\nofrills We give a complete deformation classification of real Zariski sextics, that is of
generic apparent contours of nonsingular real cubic surfaces.
As a by-product, we
observe a certain "reversion" duality in the set of deformation classes of these sextics.
\endabstract
\contents
\endcontents
\endtopmatter
\newpage

\document
{\sl \hskip1.8in ``Tout n'est qu'apparence, non ?''}

\hskip2.5in   Alberto Giacometti

\section{Introduction}

\subsection{The main problem and principal results}
There exist many ways to visualize cubic surfaces. One of those
that take into account not only the internal geometry of the
surface but also its position in the space, consists in "viewing
the surface from an external point", that is in considering
central projections of the surface onto a plane and thus
representing the surface as a 3-fold covering of the plane
branched over a certain curve, called the \emph{apparent contour} of the
surface. If the cubic surface is nonsingular, such  an apparent contour is a curve of degree 6 whose
singular locus, for a generic choice of the center of projection,
is formed by six cusps lying on a conic. Here, the condition of
``generic choice'' means certain transversality of the cubic
surface to its first two polars with respect to the selected
point, see Section \ref{complex-contours} (in fact, the same kind
of generic choice assumption is specific to the literature on
Chisini conjecture, see, for example,  \cite{Mo} and
\cite{KulCh}).

 The sextics with six cusps
lying on a conic are called below {\it Zariski sextics}, following
the  ``Arnold principle'' \cite{Arnold-teaching}.
 The converse statement that any such a sextic is an
apparent contour of some nonsingular cubic surface with respect to a generic
center of projection is a classical observation
(see \cite{Salmon},
XIV.445 for the direct statement, and \cite{Zar29},\cite{ZarAS},\cite{Segre} for the converse one).
This correspondence establishes an isomorphism
between the space of projective classes of pairs formed by a nonsingular cubic
surface with a generic center of projection, and the space of
projective classes of Zariski sextics. Therefore, studying
the former classes is reduced to studying the latter ones, and it
works equally well over $\C$ and over $\R$.

G.~Mikhalkin had undertaken an analysis of the apparent contours of
real nonsingular cubic surfaces and reported the results he obtained in \cite{M}.
Namely, he looked for a topological classification of the real apparent contours
that are enhanced by specifying the topological type of the real locus of the cubic surface
(such an enhancing can be expressed by coloring the part of the real
plane where the projection is three-to-one).
Mikhalkin listed 49 {\it enhanced isotopy classes} of apparent contours that
he constructed, mentioned 7 enhanced isotopy classes whose existence is
uncertain, and claimed that there are no others.
Our research resulted from attempts to understand
Mikhalkin's results, to find the proofs, to complete the
classification and, overall, to sharpen it by providing a classification up
to {\it equisingular deformations}.
Recall that for plane reduced curves, which is our case,
equisingular deformations can be defined in topological terms, as
continuous families of algebraic curves preserving the quantity of
singular points and their Milnor numbers. Thus,
in our setting the equisingular deformation classes
of the apparent contours of nonsingular real cubic surfaces
are nothing but the connected components
of the space of real Zariski sextics.

To obtain such a {\it deformation classification} of apparent
contours, we analyze the K3-surfaces that we obtain as double covers
of the plane branched along Zariski curves.
Kulikov's theorem on surjectivity of the period map for K3-surfaces
(see \cite{KulSurj}) allows us to reduce the deformation classification
to a classification of certain involutions on the K3-lattice
(so called "geometric involutions", see Section \ref{geometric-involutions} for precise definitions).
Then, using Nikulin's results on the arithmetics of integral lattices (see \cite{Nik1}),
and the extension of these results by Miranda-Morrison (\cite{MM1}, \cite{MM2})
we make the final classification explicit.

As a result, we prove that there exist precisely 68 deformation classes of
apparent contours. In turn, we find that 7 Mikhalkin's  uncertain enhanced isotopy classes
are actually realizable and that, moreover, there exist 6 more enhanced isotopy classes missing in his list;
the full list of enhanced isotopy classes contains 62 items.
The difference, also equal to 6, between the  number of deformation classes
and the number of enhanced isotopy classes is due to existence of 6 pairs-twins of
deformation classes, such that the twins in each pair give the
same isotopy class of apparent contours and the same topological type of the cubic surface, but differ by the complex type of the underlying Zariski sextic, which can be dividing or not dividing.

The final classification is presented at two levels.
At the first level, we establish a one-to-one correspondence
between the set of deformation classes and the set of conjugacy classes of
what we call ascending geometric involutions on the K3-lattice
(see Theorem \ref{deformation-classification})
and enumerate the latter conjugacy classes in terms of
the eigenlattices of involutions (see Theorem \ref{eigenlattices-enumeration}).
Such a classification can be viewed as a kind of "imaginary" one, since
it does not immediately disclose the topology of the corresponding apparent contours.
At the second level, we translate the information on lattices into
 a ``real'' information, which yields
a classification in terms of the ID's of the apparent contours (see Theorem \ref{main-sextics}), where
an ID is, roughly speaking, just an enhanced isotopy type with an additional
bit of information specifying whether a sextic is dividing or not.

\subsection{Partnership duality}
The list of the ID's of Zariski sextics that we give in Theorem \ref{main-sextics} is
organized so that it
emphasizes some duality resulted from the classification. This duality splits 62
 of the deformation classes into 31 pairs (so that only 6
remaining deformation classes are left without a pair).
Geometrically, in terms of arrangement of components (and their cusps) of a real Zariski
 sextic, this duality can be described as a certain reversion,
see details in Section \ref{reversion}.
 In Section \ref{explicit-partners}
 we give a ``conceptual'' explanation of this duality via
a ``wall crossing'' of special faces that can be found on
all but six exceptional fundamental period domains in the period space of our K3-surfaces.
Such special faces correspond to degeneration of Zariski sextics to a triple conic
by means of families of the form $Q^3 + tf_2Q^2+t^2f_4Q+t^3f_6$
(more precisely, of the form $Q^3+t(f_1(x,y,z) Q+ t f_3(x,y,z))^2=0$), and the
K3-surface resulting in the limit of such a family at $t=0$ is an
appropriate double cover of the ruled surface $\Sigma_4$, so that
the duality in terms of these $K3$ surfaces results in twisting
the real structure by the deck transformation of the covering, or
equivalently, in terms of the families, it results
in changing the sign of $t$, that is in switching to $Q^3 - tf_2Q^2+t^2f_4Q-t^3f_6$ (which
is equivalent to switching the sign of $Q$).

\subsection{Some related results}
The approach to studying the topology and deformation classes of real K3-surfaces and, in particular, of nonsingular real plane sextics via the period map was
developed by V.~Kharlamov \cite{Kh76} and V.~Nikulin \cite{Nik1} in
the 70th. Later on, in the 90th, the same approach was used by  I.~Itenberg \cite{Iten} for topological and deformation study of real plane sextics with one node.

Applications of a similar approach to complex K3-surfaces are much more abundant: here, we mention only the
works by T.~Urabe \cite{U}
on classification of configurations of simple singularities on complex plane sextics and
the recent works by
A.~Degtyarev \cite{Degt-Oka}-\cite{Degt-def} on deformation classification of complex plane sextics with simple singularities.

The Zariski sextics, which are the subject of our paper, belong to the so called class of plane curves of torus type: they are generic plane curves
of torus type of degree $(2,3)$. There exists a vast literature on the geometry of plane curves of torus type over the complex field. In particular, due to works by D.T.~Pho and M.~Oka a complete classification of configurations of singularities on complex sextics of torus type is actually known, see
\cite{OP} and \cite{O}.

\subsection{Contents of the paper} The paper is organized as follows. In Chapter \ref{Chapter_sextics}
we introduce Zariski sextics, analyze their relation to cubic surfaces,
and study the basic properties of Zariski sextics and of the double planes
branched along them. The chapter is concluded with our principal
deformation classification Theorem \ref{main-sextics}. In Chapter \ref{Chapter_Lattice_preliminaries}, which
is devoted to the arithmetics of integral lattices, we recall the basic definitions and
some well-known results on lattices (their discriminant finite forms,
gluing of lattices and involutions), slightly developing them and
making a special emphasis on the lattices that
have only discriminant factors $2$ and $3$, since it is
this kind of lattices that appears later on
in the proofs of the main results.
Chapter \ref{Chapter_CoveringK3} deals with our main object of investigation, the
K3-surfaces that are double coverings of the plane branched along Zariski sextics.
Here, we relate the topological and geometrical properties of these surfaces with the arithmetics of K3-lattices.
In particular, we introduce and study
{\it geometric involutions} on the $K3$-lattices
and associated with them certain {\it $T$-pairs} of eigenlattices.

Chapter \ref{Chapter_Tpairs} starts with a study of
the action of geometric involutions on the lattice generated by the exceptional divisors
and ends with proving the realizability of all the $T$-pairs via geometric involutions.
In Chapter \ref{Chapter_periods} we introduce the period space for K3-surfaces that are double coverings of the plane
branched along Zariski sextics, check that up to codimension two the periods of our K3-surfaces
fill out the fundamental domains of a group generated
by reflections, and deduce from that the bijection between the set of deformation classes of real Zariski sextics and
the set of conjugacy classes of ascending geometric involutions.
In Chapter \ref{mainProof} we classify the ascending $T$-pairs and prove
their stability. Finally, in Chapter \ref{Chapter_back} we
apply the results of preceding Chapters and enumerate
the geometric involutions, and hence
the deformation classes of real Zariski sextics, in terms of $T$-pairs. We
conclude with translating the classification into the language
of IDs of real Zariski sextics, which proves the deformation classification statement formulated in Chapter \ref{Chapter_sextics}. A few final remarks are collected in Chapter \ref{concluding}.

\subsection{Terminology conventions and notation} A {\it real algebraic variety} is always considered as a
pair $(X,c)$, where $X$ is a complex one and $c\:X\to X$ is an anti-holomorphic involution
called {\it the complex conjugation}, or {\it the real structure}.
The locus of complex points of $X$ is denoted by $X(\C)$,
and the real locus, $\Fix c$, by $X(\R)$.
In some cases, when it causes no confusion, we denote also by $c$
the induced involution in the homology.

When we speak on
singular points and on (equisingular) deformations, we take into account not only the
real points but the complex ones as well. For example,
from such a viewpoint, the real Zariski sextics having no real points can be not only different but even
non equivalent up to equisingular deformations; in
fact, such ``empty Zariski sextics'' form two distinct
deformation classes (a twin pair of dual classes in the sense mentioned above).

Working with homology or cohomology we use by default
$\Z$-coefficients, dropping them from the notation.
In the case of compact oriented even dimensional manifolds, especially in the case
of complex $K3$-surfaces, we identify the middle homology and cohomology lattices via Poincar\'e duality up to omitting, with a slight abuse of notation, the duality operator.

Several other conventions related to lattices are introduced in the beginning of Section 3.

The symbol \boxedR \,\, is used to mark the end of a remark.

\subsection{Acknowledgements}
We thank G.~Mikhalkin for sending to us his personal notes with the figures illustrating construction of the Zariski sextics listed in \cite{M}. This work was essentially done during the visits of the first author to
Strasbourg University and partially during our visits to MPIM (Bonn); we thank the both institutions for providing good working conditions. We thank also an unknown referee for many valuable suggestions.

\section{Zariski sextics}\label{Chapter_sextics}

\subsection{Generic projection in the complex setting}\label{complex-contours}
 Let $X$ be a non-singular cubic surface in $P^3$ defined by a
homogeneous cubic polynomial $f=f(x_0,\dots,x_3)$,
$\xi=[\xi_0:\dots:\xi_3]$ a point in
$P^3(\C)\sm X(\C)$, and $\p=\p_{X,\xi}\:X\to P^2$ the central
projection from $\xi$.
 The critical set of $\p$ is the curve
$B=B_{X,\xi}=X\cap X_\xi$ traced on $X$ by {\it the polar
quadric}, $X_\xi$, defined by $f_\xi=\sum_{i=0}^3\xi_i f_{x_i}$.
 We call this curve $B\subset X$ the  {\it rim-curve}.
The set of critical values, $A=\p(B)$, is called the {\it apparent
contour of $X$ with respect to $\xi$}.
The lines passing through $\xi$ and
intersecting $B$ trace on
$X$ another curve, $B'$, which we call the {\it shadow contour},
so that $\p^{-1}(A)=B\cup B'$.

We say that a point $\xi$ is {\it X-generic} if
$X_\xi$ is transverse to $X$, and $B$
is transverse to {\it the Hessian plane} $X_{\xi\xi}$ defined by
$f_{\xi\xi}=\sum_{i,j=0}^3\xi_i\xi_j f_{x_ix_j}$.
As is well known, the set of X-generic points $\xi$ is a non-empty Zariski-open subset of $P^3$ (one can find a detailed proof of the non-emptiness in \cite{CF}).

Let us choose a coordinate system so that $\xi$ turns into $[0:0:0:1]$
and $f$ appears in a depressed form (i.e., the quadratic
in $x_3$ term vanishes),
 $f=x_3^3+px_3+q$, where $p$ and $q$
are homogeneous polynomials in $x_0,x_1,x_2$ of degrees $2$ and
$3$ respectively. We say that such a coordinate system is {\it
associated with $X$ and $\xi$} (it is well-defined up to a
coordinate change in
the projection plane $P^2=\{[x_0:x_1:x_2]$).

\lemma\label{genericity} In a coordinate system associated with $X$ and $\xi$,
\roster\item the polar quadric $X_\xi$ is $\{3x_3^2+p=0\}$,
 \item the Hessian plane $X_{\xi\xi}$ is $\{x_3=0\}$,
 \item the apparent contour $A$ is
 the sextic
 defined by the
discriminant polynomial $D_f=4p^3+27q^2$ of $f$,
\item the shadow-contour $B'$ is given by equations $3x_3^2+4p=0, x_3^3+px_3+q=0$.
 \endroster
\endlemma

\proof Straightforward calculation.
\endproof

\corollary\label{generic contours} Assume that $\xi$ is X-generic.
Then, \roster \item  the rim curve $B$ and
the shadow-contour $B'$ are nonsingular and intersect each other
at the 6 points $X\cap X_\xi\cap X_{\xi\xi}$ with multiplicity 2
(i.e., $B$ and $B'$ have simple tangency at these points);
 \item the apparent contour $A$ is smooth except the
6 cuspidal points at $\p(X\cap X_\xi\cap X_{\xi\xi})$.
\qed\endroster
\endcorollary

\corollary The point $\xi$ is X-generic if and
only if the conic $p$ and cubic $q$ intersect transversely, and
the sextic $A$ has no other singularities except the 6 cusps at
$p=q=0$. \qed\endcorollary

If a sextic has six cusps
lying on a conic
and no other singular points, we call it {\it Zariski sextic}. The following fact is
also well known; its proof can be found in \cite{Segre}.

\proposition\label{torus-type} A plane sextic is a Zariski sextic
if an only if it is the apparent contour of a nonsingular cubic
surface $X$ with respect to an $X$-generic point. In particular,
each Zariski sextic can be presented by equation $4p^3+27q^2$, where $p$ and $q$
are homogeneous polynomials of degree $2$ and $3$ defining a conic and a cubic
intersecting transversely. Such a presentation is unique up to rescaling $(p,q)\mapsto (t^2p,
t^3q)$. \qed \endproposition

\corollary\label{central-correspondence} The central projection correspondence
that associates to a cubic surface $x_3^3+px_3+q=0$ and the point $[0:0:0:1]$
the sextic $4p^3+27q^2=0$ provides
a homeomorphism between the space of projective classes of pairs $(X, \xi)$,
where $X$ is a nonsingular cubic surface and  a point $\xi$ is $X$-generic, and the space
of projective classes of Zariski sextics.
\endcorollary

\proof
Bijectivity at the level of projective classes follows from Proposition \ref{torus-type}.
To construct continuous local inverse maps, it is sufficient for every
Zariski sextic $A$  to represent
a neighborhood of the $PGL(3)$-orbit of $A$
as $(G\times E)/G_A$, where $G=PGL(3)$, $G_A=Aut A$, and $E$ is a finite dimensional vector space with a linear
action of $G_A$ on it.
Afterwards, it is enough
to choose along $S=1\times E\subset (G\times E)/G_A$ the family of equations
$4p_s^3+27q_s^2=0, s\in S,$ in a way that $p_s,q_s$ depend on $s\in S$ continuously and equivariantly with respect to $G_A$.
The latter property can be achieved by simple averaging over the action of  the group $G_A$ (recall that this group is finite).\endproof

\subsection{Generic projection in the real setting}
 From now on, we restrict our considerations to cubic surfaces $X$ defined
over $\R$ and points $\xi$ in $P^3(\R)\sm X(\R)$.
However, we assume everywhere that the chosen point $\xi$ is
X-generic over $\C$, in the sense of Section \ref{complex-contours}.

Since $X,\xi$ are real, their Zariski sextic is real as well, and
vice versa. More precisely, Proposition \ref{torus-type} and
Corollary \ref{central-correspondence} imply the following.

\proposition\label{real-central correspondence} The central
projection correspondence provides an isomorphism between the three spaces:
\roster
 \item the set of real projective classes of real Zariski sextics,
 \item the set of pairs $(p,q)$ where  $p$ and $q$ are
real homogeneous polynomials of degree $2$ and $3$, respectively,
considered up to simultaneous real projective transformations of
variables and satisfying two restrictions: first, the conic $p=0$ and the cubic
$q=0$ intersect transversally, and, second, the sextic $4p^3+27q^2$
has no other singularities than 6 cusps;
 \item the set of  real projective classes of pairs $(X,
\xi)$, where $X$ stands for a nonsingular real cubic surface and
$\xi$ for a real X-generic point. \qed
\endroster
\endproposition

\rk{Remark} Recall that according to our convention, transversality is assumed at all the complex points,
and that some of the six cusps of a real Zariski curve may be imaginary.
\boxedR\endrk

 Let $\P_\pm$ denote the region $\{\pm p\le0\}$
bounded in $P^2(\R)$ by our real conic $p$.
Next, $P^2(\R)$ is divided into two ``halves'' bounded by
$A(\R)$, namely $\RR_\pm=\{[x]\in P^2(\R)\,|\,\pm D_f(x)\ge0\}$.
Note that a Zariski sextic determines uniquely the sign of its degree $6$ polynomial
 $D_f=4p^3+27q^2$, as well as the sign of $p$, and therefore
determines the signs of the above regions, $\P_\pm$ and $\RR_\pm$.

At a neighborhood of a real cusp the curve divides the real plane
into an ``acute region'' between the branches of the curves
forming angle zero at the cusp, and the complementary ``reflexive region'', so that
we can speak on {\it the acute and the reflexive sides} of a curve at its real cusp.
The definitions immediately imply the following.

\lemma\label{o_n}
\roster
\item
The real part of the apparent contour, $A(\R)$, lies entirely in
$\P_+$. \item The projection $p_\R=p|_{X(\R)}$ is three-to-one
over the interior points of $\RR_-$ and one-to-one over the
interior of $\RR_+$. \item The region $\RR_-$ lies entirely in
$\P_+$ and bounds the real cusps of $A$ from the acute side,
see Figure \ref{Zariski-example}.\qed
\endroster
\endlemma

\midinsert
\epsfbox{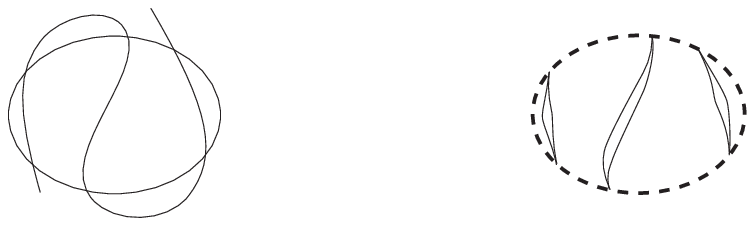}
\figure{Zariski sextic  on the right is obtained by a small perturbation $4(\epsilon p)^3 + 27 q^2=0$
of the conic $p=0$ and cubic $q=0$ shown on the left; $p<0$ inside the conic and $\epsilon>0$.}
\label{Zariski-example}
\endfigure
\endinsert


Note that one of the regions $\RR_\pm$ is orientable, and the
other is not; let us denote them respectively $\RR_o$ and $\RR_n$,
thinking of $o$ and $n$ as functions of
$(X,\xi)$
taking opposite
values $o,n\in\{+,-\}$.
Then, in accordance with Lemma \ref{o_n}, $o=+$ if and only if the projection $p_\R=p|_{X(\R)}$ is three-to-one
over the interior points of the non-orientable component of $P^2(\R)\sm A(\R)$. Due to Proposition \ref{torus-type},
these functions $o$ and $n$  can be also considered as
functions of $A$.

 The connected components of $A(\R)$ will be called {\it ovals}.
 An oval may have real cusps; it is called {\it cuspidal} in this case
 and {\it smooth} otherwise.

\subsection{Restrictions on the arrangements of ovals and real cusps}

 \lemma\label{cusp-parity}
Every cuspidal oval of a real Zariski sextic has an even number of cusps.
 \endlemma

\proof Every component of the real rim-curve $B(\R)$ is two-sided on $X(\R)$
(since the polar quadric $X_\xi$ is two-sided), hence, null-homologous in $P^3(\R)$.
Therefore, $B(\R)$ intersects the plane $X_{\xi\xi}$ at an
even number of points.
But it is these points that are
projected into the cusps of the corresponding oval, see Corollary \ref{generic contours}.
\endproof

We denote by $2\n_i$ the number of the imaginary cusps of $A$, and
by $2\n_r$ the number of real ones, so that $\n_i,\n_r\ge 0$, and
$\n_i+\n_r=3$.

Every oval, $O$, of a Zariski sextic $A$ is obviously
null-homologous in $P^2(\R)$ and, thus, divides the
latter into {\it the interior of $O$} (homeomorphic to a
disc), and {\it the exterior of $O$} (homeomorphic to a M\"obius band).
Another oval lying in the interior (respectively, exterior) of $O$
is called its {\it internal} (respectively, {\it external}) oval.
If $O$ has no internal ovals, then it is called {\it empty oval}.
An oval $O$ is called {\it ambient one} with respect to its internal ovals.

A cusp on $O$ may be
directed towards the interior, or the exterior, of $O$ and we call it
an {\it inward cusp}, or {\it outward cusp}, respectively.
The following statement follows directly from the definition of $\P_\pm$ and $\Cal A_\pm$.

\lemma\label{cusp-direction}
Any cusp of an oval $O\subset A(\R)$ is inward, if the region $\P_+$ lies
inside the interior of $O$, and outward, if it lies inside the exterior.
 All the real cusps are directed from $\Cal A_-$ to $\Cal A_+$.
\qed\endlemma

\corollary All the cusps on an oval of a real
Zariski curve are alike: either all inward,
or all outward. \qed
\endcorollary

\corollary\label{cusps-direction} For any real
Zariski sextic the following properties hold:
\roster\item there
can not be more than one oval with
inward cusps;
 \item
all the ovals in the exterior of an oval with inward cusps
 are smooth;
\item all the ovals in the interior of an oval with outward cusps are smooth;
 \item
if one of external ovals has an outward $($inward$)$ cusp, or one of
internal ovals has an inward $($resp. outward$)$ cusp, then $o=-$
$($resp. $o=+)$. \qed\endroster
\endcorollary

\lemma\label{polycuspidal-ovals}
If some oval of a real Zariski sextic is non-empty, then no other its oval may contain
{\eightrm(a)} an inward cusp, {\eightrm(b)} more than two outward cusps.
\endlemma

\proof
(a) If an inward cusp is on an external oval, then we take a line passing through this cusp and a point inside
an internal oval. If an inward cusp is on an internal oval, then it contains another cusp and we take a line through both of them. In each of the cases there will be a
contradiction to the Bezout theorem.

(b) If there is more than two outward cusps on an external oval, $\om$, then we take a line passing through an internal oval and one of the four cusps on $\om$ chosen so that this line intersects $\om$ at
some other point (it is possible, since it contains $>2$ cusps).
 If more then two, and thus, at least four outward cusps lie on an internal oval, $\om$, then we can find a line passing through two of the cusps and intersecting $\om$ in at least one more point. This will also contradict to
the Bezout theorem.
\endproof

\lemma\label{few-smooth-ovals}
There cannot be more than one smooth empty oval bounding a disc contained in $\Cal A_-$. In particular,
there cannot be more than one smooth empty oval
in the interior of an oval with inward cusps, as well as
in the exterior of a non-empty oval with outward cusps.
\endlemma

\proof
The projection $p_\R$ is 3-fold over $\Cal A_-$, so,
smoothness of the ovals implies that $X(\R)$ must have spherical components projecting to the discs in $\Cal A_-$
bounded by these ovals.
 On the other hand, $X(\R)$ has at most one spherical components.
\endproof

\subsection{The code of a Zariski sextic}
\label{codes}
The ambient differential-topological type of $A(\R)$ in
 $P^2(\R)$, $A$ being a Zariski sextic, is characterized by a certain code of $A$,
which is defined as follows.

An oval with $2\le2k\le6$ cusps has code $1_k$ if the cusps are outward,
and $1_{-k}$ if the cusps are inward; a smooth oval has code $1$.
The codes $n$ and $n_1$ are abbreviations for $1\+\dots\+1$ and
$1_1\+\dots\+1_1$ that
denote a group of $n$ empty
ovals, which are all smooth, or respectively, all have
 $2$ outward cusps (recall that according to Corollary \ref{cusps-direction}
 groups of $n\ge 2$ ovals with  inward cusps are impossible).
  For an arrangement in which an oval $O$ contains inside a set of
 ovals with the code $S$, we use the code $1_k\la S\ra$,  $1_{-k}\la S\ra$, or $1\la
 S\ra$ (depending on the number and the direction of cusps on $O$).
A few examples are shown on Figure \ref{codes-examples}.
 \midinsert\hskip10mm\epsfbox{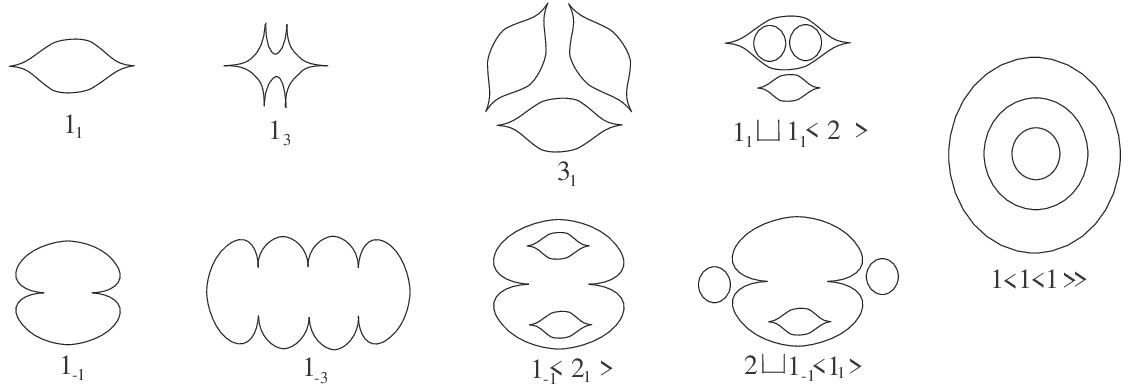} \figure
\label{codes-examples}
\endfigure
\endinsert
We may also ignore the cusps and describe the purely topological
(i.e., class $C^0$) arrangements of ovals by dropping the
subscripts from the codes.
 The result is called {\it the simple code of $A$}:
 it looks like $\a$ if all the ovals are empty, like
 $\a\+1\la\b\ra$ if
 one of the ovals contains $\b$ ovals inside and $\a$ ovals
 outside, or like
 $1\la1\la1\ra\ra$, if there is a nest of three ovals.
 It is convenient to allow sometimes an alternative form of this code,
and write $\a\+1\la0\ra$ for an arrangement of
 $\a+1$ empty ovals.
 For the empty set of ovals ($A(\R)=\oo$), we use the code $0$
 and call it the
 {\it null-code}. The code $1\la1\la1\ra\ra$ is called
 the
 {\it 3-nest code}.

 \lemma\label{real-schemes}
 Assume that a real sextic $A$ has six cusps and no other singular
 points. Then the arrangement of its ovals has one of the
 following simple codes:
 \roster\item
 $\a\+1\la\b\ra$, $0\le\a,\b\le4$, $\a+\b\le4$,
 \item $0$,
\item $1\la1\la1\ra\ra$.
\endroster
 \endlemma

\proof It follows from the Bezout theorem and Harnack's estimate of
the number of ovals (see Subsection \ref{types-of-curves-section}), like in the
well-studied case of non-singular sextics, see, for example, Gudkov's survey \cite{Gudkov}.
\endproof

\lemma\label{nested-no-cusps}
A real sextic cannot have real singular points, if its
real locus is a union of three nested ovals.
\endlemma

\proof
A line through a point inside the internal oval of the nest intersects such a sextic in 6 real points, which should be all non-singular not to contradict to the Bezout theorem.
\endproof

\subsection{The types I and II of (M-d)-curves}\label{types-of-curves-section}
Recall that the {\it Harnack inequality} $\ell(A)\le g(A)+1$, which bounds the
number of connected components, $\ell(A)$, of $A(\R)$ for a
non-singular real curve $A$ in terms of its genus, $g(A)$, extends
to singular irreducible curves as soon as one understands by
$\ell(A)$ and $g(A)$ the number of components and the genus after
normalization.
 Real curves with  $\ell(A)=g(A)+1$
 are  called {\it M-curves}. Let $d(A)=g(A)+1-\ell(A)$.
 If $d=d(A)>0$, then a curve is called an {\it (M-d)-curve}.

 An irreducible curve $A$ is said to be {\it of type} I\  if
$A(\C)\sm A(\R)$ is disconnected (i.e., has two components), and {\it
of type} II\
otherwise.
The following is well-known.

\lemma\label{type-lemma} Assume that $A$ is a real irreducible curve.
\roster\item If $A$ is an M-curve, then it is of type I.
 \item If $A$ is an (M-d)-curve, where $d$ is odd, then it is of type II.
 \endroster
If $A$ is a real sextic whose real locus is a nest of three ovals, then $A$ is of type I.
\qed\endlemma

\corollary\label{type} A real Zariski sextic has at most 5
ovals. It is of type I if it has 5 ovals, and of type II if it has
4,2, or no ovals. \qed\endcorollary

\subsection{Reversion of Zariski sextics}\label{reversion}
Given an arrangement $D\subset\Rp2$ of ovals (smooth, or cuspidal), a point $p\in\Rp2\smallsetminus D$,
and a simple closed one-sided curve $L\subset\Rp2\smallsetminus(D\cup\{p\})$,
one can obtain another arrangement of ovals,
$D'$ on the projective plane $\Rp2_{p,L}$
obtained from the given one by contracting $L$ and blowing up $p$.
 If we pass to the isotopy classes of $D$ and $D'$, one can identify  $\Rp2_{p,L}$ with the initial $\Rp2$
 and observe that the result of this operation depends only on the choice of
the component of $\Rp2\smallsetminus D$ containing $p$.
 An alternative description of this operation is
 to pick up an annulus neighborhood of $D$,
 $D\subset N\subset\Rp2$, $N\cong S^1\times[-1,1]$, and let $D'\subset\Rp2$ be the image of $D$ under the reversion
mapping $N\to N$, $(x,t)\mapsto(x,-t)$.
 We say that $D$ and $D'$ (considered up to isotopy)
 are {\it in reverse position} with respect to annulus $N$, or with respect to point $p$, and
 call this operation {\it reversion of $D$}.
 An oval containing the
 point $p$ inside, or equivalently, homologically non-trivial in $N$, will be called
{\it a principal oval of a reversion}.

 We say that $D\subset N\cong S^1\times[-1,1]$ is {\it trigonal with respect to $N$}, if each segment
$x\times [-1,1]\subset N$ intersects transversely $D$ at one or three points, except finitely many critical values
of $x$ for which segment $x\times [-1,1]$ is tangent to $D$ or contains a cusp whose tangent is transverse
to this segment. The following observations are trivial.

\lemma\label{trigonal-properties} If an arrangement of ovals $D$ is trigonal with respect to an annulus $N$, then
it has either one or three principal ovals. Moreover
 \roster\item
non-principal ovals are empty;
 \item
a non-principal oval cannot have inward cusps,
and may have maximum two outward cusps;
 \item
the ovals lying inside the principal oval, $\om$, of $D$ after reversion will lie outside the image of $\om$, and vice versa, the ones lying outside $\om$ will lie inside its image after reversion;
 \item
if $D$ has three principal ovals, then there is no other ovals, there is no cusps on the ovals, and
 $D$ as well as its reversion $D'$ form a 3-nest arrangement $1\la1\la1\ra\ra$.
\qed\endroster
\endlemma

We say that non-empty real Zariski sextics $A$ and $A'$ are {\it reversion partners} if
(a) $A(\R)$ is trigonal with respect to
some annulus $N\supset A(\R)$ and $A'(\R)$ is obtained by reversion of $A(\R)$,
(b) $A$ and $A'$ are both of the same type, and
(c) the signs $o(A)$ and $o(A')$ are opposite.

Note furthermore that the definition of reversion implies that
the cusps on the ambient ovals change their shape after reversion: the inward cusps
become outward and vice versa. And by contrary, the cusps on the non-principal ovals do not change their shape.

\corollary\label{reversion-cor} Assume that  real Zariski sextics $A$ and $A'$ are reversion partners. Then:
 \roster\item
if $A$ has 3-nest code $1\la1\la1\ra\ra$, then all three ovals are principal and $A'$ has the same code;
 \item
if an oval of $A$ has an inward cusp, or more than two outward cusps, then it is principal;
 \item
if the simple code of $A$ is $\a\+1\la\b\ra$, $\b>0$,
then the ambient oval must be principal and $A'$ must have simple code $\b\+1\la\a\ra$;
 \item
 if the complete code of $A$ is $\a_k\+1_m\la\b_n\ra$,
 $\b>0$, $k,n,m\in\Z$ (zero index here means that it should be dropped),
then the complete code of $A'$ is $\b_n\+1_{-m}\la\a_k\ra$.
\qed\endroster\endcorollary

Rule (4) in Corollary \ref{reversion-cor} can be extended to $\beta=0$
in two cases: if by using
\ref{reversion-cor}(2) we can determine which of the ovals is principal, or if all the ovals look alike.
For instance, if the complete code of $A$ is $1\+1_n$, where $n<0$, or $n>1$, then
 $A'$ has complete code $1_{-n}\la1\ra$. If $A$ has code $n_1$, then $A'$ has code $1_{-1}\la(n-1)_1\ra$.

\remark{Remark}
 If $A$ has, for instance, code $1\+1_1$ then the rules stated above do not give an answer, which
of the ovals is principal, and so, one could question if the complete code of a reversion partner
of $A$ is determined uniquely by the complete code of $A$? Is it possible for $A$ to have several reversion partners?
 In what follows we will prove that a real Zariski sextic $A$ cannot have more than
one reversion partner. For example, the partner for $A$ with the code $1\+1_1$ has code $1\la1_1\ra$, and there is no
Zariski sextic with the code $1_{-1}\la1\ra$.
\boxedR\endremark

\subsection{The relation between the topology of
cubic surfaces and their Zariski sextics}
Here, we consider a real nonsingular cubic surface $X$ and one of its real Zariski sextics, $A$.
Recall that the real locus $X(\R)$ of
$X$ may consist of
two components, one homeomorphic to $S^2$
and another to $P^2(\R)$. Otherwise, $X(\R)$ is homeomorphic to
$P^2(\R)$ with $\h\le3$ handles. In particular, $X(\R)$ is
determined, up to homeomorphism, by its Euler characteristic $\chi(X(\R))=1-2\h$, or
equivalently, by $\h\in\{-1,0,1,2,3\}$, where $\h=-1$ corresponds to
the case of disconnected $X(\R)$.

\lemma\label{cubic-via-code} Assume that $A$ has code $\a\+1\la\b\ra$, so that
$d(A)=4-\a-\b$. Then
$$\chi(X(\R))=\cases &3+2(\a-\b)-2\n_r=4\a+2d(A)+2\n_i-11
\text{ \ \ if $o(A)=-$,}
\\
&1+2(\b-\a)-2\n_r=4\b+2d(A)+2\n_i-13 \text{ \ \ if $o(A)=+$.}
\endcases$$
In the case of null-code, $\chi(X(\R))$ equals $1$ if $o(A)=-$ and $3$ if $o(A)=+$.

\noindent
In the case of 3-nest code,
$\chi(X(\R))$ equals $3$ if $o(A)=-$ and $1$ if $o(A)=+$.
\endlemma

\proof
As it follows from Corollary \ref{generic contours},
$\chi(X(\R))=3\chi(\RR_-)+\chi(\RR_+)-2\n_r=1+2\chi(\RR_-)-2\n_r$.
\endproof

\corollary\label{h} If $A$ has code $\a\+1\la\b\ra$, then
$$\h(X)=\cases
\n_r+\b-\a-1&=6-(2\a+d(A)+\n_i) \text{ \ \ if $o(A)=-$,}\\
\n_r+\a-\b&=7-(2\b+d(A)+\n_i) \text{ \ \ if $o(A)=+$.}\endcases$$
If $A$ has null-code, then $\h(X)=0$ if $o(A)=-$, and $\h(X)=-1$ if $o(A)=+$.

\noindent If $A$ has 3-nest code, then $\h(X)=-1$ if $o(A)=-$, and $\h(X)=0$ if $o(A)=+$.
\qed\endcorollary

The following observation is a kind of refinement of Lemma \ref{few-smooth-ovals}.

\lemma\label{one-smooth-oval}
If the disc bounded by a smooth empty oval of $A(\R)$ lies in $\Cal A_-$, then $\h(X)=-1$.
In particular, for any $A$ with code $\a\+1\la\b\ra$, $\b\ge0$,
\roster\item
if either some external oval or the ambient oval has an outward cusp and
one of the external ovals is smooth, then $\nu_r=\a-\b$;
\item
if either an internal oval has an outward cusp,
or the ambient oval has an inward cusp and one of
the internal ovals is smooth, then $\nu_r=\b-\a-1$.
\endroster
\endlemma

\proof
Smoothness of an oval bounding a disc in $\Cal A_-$,
over which the projection $X(\R)\to P(\R)$ is three-to-one,
implies that $X(\R)$ must contain a spherical component projecting to this disc, and thus, $\h(X)=-1$.
Assumption (1), as well as (2), guarantees that the smooth oval bounds a disc component of $\Cal A_-$,
and the conclusion $\h(X)=-1$ is reformulated via Corollary \ref{h}.
\endproof

\corollary\label{empty-ovals-cor}
If all the ovals of $A(\R)$ are empty and $o(A)=-$,
then either (a) all ovals have outward cusps, or (b) one oval is smooth and each of the others has precisely two outward cusps.
\endcorollary

\proof
Since $o(A)=-$, the cusps are outward.
By Lemma \ref{few-smooth-ovals} there cannot be more than one smooth oval, and letting $\b=0$ in Lemma \ref{one-smooth-oval}(1), we see that in the presence of a smooth oval the others cannot have more than one
pair of cusps.
\endproof

\subsection{Relation to the double covering cuspidal K3-surface $Y$}\label{cuspidal-K3}
By taking double covering $\p_Y\:Y\to P^2$ ramified along Zariski sextic
$A$, we obtain a K3-surface
$Y$ which has six cuspidal singular points. There exist two
liftings of the complex conjugation in $P^2$ to an involution in
$Y$; they differ by the deck transformation of the covering and
both are anti-holomorphic.
Let us choose and denote by $\conj_Y\:Y\to Y$ the
one whose real locus, $Y(\R)=\Fix(\conj_Y)$, is projected by $\p_Y$
to $\RR_n$.  We call such $\conj_Y$ the {\it M\"obius involution}.

The group $H_2(Y)/\operatorname{Tors}$ is a
free abelian
group of rank $b_2(Y)=22-12=10$.
We consider the involution $(\conj_Y)_*$ induced in $H_2(Y)/\operatorname{Tors}$
and denote by $r_\pm$
the ranks of the $\pm1$-eigengroups, $\{x\in H_2(Y)/\operatorname{Tors}\,|\,(\conj_Y)_*(x)=\pm x\}$.

\lemma\label{codes-and-ranks}
For any real Zariski sextic $A$ we have $\chi(\Cal A_-)=1+\frac{r_+-r_-}2$.
Specifically, if $A$ has
 code $\a\+1\la\b\ra$, $\a,\b\ge0$, so that
$d(A)=4-(\a+\b)$, then
 $$\aligned
r_+=&(\b-\a)+4=2\b+d(A), \\ 
r_-=&(\a-\b)+6=2\a+d(A)+2. 
\endaligned$$
 \newline
If $A$ has null code,
then $r_+=r_-=d(A)=5$.
\newline
If $A$ has 3-nest code, then
$r_+=d(A)+2=4$ and $r_-=6$.
\endlemma

\proof We have obviously $r_++r_-=b_2(X)=10$, and the
Lefschetz fixed-point formula applied to $\conj_Y$ yields
$r_+-r_-=\chi(Y(\R))-2=2(\chi(\RR_-)-1)$, where
$\chi(\RR_-)=\b-\a$ in the case
of code $\a\+1\la\b\ra$. For the null-code and 3-nest code
one has $\chi(\RR_-)=\chi(\Rp2)=1$ and $\chi(\RR_-)=0$ respectively.
\endproof

\corollary\label{empty-to-latticecoding}
For any real Zariski sextic
the following three conditions are equivalent:
\roster\item
the sextic has null-code;
\item $r_+=d(A)=5$;
\item $d(A)>r_--2$. \qed
\endroster
\endcorollary

Next, we obtain the following relation between
$Y(\R)$ and $X(\R)$, or more precisely, between $r_\pm$ and  $\h=\frac{1-\chi(X(\R))}2$.

\corollary\label{codes-and-cubic}
For any real Zariski sextic $A$ with
code $\a\+1\la\b\ra$
the following identities hold:
$$\aligned \text{ if $o(A)=-$ then } &{\cases
r_+=&\n_i+2+\h,\\ 
r_-=&(4-\n_i)+2+(2-\h); 
\endcases}\\
\text{ if $o(A)=+$ then }&{\cases
r_-=&\n_i+2+(\h+1),\\ 
r_+=&(4-\n_i)+2+(1-\h). 
\endcases}
\endaligned$$
\endcorollary

\proof
 It follows from Lemma \ref{codes-and-ranks} and Corollary \ref{h}.
\endproof

\subsection{Deformation classification statement}\label{def-class}
By the {\it ID of a real Zariski sextic $A$} we mean the triple
$\ID(A)=(\text{complete code of $A$},{\text{type of $A$}},o(A))$, i.e.,
its complete code enhanced with two additional bits of information:
the type (I or II), and the sign (+ or -).

Let us recall that two real Zariski sextics are said {\it to be equivalent up to equivariant deformation}, or shortly {\it belonging to the same deformation class}
if and only if they belong to the same component of the space of real Zariski sextics.

\midinsert
 $$\matrix \\ \\ 1\\2\\3\\4\\5\\6\endmatrix
\matrix\text{Table 1A. Zariski sextics without a
partner.}\label{table-without-partner}
\\
 \boxed{\matrix
\text{simple codes}&\n_r(A)&o(A)&\text{complete codes}
&\text{types}\\
  1\la4\ra&0&-&1\la4\ra&I\\
  1\la3\ra&1&-&1_1\la3\ra&II\\
  1\la2\ra&1&+&1\la1_1\+1\ra&I\\
  1\la2\ra&2&-&1_2\la2\ra&II\\
  1\la1\ra&3&-&1_3\la1\ra&II\\
  1\la1\ra&0&+&1\la1\ra&II\\
\endmatrix}\endmatrix$$
\endinsert

\midinsert\topcaption{Zariski sextics with a partner}\endcaption
 $\matrix \\ \\ \\
1\\2\\3\\ \\4\\5\\6\\7\\8\\9\\ \\ 10\\11\\12\\13\\14\\15\\16\\17\\18\\19\\20\\ \\ 21\\22\\23\\24\\25\\26\\ \\27\\28\\29\\30 \\ \\31\\
\endmatrix
\matrix\text{Table 1B. The case of
$o(A)=-$.}\label{partner-minus}
\\
\boxed{\matrix
\text{simple}&\n_r(A)&\text{complete}&\text{types}\\
\text{codes}&&\text{codes}\\
 3\+1\la1\ra&3&3_1\+1\la1\ra&I\\
 2\+1\la2\ra&2&2_1\+1\la2\ra&I\\
 1\+1\la3\ra&1&1_1\+1\la3\ra&I\\
 \\
 4&3&3_1\+1&II\\
 2\+1\la1\ra&3&2_1\+1_1\la1\ra&II\\
 2\+1\la1\ra&2&2_1\+1\la1\ra&II\\
 1\+1\la2\ra&2&1_1\+1_1\la2\ra&II\\
 1\+1\la2\ra&1&1_1\+1\la2\ra&II\\
 1\la3\ra&0&1\la3\ra&II\\
 \\
 3&3&3_1&I\\
 3&3&3_1&II\\
 3&2&2_1\+1&II\\
 1\+1\la1\ra&3&1_1\+1_2\la1\ra&II\\
 1\+1\la1\ra&2&1_1\+1_1\la1\ra&I\\
 1\+1\la1\ra&2&1_1\+1_1\la1\ra&II\\
 1\+1\la1\ra&1&1_1\+1\la1\ra&II\\
 1\la1\la1\ra\ra&0&1\la1\la1\ra\ra&I\\
 1\la2\ra&1&1_1\la2\ra&I\\
 1\la2\ra&1&1_1\la2\ra&II\\
 1\la2\ra&0&1\la2\ra&II\\
 \\
 2&3&1_2\+1_1&II\\
 2&2&2_1&II\\
 2&1&1_1\+1&II\\
 1\la1\ra&2&1_2\la1\ra&II\\
 1\la1\ra&1&1_1\la1\ra&II\\
 1\la1\ra&0&1\la1\ra&II\\
 \\
 1&3&1_3&II\\
 1&2&1_2&II\\
 1&1&1_1&II\\
 1&0&1&II\\
 \\
 0&0&0&II\\
\endmatrix}\endmatrix$
\hskip-2mm $\matrix\text{Table 1C. The case of
$o(A)=+$.}\label{partner-plus}
\\
\boxed{\matrix
\text{simple}&\n_r(A)&\text{complete}&\text{types}\\
\text{codes}&&\text{codes}\\
 1\+1\la3\ra&3&1\+1\la3_1\ra&I\\
 2\+1\la2\ra&2&2\+1\la2_1\ra&I\\
 3\+1\la1\ra&1&3\+1\la1_1\ra&I\\
 \\
 1\la3\ra&3&1\la3_1\ra&II\\
 1\+1\la2\ra&3&1\+1_{-1}\la2_1\ra&II\\
 1\+1\la2\ra&2&1\+1\la2_1\ra&II\\
 2\+1\la1\ra&2&2\+1_{-1}\la1_1\ra&II\\
 2\+1\la1\ra&1&2\+1\la1_1\ra&II\\
 4&0&4&II\\
 \\
 1\la2\ra&3&1_{-1}\la2_1\ra&I\\
 1\la2\ra&3&1_{-1}\la2_1\ra&II\\
 1\la2\ra&2&1\la2_1\ra&II\\
 1\+1\la1\ra&3&1\+1_{-2}\la1_1\ra&II\\
 1\+1\la1\ra&2&1\+1_{-1}\la1_1\ra&I\\
 1\+1\la1\ra&2&1\+1_{-1}\la1_1\ra&II\\
 1\+1\la1\ra&1&1\+1\la1_1\ra&II\\
 1\la1\la1\ra\ra&0&1\la1\la1\ra\ra&I\\
 3&1&1_{-1}\+2&I\\
 3&1&1_{-1}\+2&II\\
 3&0&3&II\\
 \\
 1\la1\ra&3&1_{-2}\la1_1\ra&II\\
 1\la1\ra&2&1_{-1}\la1_1\ra&II\\
 1\la1\ra&1&1\la1_1\ra&II\\
 2&2&1_{-2}\+1&II\\
 2&1&1_{-1}\+1&II\\
 2&0&2&II\\
 \\
 1&3&1_{-3}&II\\
 1&2&1_{-2}&II\\
 1&1&1_{-1}&II\\
 1&0&1&II\\
 \\
 0&0&0&II\\
\endmatrix}\endmatrix$
\vskip2mm
\endinsert

\theorem\label{main-sextics}
Each real Zariski sextic is determined up to equisingular
deformation by its $\ID$. The list of $\ID s$ of real Zariski
sextics is given in Tables 1A-C.
\endtheorem

The proof of Theorem \ref{main-sextics} is one of our main goals. It
is given in the very end of the paper, see Section \ref{proof_main}.

\section{Preliminaries on lattice theory}\label{Chapter_Lattice_preliminaries}

\subsection{Even lattices and their discriminants}
 By a {\it lattice} we mean a free abelian
 group of finite rank endowed with a {\sl non-degenerate} bilinear
 symmetric $\Z$-valued pairing, called also the {\it inner product of a
 lattice}.

We denote by $\la n\ra$, $n\in\Z\sm\{0\}$,
the  lattice of rank $1$ whose generator has square $n$, by $\U$ the
lattice of rank $2$ defined by matrix
 $\left[\matrix 0&1\\1&0\endmatrix\right]$,
 and by $\A_n$ ($n\ge1$), $\DD_n$ ($n\ge4$), $\E_n$ ($6\le n\le8$)
 the lattices generated by the corresponding
 {\sl negative definite} root systems (the same notation, $\A_n$,  $\DD_n$, and $\E_n$,
 is used also to refer to the corresponding types of simple singularities).
 Given two lattices $L$ and $M$, we
 denote by $L + M$ their direct sum (this notation does not lead to confusions, since
non-direct sums of lattices are never considered).
 Notation $nL$, $n\ge1$,
 stands for the direct sum of
 $n$ copies of $L$, while $L(n)$, $n\in\Z\sm\{0\}$, denotes the
 result of rescaling the lattice $L$ with a scale factor $n$, i.e., $L(n)$ as a group coincides with $L$
 but the product of elements in $L(n)$ is $n$ times greater than in $L$.
 An isomorphism of lattices is indicated by writing $L=M$.

In this paper,
we deal generally with {\it even lattices}, unless it is stated
otherwise (although some of the techniques that we use or
develop can be adapted to the case of odd lattices as well).
Recall that a lattice $L$ is called even, if
$x^2$ is even for any $x\in L$.
When a lattice $L$ is equipped with a lattice involution $c\:L\to L$,
we introduce a {\it c-twisted inner product} $\la x,y\ra_c=x \cdot c(y)$,
which is obviously bilinear, symmetric and non-degenerate,
and say that the pair $(L,c)$ is of {\it type I} if the c-twisted product is even,
and of {\it type II} otherwise.

Given a lattice, $L$, we denote by $\discr L$ its (finite abelian) discriminant group, $L^*/L$, where $L^*=Hom(L,\Z)$ is identified with a subgroup in $L\otimes\Q$ by means of the lattice pairing.
 The group $\discr L$ is endowed with the nondegenerate $\Q/\Z$-valued
 inner product $[x][y]=xy\mod\Z$, where $x,y\in L^*$
and $[x],[y]$ stand for the cosets. Our assumption that
the lattice $L$ is even endows $\discr L$ with
a $\Q/2\Z$-valued quadratic  refinement
of the inner product.
This refinement, $\q_L\:\discr L\to\Q/2\Z$,
is given
by $\q_L([x])=x^2\mod2\Z$; it is called the {\it discriminant form} of $L$.
 The relation $\q_L([x]+[y])=\q_L([x])+\q_L([y])+2[x][y]$
implies that $\q_L$ determines the inner product of $\discr L$.

The pair $(\discr L, \q_L)$ is called the {\it discriminant} of
$L$, and when it does not lead to a confusion is denoted simply by
$\discr L$ or $\q_L$.

\subsection{Groups with inner products and quadratic refinements}\label{enhanced-groups}
A finite abelian group $G$ endowed with a non-degenerate symmetric
bilinear form $G\times G\to\Q/\Z$ will be called a {\it finite
inner product group}.
We denote by $ab$ the inner product of elements $a,b\in G$
 and use notation $\la\frac{m}n\ra$
(with coprime $m$ and $n$) for a finite inner product
cyclic group $\Z/n$ that has
$a^2=\frac{m}n\in\Q/\Z$ for one of generators $a\in \Z/n$.

If a finite inner product group
$G$ is endowed additionally with a {\it quadratic
refinement}, $\q\:G\to\Q/2\Z$, then $G$ will be called
an {\it enhanced group}.
By definition, $\q$ must be related to the inner product as
follows: for all $a,b\in G$, $n\in\Z$,
 \roster\item
 $\q(a+b)=\q(a)+\q(b)+2ab$,
 \item
 $\q(na)=n^2\q(a)$.
 \endroster

These relations imply
$\q(-a)=\q(a)$, $\q(a)=a^2\mod\Z$,
 and  $\q(2a)=4a^2\mod2\Z$ (note that $4a^2$ is well-defined modulo 4 and, thus, modulo 2).
 The latter relation shows that $\q(a)$ is defined
 uniquely  by the inner product as soon as $a$ is divisible by 2.
In the context of enhanced groups, notation $\la\frac{m}n\ra$ (with coprime $m$ and $n$) has a meaning that $\q(a)=\frac{m}n\in\Q/2\Z$ for a generator $a\in \Z/n$.
Note that in this case either $m$ or $n$ must be even, since $\q((n+1)a)=\q(a)$ if and only if $\frac{m}n((n+1)^2-1)$ is even.

\lemma\label{RIP=IP} Any finite inner product group G of odd order
has a canonical quadratic refinement and, thus, is an enhanced
group.
\endlemma

\proof
Let $\q(a)=4(\frac{a}2)^2\mod2\Z$, then all the required properties are satisfied.
\endproof

\example\label{2/3-example}
A finite inner product group $\la\frac13\ra$ being enhanced takes notation $\la-\frac23\ra$ according to our conventions.
An enhanced group $\la\frac16\ra$ splits into a direct (orthogonal) sum $\la-\frac12\ra+\la\frac23\ra$.
\endexample

In general, the set of quadratic refinements of a given inner product
in a group G form an affine space over $Hom(G,\Z/2\Z)$.

We say that a subgroup $K\subset G$
of a finite inner product group $G$ is {\it non-degenerate} if its
{\it kernel} $\{x\in K\,|\,xK=0\}$ is trivial. It is
straightforward to check the following.

\lemma\label{split-out} If $K$ is a non-degenerate subgroup of a
finite inner product group $G$, then $K$  splits out as a direct
summand, $G=K+K^\perp$. \qed\endlemma

\subsection{The p-components}
By a {\it p-group}, where $p$ is prime, we mean a finite abelian group $G$, such that the order, $\ord(x)$, of any
element $x\in G$ is a power of $p$. Note that any p-group $G$ can be presented as
a direct sum of cyclic subgroups of the form $\Z/p^k$, $k\ge1$, and that the number of summands is an invariant
of $G$ independent of a decomposition; it is called the
{\it rank of G}.
 If $\ord(x)=p$ for all non-identity $x\in G$, then
$G$ is called an {\it elementary p-group}; such a group
$G$ can be viewed as a vector space over $\Z/p$.

Any finite abelian group $G$ splits into a direct sum of its maximal p-subgroups $G_p$
called {\it prime components}, or {\it $p$-components} of $G$.
The $p$-primary component of $G$
is non-trivial if and only if $p$ divides the order of $G$ and coincides with
the subgroup formed in $G$ by elements whose order is a power of p.
 With respect to any inner product in $G$ the subgroups $G_p$ must be orthogonal,
 which reduces studying of inner products to the case of p-groups.

\lemma\label{components-splitting} Assume that $G$ is a finite
inner product group (for instance, an enhanced group).
\roster\item  If $x\in G$ has order $n$, then $x^2\in
\frac1n\Z/\Z$ (respectively, $\q(x)\in\frac1n\Z/2\Z$).
 \item
 Any two prime components, $G_{p_1}$ and $G_{p_2}$, $p_1\ne p_2$,
are orthogonal with respect to the inner product in $G$.
\endroster
\endlemma

\proof If $x,y\in\discr L$ and $x$ has order $n$, then
$xy\in\frac1n\Z/\Z$ and $x$ is orthogonal to $ny$, since $nxy= x(ny)=(nx)y=0\in\Q/\Z$. This implies
both (1) and (2), since $\q(x)=x^2\mod\Z$.
\endproof

\subsection{The discriminant p-ranks of lattices}
In a particular case of a lattice $L$ (which can be odd here) and $G=\discr L$ (viewed here only as a group),
we call the prime components $G_p$  the {\it discriminant p-components} and their ranks
the {\it discriminant p-ranks of $L$}.
We denote them $\discr_pL$ and  $r_p=r_p(L)$, respectively.
 The primes $p$ for which $\discr_pL\ne0$ (i.e., the
prime divisors of $|\discr L|$) are called
the {\it discriminant factors of $L$}.
 A discriminant factor $p$ is said to be {\it elementary} if $\discr_pL$ is an elementary p-group, we say also
that $L$ is {\it p-elementary} in this case.
A lattice $L$ is said to be {\it divisible by $p$} if $xy$ is divisible by $p$ for all $x,y\in L$, or equivalently,
if $L'=L(\frac1p)$ is a lattice.

Consider an endomorphism
$m_p\:\discr L\to\discr L$,
$m_p(x)=px$, then its kernel $K_p\subset \discr L$ is an elementary $p$-group of rank $r_p$.
Let $L^*_p=\psi^{-1}(K_p)\subset L^*$ be the
pull-back with respect to the projection $\psi\:L^*\to\discr
L=L^*/L$.

\proposition\label{p-divisible}
For any lattice $L$ and
prime $p$, the $p$-rank $r_p$ of $L$ is not greater than the rank $r$ of $L$. Moreover, $L$ is divisible by $p$ if and only if $r_p=r$.
In the latter case
the following properties hold:
\roster
\item there is a canonical exact sequence $0\to\discr L'@>f>> \discr L\to
L^*/pL^*\to 0$, where $L^*/pL^*=(\Z/p)^r$, and the restriction of $f$ yields
a group isomorphism $\discr_qL'= \discr_qL$
for any prime $q\ne p$, as well as
an isomorphism between $\discr_pL'$
and  the subgroup
 $p\,\discr_pL\subset\discr_pL$;
 \item
the group isomorphism $\discr L'\to f(\discr L')\subset\discr L$ identifies
the discriminant inner product in $\discr L$ restricted to $f(\discr L')$ with
the discriminant inner product in $\discr L'$ multiplied by $p$,
and if $L'$ is even, then $\q_L|_{f(\discr L')}$ is identified with $p\q_{L'}$.
\endroster
\endproposition

\proof
Note that $K_p=(\Z/p)^{r_p}\subset \discr L$ is the
maximal subgroup of exponent $p$ in $\discr_pL$ (and thus, in
$\discr L$).
As it follows from
definition, $pL^*_p\subset L$, and thus,
 $$K_p=L^*_p/L\subset L^*_p/pL^*_p,$$
where the latter group is isomorphic to $(\Z/p)^r$, since $L^*_p$
is a free abelian group
of rank $r$ (as it contains $L$). Thus, $r_p\le r$.
In the case of equality $r_p=r$, we have
$K_p=L^*_p/pL^*_p$, that is $L=pL^*_p$. So, for any $x,y\in L$,
$x(\frac1py)\in\Z$, since $\frac1py\in L^*_p\subset L^*$, and
thus, $xy$ is divisible by $p$, and so $L$ is divisibly by $p$.
Conversely, if $L$ is divisibly by $p$, then $\frac{x}p\in L_p^*$ for any $x\in L$, and so
$L=pL_p^*$ and $K_p=L_p^*/pL_p^*$ has rank $r_p=r$.

The exact sequence in (1), with its properties, follow from the
observation that under our identification of $L'$ with $L$ as a
group, the dual $(L')^*$ of $L'$ is identified with $pL^*$.
\endproof

\remark{Remark} Given a free abelian group $M$ of rank $m$, and
its subgroup $L\subset M$ of the same rank, one can find a basis
$e_1,\dots,e_m$ of $M$, such that $L$ is spanned by some multiples
 $k_1e_1,\dots,k_me_m$, $k_i\ge1$, and $k_i$ divides $k_{i+1}$ for
 $i=1,\dots,m-1$. This yields a direct sum decomposition
 $M/L=\Z/k_1+\dots+\Z/k_m$.
If we apply it to $L\subset L^*$, then we get another proof of Proposition \ref{p-divisible}.
\boxedR\endremark

\lemma\label{r2=r-mod2}
For every even lattice $L$, we have
$r_2(L)=r(L)\mod2$.
\endlemma

\proof
In accordance with the definition of the $2$-rank, $r_2$ is equal to the rank of the
subgroup of $\discr L$ that is formed by elements of order $2$.
In its turn, the rank of this subgroup is equal to
the rank of $\discr(L)/2\discr(L)$, and thus to the rank
of the radical of the $\operatorname{mod} 2$ reduction of $L$ (as a quadratic
space). This reduction is even, since $L$ is even. Therefore,
 $r(L)-r_2(L)$ as the rank of a non degenerate even $\Z/2$-valued form is even.
\endproof

\subsection{The Brown invariant}
As is known, for any enhanced group $(G,\q)$, the Gaussian sum
$\Cal G(\q)=\sum_{x\in G}e^{\pi i\q(x) }$ has absolute value
$\sqrt{|G|}$, whereas
$\frac{\Cal G(\q)}{\sqrt{|G|}}$
belongs to the group $\mu_8$ of eight's roots of $1$. In modern
terminology one speaks on the (generalized) {\it Brown invariant},
$\Br(\q)\in\Z/8\cong \mu_8$. This Gaussian sum appears in Van der
Blij's famous formula that relates, in the case of even lattices
and their discriminants, the argument of $\Cal G(\q)$ with
the signature of the lattice.

\theorem\label{van-der-blij} (Van der Blij \cite{vdB})
If $L$ is an even lattice, then $\Br(\q_L)$
is equal to the mod 8 residue of the signature $\s(L)$ of $L$.\qed
\endtheorem

 Note that the definition via Gauss sums immediately implies additivity,
 $$\Br(\q_1\oplus \q_2)=\Br(\q_1)+\Br(\q_2).$$

Speaking on even lattices $L$, we let $\Br_p(L)=\Br(\q_L|_{\discr_p(L)})$ and use
 an alternative notation $\Br(L)$ for
$\Br(\q_L)=\s(L)\mod8$ (since
in certain formulas
it is more instructive and convenient
than $\s(L)\mod8$).
 Then Lemma \ref{components-splitting} and the above additivity formula imply
the following statement.

\proposition\label{Br_additivity}
For any even lattice $L$, we have $\Br(L)=\sum_{\text{prime } p}\Br_p(L)$.
\qed\endproposition

\subsection{Elementary enhanced 2-groups}
Assume that $G$ is an enhanced 2-group
with the quadratic refinement
$\q$. The finite inner product of $G$ is said to be {\it even} if $x^2=0$ for all $x\in G$, or in terms of $\q$, if
the values of $\q(x)$ are integer (and thus, $\q(x)\in\Z/2\Z$).
Otherwise, the inner product is
called {\it odd}. We encode this parity by putting $\delta_2(\q)=0$ in the even case, and $\delta_2(\q)=1$ otherwise.
For each 2-elementary lattice $L$, we define its {\it discriminant parity}
$\delta_2(L)\in\Z/2$ by letting
$\delta_2(L)=\delta_2(\q_L)$.

\example\label{2-rank-1} The group $\Z/2=\{[0],[1]\}$ has a unique
inner product, $[1][1]=\frac12$, $[*][0]=[0][*]=0$, but two
possible quadratic refinements, $\q([1])=\pm\frac12$, $\q([0])=0$,
denoted in accordance with our convention in Section \ref{enhanced-groups} by $\la\pm\frac12\ra$.
 Note that an even enhanced 2-group cannot contain $\la\pm\frac12\ra$ as a direct summand.

 Observing that $\discr(\la\pm2\ra)=\la\pm\frac12\ra$ we
obtain $\delta_2(\la\pm2\ra)=1$ and $\Br(\la\pm2\ra)=\pm1$.
\endexample

 \example\label{2-rank-2}
The discriminant groups of lattices $\U(2)$ and $\DD_4$ are both
isomorphic to $G=\Z/2+\Z/2$, with the inner product
$\left[\matrix 0&\frac12\\
\frac12&0\endmatrix\right]$. The quadratic refinements are
 $\q_{\U(2)}(1,0)=\q_{\U(2)}(0,1)=0$,
$\q_{\U(2)}(1,1)=1$, and $\q_{\DD_4}(1,1)=\q_{\DD_4}(1,0)=\q_{\DD_4}(0,1)=1$.
We denote the corresponding enhanced groups by $\uu$ and $\vv$
respectively. Here we have $\delta_2(\U(2))=\delta_2(\DD_4)=0$,
whereas $\Br(\U(2))=0$ but $\Br(\DD_4)=4$.
\endexample

The following  statement is a straightforward consequence of well
known "uniqueness"
results, see \cite{W1} (cf., \cite{GM}).

\theorem\label{2-classification}
Any elementary enhanced 2-group $(G,\q)$ is characterized up to
isomorphism by its rank, the Brown invariant, and the parity. The
only non-split enhanced 2-groups are $\la\pm\frac12\ra$,  $\uu$,
and $\vv$. Moreover:
 \roster\item If $\delta_2(\q)=0$, then
 $G=a\uu+b\vv$, for some $a,b\ge0$.
 In particular, the rank of $G$ is even and
 $\Br(\q)=a\Br(\uu)+b\Br(\vv)= 4b\mod8$
 is divisible by $4$.
 \item  $a\uu+b\vv$ is isomorphic to  $a'\uu+b'\vv$
if and only if $a+b=a'+b'$ and $a=a'\mod2$.
 \item
If $\delta_2(\q)=1$, then $G=a\la\frac12\ra+b\la-\frac12\ra$, for some $a,b\ge0$, $a+b=r_2$.
In particular, $\Br(\q)=a-b\mod8$.
 \item $a\la\frac12\ra+b\la-\frac12\ra$ is
isomorphic to $a'\la\frac12\ra+b'\la-\frac12\ra$ if and only if $a+b=a'+b'$
and $a=a'\mod 4$. \qed
\endroster
\endtheorem

If $G$ is a finite inner product elementary 2-group, then the map
$x\mapsto x^2$ is a homomorphism $G\to\frac12\Z/\Z=\Z/2$ and,
therefore, there exists one and only one element $v\in G$ such that $vx=x^2$ for all
$x\in G$; this $v$ is called the {\it characteristic element of $G$}.

\lemma\label{Wu-even} Assume that $(G,\q)$ is an elementary
enhanced 2-group and $v\in G$ such that $\q(v)\ne0$. Then:
\roster\item The orthogonal complement $v^\perp$ is an enhanced
subgroup of $G$. \item
 $v^\perp$ is even
if and only if $v$ is the characteristic element
of $G$.
\item
 If $v$ is characteristic, then
$2\q(v)\in\Z/4\Z$ equals to the $\roman{mod}\,4$-residue of
$\Br(\q)\in\Z/8$.
\item
If $\q(v)=\pm\frac12$, then the Brown invariant $\Br\la v\ra$ of the subgroup
$\la v\ra\subset G$ spanned by $v$ is respectively $\pm1$.
\endroster\endlemma

\proof Item (4) is already seen in Example \ref{2-rank-1}.
Item (1) follows from non-degeneracy of $\q|_{v^\perp}$, which is due to $\q(v)\ne0$.
Item (2) follows from $\q(x)=xv\mod\Z$, by the definition of the characteristic element $v$.
To prove item (3) it is sufficient,
due to additivity of $\q(v)$ and $\Br$ with respect to the direct sums
and the classification given in Theorem \ref{2-classification},
to check the required relation
on the
four non-split enhanced groups
$\la\frac12\ra$, $\la-\frac12\ra$, $\uu$, and $\vv$
(cf., Examples \ref{2-rank-1} and \ref{2-rank-2}).
\endproof

\subsection{Elementary 3-groups with an inner product}
The group $G=\Z/3=\{[0],[1],[2]\}$
admits two different inner product structures, $\la\frac13\ra$, and  $\la-\frac13\ra$.
For the first one, $[1][1]=[2][2]=\frac13$, and for the other, $[1][1]=[2][2]=-\frac13$.
 The quadratic refinements take values $\q([1])=\q([2])=-\frac23 $ and
$\q([1])=\q([2])=\frac23 $,
respectively (cf. Example \ref{2/3-example});
that is why to indicate that we deal with these enhanced structures we use notation $\din$ and $\dip$.

\example
The discriminant of the lattice $\A_2$ is
$\din$, while $\discr \A_2(-1)$ and $\discr \E_6$
 are both isomorphic to $\dip$.
 It is easy to check also that $\discr_3\la\pm6\ra=\discr_3\A_2(\pm2)=\di\pm23$.
\endexample

The  statements of the following Lemma are well known and can be
easily extracted, for example, from \cite{W1}.

 \lemma\label{3-classification}
  \roster\item
Any finite inner product 3-group is isomorphic to $a\dip+b\din$,
for some $a,b\ge0$, $a+b=r_3$.
 \item
 $\Br(a\dip+b\din)=2(a-b)\mod8$.
 \item
 $a\dip+b\din$ is isomorphic to
 $a'\dip+b'\din$
 if and only if
$a+b=a'+b'$ and $a=a'\mod2$. In the other words, such
inner product groups are isomorphic if and only if
their ranks and Brown invariants coincide. \qed
 \endroster
 \endlemma

\corollary Finite inner product 3-groups $(G,\q)$ are characterized up
to isomorphism by pairs $(a,b)$, where $a\in\{0,1\}$ and $b\ge0$,
such that $G=a\dip+b\din$. \qed
\endcorollary

Theorem \ref{van-der-blij} and Proposition \ref{Br_additivity} with Lemma \ref{3-classification} imply the following.

\corollary\label{vanderblij-refined} Assume that a lattice $L$ is
even and has only discriminant factors 2 and 3, which are both
elementary. Then
$$\Br_2(L)+\Br_3(L)=\sigma (L)\mod8,$$
where $\Br_3(L)=2(a-b)\mod8$ if $\discr_3L=a\dip+b\din$. \qed
\endcorollary

\subsection{Extensions of lattices}\label{Extensions}
Consider an even lattice $L$ and its {\it extension} $M\supset
 L$, that is another even lattice of finite index $[M:L]$.  By means of the lattice pairing every such
$M$ is canonically embedded in $L\otimes Q$, and we consider
as equivalent the extensions having the same image in $L\otimes Q$.
Furthermore,
$L\subset M\subset M^*\subset L^*$ and the subgroup $H=M/L\subset\discr L$
is {\it isotropic}, i.e., $\q_L$ vanishes on $H$.
 Conversely, for any isotropic subgroup $H\subset\discr L$, the
 preimage, $L_H\subset L^*$, of $H$ under the quotient map
 $L^*\to L^*/L$ is an even lattice. This implies the following
(where \ref{discr-extesion}(1), \ref{extensions}, and \ref{extension-automorphism} are shown in \cite{Nik1},
while \ref{discr-extesion}(2) follows from the van der Blij Theorem \ref{van-der-blij}).

\lemma\label{discr-extesion} If  $M$ is an extension of $L$ associated with an isotropic subgroup
$H$ of $\discr L$, then:
\roster\item $\discr M=H^\perp/H$,
where $H^\perp=M^*/L=\{x\in
 \discr L\,|\,xH=0\}$ is {\it the orthogonal complement of $H$};
\item $\Br(L)=\Br(M)$.
 \qed\endroster
\endlemma

\lemma\label{extensions} For any even lattice $L$ the correspondence
between the isotropic subgroups of $H\subset\discr L$ and
the extensions of $L$ is one-to-one, and $\discr L_H=H^\perp/H$.
\qed\endlemma

\lemma\label{extension-automorphism}
Let $M_1, M_2$ be two extensions of $L$ associated with isotropic subgroups $H_1, H_2$ of $\discr L$.
An automorphism $f\:L\to L$ can be extended
to an isomorphism $M_1\to M_2$ if and
only if
the induced automorphism of $\discr L$ maps isomorphically $H_1$ onto $H_2$. \qed\endlemma

\subsection{Gluing of lattices}\label{gol}
Consider
a pair of even lattices $L_1$, $L_2$ and subgroups
$K_i\subset\discr L_i$, $i=1,2$. We say that $\phi\: K_1\to K_2$
is an {\it anti-isomorphism} if it is a group isomorphism such that
$\q_{L_1}(x)=-\q_{L_2}(\phi(x))$ for all $x\in K_1$.
 For any such an anti-isomorphism $\phi$ the
 graph-subgroup $H_\phi=\{x+\phi(x)\,|\,x\in K_1\}\subset K_1+ K_2$
is isotropic in $discr (L_1+ L_2)=\discr L_1+\discr L_2$ and thus
defines an extension of $L_1+ L_2$.
We denote this extension by $L_1+_{\phi}L_2$
and say that the latter is the result of {\it gluing $L_1$ with
$L_2$ along $\phi$}.

 Recall that a sublattice $L_1\subset L$ of a lattice $L$ is
 called {\it primitive}, if the group $L/L_1$ contains no torsion.
 It is trivial for instance, that the orthogonal complement
 $L_2=L_1^\perp=\{x\in L\,|
 xy=0\ \,\text{for all}\, y\in L_1\}$ of any sublattice
 $L_1\subset L$ is primitive, and primitivity of $L_1$ is
 equivalent to that
$L_1=L_2^\perp$. Note that for each $ i=1,2$
the image of $L$ by the orthogonal projection
to $L_i\otimes\Q $
is contained in $L_i^*\subset L\otimes\Q$ and
the kernel of the composition $L\to L_i^*\to L^*_i/L_i$ is
$L_1+ L_2$, so that
there appear two well defined induced monomorphisms
$p_i\:L/(L_1+ L_2)\to\discr L_i$. Thus, we get the following (see  \cite{Nik1}).

\proposition\label{gluing-lattices}
For any gluing $L=L_1+_{\phi}L_2$, the lattices $L_1$ and $L_2$ are orthogonal complements of
each other, and thus primitive, in $L$. Conversely, if two even sublattices $L_1, L_2$ of a lattice $L$ are orthogonal complements
of each other, then:
 \roster
 \item
they determine canonically
subgroups $H_i\subset\discr L_i$, $i=1,2$, and an anti-isomorphism
$\phi\:H_1\to H_2$,
so that $L$ can be identified with $L_1+_\phi L_2$ by an isomorphism
identical on $L_1$ and $L_2$;
 \item the above subgroups $H_i$ are
 nothing but the images,
 $p_i(L/(L_1+ L_2))$, and $\phi(p_1(x))=p_2(x)$ for all $x\in
 L/(L_1+ L_2)$. \qed
\endroster
\endproposition

\corollary\label{discr-throughH} Gluing $L=L_1+_\phi L_2$ is an even lattice with
$\discr L=H_\phi^\perp/H_\phi$. \qed\endcorollary

Corollary \ref{discr-throughH} together with Lemma \ref{split-out} imply the following.

\proposition\label{discr-gluing}
Let  $L_1+_\phi L_2$ be a gluing along $\phi :K_1\to K_2$.
If  $K_i\subset\discr L_i$, $i=1,2$, are non-degenerate
then
$\discr(L_1+_\phi L_2)=K_1^\perp+K_2^\perp$ . \qed\endproposition

Consider two gluings: $L^{(1)}_1+_{\phi^{(1)}} L^{(1)}_2$ along $\phi^{(1)} :K^{(1)}_1\to K^{(1)}_2$
and $L^{(2)}_1+_{\phi^{(2)}} L^{(2)}_2$ along $\phi^{(2)} :K^{(2)}_1\to K^{(2)}_2$.
We say that
homomorphisms
$f_i\:L^{(1)}_i\to L^{(2)}_i$, $i=1,2$, are
$(\phi^{(1)}, \phi^{(2)})$-compatible, if the induced
homomorphisms $f_i^{\discr}\:\discr
L^{(1)}_i\to \discr L^{(2)}_i$ restricted to $K^{(j)}_i$ commute with $\phi^{(j)}$:
$$\CD K^{(1)}_1 @>f_1^{\discr}>> K^{(2)}_1 \\
@VV\phi^{(1)} V  @V\phi^{(2)} VV \\
K^{(1)}_2 @>f_2^{\discr}>> K^{(2)}_2 \, .
\endCD$$
Lemma \ref{extension-automorphism} immediately implies the
following.

\lemma\label{automorphisms-gluing}
Homomorphisms
$f_i\:L^{(1)}_i\to L^{(2)}_i$, $i=1,2$, can be
extended in a unique way to a homomorphism
$f\:L^{(1)}_1+_{\phi^{(1)}} L^{(1)}_2 \to L^{(2)}_1+_{\phi^{(2)}} L^{(2)}_2$
if and only if $f_i$ are
$(\phi^{(1)}, \phi^{(2)})$-compatible
 If $f_i$ are
$(\phi^{(1)}, \phi^{(2)})$-compatible\
 isomorphisms and $f_i^{\discr}$ are isomorphisms, then
$f$ is also an isomorphism.
\qed
\endlemma

\subsection{The orthogonal complement of
$\pm2$-elements}\label{2-extensions}
Primitive lattice elements $v\in L$
with $v^2=n$ are called {\it
$n$-elements}. For $n\ne0$,  the sublattice $\Z v$ generated by such a
$v$ is isomorphic to $\la n\ra$ and its orthogonal complement
is denoted by $v^\perp$ or $L^v$.

For a given $n$-element $v$ in an even lattice $L$ we have
$L=\la n\ra+_\phi L^v$ with $\phi\:K_v\to K^v$,
$K_v\subset\discr\la n\ra$ and $K^v\subset\discr L^v$. If
$n=\pm2$, there are two cases: both $K_v$ and $K^v$ are trivial, or $K_v=\discr\la\pm2\ra=\la\pm\frac12\ra$ and
 $K^v=\la\mp\frac12\ra$. We say that {\it $v$ is even} in the first case, and that {\it $v$
is odd} in the second.

\lemma\label{evenness-criteria} For any $\pm2$-element $v$ in an even lattice $L$
the conditions below are equivalent:
 \roster\item
 $v$ is even,
 \item
 $\discr L=\la\pm\frac12\ra+\discr L^v$,
 \item
 $L=\Z v+ L^v$,
 \item
 the product $vx$ is even for all $x\in L$,
 \item
 $\frac{v}2\in L\otimes\Q$ lies in
 $L^*$ and its coset $[\frac{v}2]$ is a non-trivial
 element of $\discr L$.
 \endroster\endlemma

\proof Equivalences $(1)\leftrightarrow (2) \leftrightarrow (3)$
and $(4)\leftrightarrow (5)$ are evident. The remaining
equivalence follows from the orthogonal projection formula
$\operatorname{proj}_v x= \frac{vx}{vv}v$.
 \endproof

\lemma\label{h-discriminant}
For any odd
$\pm2$-element $v$ in an even lattice $L$, the complementary discriminant $\discr
L^v$ is isomorphic to $\discr L+\la\mp\frac12\ra$.
\endlemma

\proof
Since $K^v=\la\mp\frac12\ra$ is non-degenerate, it splits off as a direct summand of
$\discr L^v$ by Lemma \ref{split-out}.
\endproof

If an even lattice $L$ is 2-elementary, then its
even $\pm2$-elements are subdivided into two species: {\it
ordinary even elements} and {\it Wu elements}. By definition, an
even element
$v$ is a Wu element if $[\frac{v}2]$ is the characteristic element of
$\discr_2 L$, i.e., $\q_L(x)=x[\frac{v}2] \mod\Z$ for all $x\in\discr_2 L$;
otherwise $v$ is ordinary.

\lemma\label{Wu-and-even} An even $\pm2$-element $v\in L$ in a 2-elementary
even lattice $L$
is a Wu element if and only
if the complementary discriminant $\discr L^v$ is even, i.e.,
$\q_L(x)\in\Z/2\Z\subset \Q/2\Z$ for all $x\in \discr L^v$.
\endlemma

\proof Follows from Lemma \ref{evenness-criteria}(3)
and Lemma \ref{Wu-even}.
\endproof

\subsection{Involutions via gluing}\label{via_gluing}
Consider an even lattice $L$, a lattice involution $c\:L\to L$,
and its {\it eigenlattices} $L_\pm=\{x\in L\,|\,c(x)=\pm x\}$.
Note that $L/(L_++L_-)$ is an elementary 2-group, since $2x=(x+c(x))+(x-c(x))$, $x\pm c(x)\in L_\pm$, for all $x\in L$. Let $r_2(L,c)$ denote the 2-rank of $L/(L_++L_-)$.
Consider the projection $L\to L_\pm^*$ sending $x\in L$ to $y\mapsto xy$, $y\in L_\pm$,
denote by $q_\pm\:L\to\discr L_\pm$ its composition with the quotient map $L_\pm^*\mapsto\discr L_\pm$
and put $K_\pm=q_\pm(L)\subset\discr L_\pm$.
Clearly, each of $q_\pm$ induces a group isomorphism between $L/(L_++L_-)$ and $K_\pm$;
in particular, they give rise to a canonical isomorphism $\phi\:K_+\to K_-$.

The following
proposition showing how $c$ can be described in terms of 2-elementary subgroups of $\discr L_\pm$ is essentially Proposition 1.2.1 in \cite{Nik83}.

\proposition\label{gluing-involutions}
\roster
\item
A lattice involution $c\:L\to L$
yields a presentation of $L$ as a result of gluing $L=L_++_\phi L_-$
 of its eigenlattices
along an anti-isomorphism $\phi\:K_+\to K_-$ between 2-elementary subgroups $K_\pm\subset\discr L_\pm$.
 \item
 Conversely, if a lattice $L$ is glued from
 even lattices, $L=L_++_\phi L_-$, along an anti-isomorphism $\phi\:K_+\to K_-$
between  2-elementary subgroups $K_\pm\subset\discr L_\pm$, then
there exists
a lattice
involution $c\:L\to L$, for which $L_\pm$ are the
$(\pm1)$-eigenlattices.
\qed\endroster\endproposition

\proposition\label{2-elementaryL-pm}
Assume that $L$ is an even lattice with $\discr_2L=\la\pm\frac12\ra$ and $c\:L\to L$ is a lattice involution.
Then the eigenlattices $L_\pm$ of $c$ are 2-elementary, and $|r_2(L_+)-r_2(L_-)|=1$.
\endproposition

 \proof
Let us first extend $c$ to the
 lattice $L'=L+_\phi\la\mp2\ra$,
 where $\phi\:\discr_2 L\to\discr\la\mp2\ra=\la\mp\frac12\ra$ is the anti-isomorphism.
 Namely, choosing one of the two generators, $h\in \la\mp2\ra$ we consider
 the involution defined in $\la\mp2\ra$
by $h\mapsto \e h$, where $\e\in\{+,-\}$ (i.e., the identity or the anti-identity).
Lemma \ref{automorphisms-gluing} allows to glue the
involution $c$ with the latter to obtain
$c_\e\:L'\to L'$, which is an involution as well.
Since $\discr_2L'=0$, Propositions \ref{gluing-involutions} and \ref{discr-gluing} implies that
the eigenlattices $L'_\pm$ of $c_\e$ are glued into $L'$ along $\psi\:\discr_2L_+'\to\discr_2L'_-$,
and that $\discr_2L'_\pm$ are 2-elementary.
 On the other hand, $h\in L'_\e$, thus $L_{-\e}=L'_{-\e}$, and applying Lemmas \ref{evenness-criteria}(2) and
\ref{h-discriminant} we conclude that either $\discr L_\e=\discr L_\e+\la\pm\frac12\ra$ (if $h$ is odd)
or $\discr L_\e+\la\mp\frac12\ra=\discr L_\e$ (if $h$ is even).
\endproof

Given a lattice involution $c$ on an even lattice $L$,
denote by $c_2$ the involution induced in $L_2=L\otimes\Z/2$ by $c$ and put
 $$ r_2(L_2,c_2)=\rank(L_2/L_2^{c_2}),  \text{ where } \   L_2^{c_2}=\{x\in L_2\,|\,c(x)=x\}. $$

\proposition\label{even_twist}
If $L$ has odd discriminant, then
\roster\item
$\delta_2(L_+)=0$ if and only if $(L,c)$ is of type I.
\item
$r_2(L,c)=r_2(L_2,c_2)=r_2(L_\pm)$.
\endroster
\endproposition

\proof
The preimage of $L_2^{c_2}$ under the reduction homomorphism $L\to L\otimes\Z/2=L_2$
is
$L_+\oplus L_-$, which implies $r_2(L,c)=r_2(L_2,c_2)$.
 Each element $a\in\discr_2 L_+$ is represented by an element of the form
$\frac12(x\pm cx), x\in L$,
 which implies $r_2(L,c)=r_2(L_\pm)$. To prove (1), it remains to notice that
$\q(a)=(\frac12(x+ cx))^2=\frac12(x^2+ x\cdot cx)=\frac12 x\cdot cx\in \Q/\Z$.
\endproof

\subsection{Stability}
An even lattice $L$  is called  {\it stable} if any other even lattice $L'$ having the same
inertia indices, $r_+(L')=r_+(L)$ and $r_-(L')=r_-(L)$, and
isometric discriminants, $(\discr L', \q_{L'})=(\discr L, \q_L)$,
is isomorphic to $L$
(such a stability is equivalent to what is also phrased as "uniqueness in its genus", see \cite{Nik83}).

Let us call an even lattice $L$ {\it
epistable} (respectively, {\it $p$-epistable}) if any
automorphism of $\discr L$ (respectively, of $\discr_pL$) is
induced by some automorphism of $L$.

The following Nikulin's criterion
shows that such complications as non-stability or
non-epistability may happen only
if the lattice has the extremal or next to extremal value, $r$ or $r-1$,
(cf. Proposition \ref{p-divisible}) of the ranks $r_p$ for some prime $p$.

\theorem\label{stability-criterion} {\rm (Nikulin \cite{Nik1})} Assume that a lattice
$L$ of rank $r$
is even, indefinite, and the
ranks $r_p$ satisfy the following conditions:
 \roster\item
$r_p\le r-2$ for all primes $p\ne2$;
 \item
if $r_2=r$, then $\discr_2 L$ contains
$\uu$ or $\vv$ as a direct summand.
\endroster
 Then $L$ is both stable and epistable (and in particular, p-epistable for all $p$).\qed
\endtheorem

Note that, as it follows from Theorem \ref{2-classification},
the condition (2) is satisfied, if $\discr_2 L$ is 2-elementary
and $r_2=r>2$, so the condition (2) requires analysis only for $r_2=r=2$.
 Note also that any even lattice $L$ of rank $1$ is stable, since $L=\la 2n\ra$, where $n$
is determined by $\discr L$.

R. Miranda and D. Morrison \cite{MM1}, \cite{MM2} developed further Theorem \ref{stability-criterion} and gave
a necessary and sufficient criterion of stability and epistability, which is in the special case of our interest looks
as follows.

\proposition\label{MM-23elementary}
Suppose that $L$ is an indefinite even lattice of rank $r\ge3$, which has only discriminant factors 2 and 3, and the
latter ones are elementary. Then $L$ is stable and epistable, except possibly the case $r_2=r_3=r$.
\qed\endproposition

\subsection{Different involutions with the same eigenlattices}
Consider a pair of even lattices $T^1$ and $T^2$ with
lattice involutions, $c_i\:T^i\to T^i$, whose eigenlattices
are isomorphic, $T^1_\pm=T^2_\pm$.
By Proposition \ref{gluing-involutions}, $T^i=T^i_++_{\phi^i}T^i_-$, where $\phi^i\:K_+^i\to K_-^i$, $i=1,2$,
are anti-isomorphisms between elementary 2-groups
$K_\pm^i\subset\discr_2 T^i_\pm$.
 We say that involutions $c_1$ and $c_2$ are {\it conjugate via some isomorphism $f\:T^1\to T^2$} if
$f\circ c_1=c_2\circ f$.
 The description of isomorphisms via gluing in Lemma \ref{automorphisms-gluing} yields easily the following.

\proposition\label{conjugate-involutions}
The following conditions are equivalent:
\roster\item
$c_1$ and $c_2$ are conjugate via some isomorphism $f\:T^1\to T^2$;
\item
there exists an isomorphism $f\:T^1\to T^2$ which
maps isomorphically $T^1_\pm$ onto $T^2_\pm$;
\item
there exist isomorphisms $f_\pm\:T^1_\pm\to T^2_\pm$ such that the induced ones,
$f_\pm^{\discr_2}\: \discr_2T^1_\pm\to\discr_2 T^2_\pm$,
map isomorphically $K_\pm^1$ onto $K_\pm^2$ so that the following diagram commutes.
$$\CD
K_+^1    @>f_+^{\discr_2}>>           K_+^2\\
@V\phi^1VV                 @V\phi^2VV \\
K_-^1    @>f_-^{\discr_2}>>           K_-^2 &\ \ \ \qed\\
\endCD$$
\endroster
\endproposition

The next Proposition gives a criterion for involutions
to be conjugate under additional assumption that $\discr_2T^i$ is $\Z/2$ as a group.

\proposition\label{conjugate-prop}
Assume that $c_i\:T^i\to T^i$, $i=1,2$, are
lattice
involutions, whose eigenlattices $T^i_\pm$ are respectively isomorphic,
namely $T^1_\pm=T^2_\pm$. Assume, in addition, that
$r_2(T_+^i)<r_2(T_-^i)$ and
$\discr_2 T^i=\la\e\frac12\ra$,$\e\in\{+,-\}$,
for each  $i=1,2$.
Then $c_1$ and $c_2$ are conjugate via some isomorphism $f\:T^1\to T^2$ if
any one of the following conditions is satisfied:
\roster\item
the lattice $T_-^1$ is 2-epistable;
\item
the lattice $T_+^1$ is 2-epistable and $\Aut(T_-^1)$ acts transitively on
the subgroups of $\discr_2 T_-^1$ anti-isomorphic to $\discr_2T_+^1$.
\endroster
\endproposition

\proof
By Proposition \ref{gluing-involutions}, we have
$T^i=T^i_++_{\phi^i}T^i_-$ with $\phi^i\:K^i_+\to K^i_-$, $i=1,2$.
According to Proposition \ref{2-elementaryL-pm},
$T^i_\pm$ are $2$-elementary, $K^i_+=T^i_+$ (since $r_2(T_+^i)<r_2(T_-^i)$
and $r_2(T^i)=1$), and $K^i_-\subset T^i_-$ is a subgroup of corank 1.
Note also that the orthogonal complements of $K_-^i$ in
$\discr_2T^i_-$ are isomorphic, since the enhanced-group
structures on $\Z/2$ are determined by the Brown invariant.
Thus, $K_-^1$ is sent to $K_-^2$ by some isomorphism
$f_-^{\discr_2}\:\discr_2 T_-^1\to \discr_2 T_-^2$. Moreover, we
can choose it so that the diagram in
Proposition \ref{conjugate-involutions}(3) commutes.
Namely, in the case (1), the  2-epistability of $T_-^1\cong T^2_-$ implies existence of an isomorphism
$f_-\:T_-^1\to T^2_-$ inducing $f_-^{\discr_2}$. Then, Proposition
\ref{conjugate-involutions} implies that $c_1$ and $c_2$ are
conjugate via $f$ defined by $f_+=\id$ and $f_-$ constructed
above.
 In the case (2), we can find $f_-\:T_-^1\to T^2_-$ such that the induced map in $\discr_2(T_-^1)$ sends
  $K_-^1$ to $K_-^2$. Then we can use the epistability of $T_+^1$
 to construct $f_+\:T_+^1\to T^2_+$ which is compatible with $f_-$, that is the diagram like in Proposition \ref{conjugate-involutions} commutes, and we again conclude that $c_1$ and $c_2$ are conjugate.
\endproof

\section{Topology and arithmetics of the covering K3-surfaces}\label{Chapter_CoveringK3}

\subsection{Covering K3 after desingularization}\label{Covering_K3}
In addition to the cuspidal K3-surface $Y$ introduced in Section
\ref{cuspidal-K3} and obtained by taking the double covering of
$P^2$ ramified along a Zariski curve $A$, we take also into
consideration
the non-singular K3 surface $\til Y$ obtained by the minimal
resolution of the six cusps of $Y$. Note that $\til Y$ inherits
from $Y$
a pair of complex conjugations
that differ by
the deck transformation of the double covering $\til Y\to
\til P^2$ of the plane blown-up at the six cusps,
$\til P^2\to P^2$.
This covering is ramified along the proper transform of the Zariski
curve
and fits in a commutative diagram
$$
\CD
\til Y   @> >>       Y\\
@V\til\pi VV                 @V\pi VV \\
\til P^2  @> >>           P^2 \, .\\
\endCD
$$
 Like for $Y$, we give preference to that complex conjugation
whose real locus, $\til Y(\R)$, is projected to $\RR_n$,
call it the {\it M\"obius involution {\rm(or {\it M\"obius real structure})} in $\til Y$} and denote by $\conj=\conj_{\til Y}$.

 The K3-lattice $L=H_2(\til Y)$  contains a sublattice $6\A_2$
spanned by the twelve exceptional
 divisors of the resolution and
the {\it polarization class} $h$, $h^2=2$, which is represented by
the pull-back of a line in $P^2$. In what follows, we  work with
a natural "gluing" reconstruction
of $L$
from two complementary sublattices:
the primitive closure $S\subset L$ of the sublattice spanned by
$6\A_2$ and $h$, and the orthogonal complement
$T=S^\perp=\{x\in
L\,|\,xS=0\}$. We will consider also the sublattice $S^0=\{x\in
S\,|\,xh=0\}\subset S$
and its orthogonal complement
$T'=(S^0)^\perp$.
 All the sublattices, $S$, $T$, $S^0$, and $T'$
are invariant with respect to the complex conjugation involution
$c=\conj_*\:L\to L$, and we denote by $L_\pm=\{x\in L\,|\,c(x)=\pm
x\}$, $S_\pm=L_\pm\cap S$, $T_\pm=L_\pm\cap T$, $S^0_\pm=S^0\cap
L_\pm$, and $T'_\pm=T'\cap L_\pm$ the corresponding eigenlattices.
 Note that $T'_+=T_+$ and $S^0_+=S_+$,
 since $c(h)=-h$ and thus $h\in T'_-\cap S_-$.
It follows also that $T_-$ (respectively, $S_-^0$) is the orthogonal complement of $h$ in $T_-'$
(respectively, in $S_-$).

In particular, we obtain a relation to the ranks $r_\pm$ introduced in
Section \ref{cuspidal-K3}.

\lemma\label{ranks-sing-nonsing}
 $$\align
 r_+=&\rank T_+,\\
 r_-=&\rank T_-+1,\\
 \endalign$$
\endlemma
\proof
This follows
from that  $H_*(Y;\Q)=H_*(\til Y;\Q)/(6\A_2\otimes\Q)= \Q h + T\otimes \Q$ and $h\in L_-$.
\endproof

\subsection{Deficiency in the Smith inequalities}
As it follows from the Smith theory applied to the complex
conjugation involution, the relations
$b_*(X(\R);\Z/2)\le b_*(X;\Z/2)$ and $b_*(X(\R);\Z/2) =
b_*(X;\Z/2) \mod 2$ (here, as before, $X$ and $X(\R)$ stand for
the set of complex and real points respectively) hold for any complex algebraic variety $X$ defined over $\R$.
 The variety $X$ is called {\it M-variety} if  $b_*(X(\R);\Z/2)=b_*(X;\Z/2)$,
 and {\it (M-d)-variety} where $d=\frac12(b_*(X;\Z/2)-b_*(X(\R);\Z/2))$
 otherwise (for instance, one speaks on (M-d)-curves, (M-d)-surfaces, etc.).

First of all, let us compare $d(A)=5-\ell(A)$ of a real
Zariski sextic $A$ and
$d(Y)=\frac12(12-b_*(Y(\R);\Z/2))$ of the double covering
K3-surface $Y$.

 \lemma\label{defects-are-the-same}
 If we choose in $Y$ the M\"obius real structure, then $d(Y)=5-\ell(A)$.
Otherwise (for the non-M\"obius real structure), $d(Y)=6-\ell(A)$.
 \endlemma

 \proof Since $Y(\R)$ projects to $\RR_\pm$ as an orientation double covering with boundary glued
 to itself via deck transformation, it follows that $b_*(Y(\R);\Z/2)=2b_*(\RR_\pm;\Z/2)$, which is
$2+2\ell(A)$ in the case of $\RR_n$
(i.e., M\"obius real structure), and $2\ell(A)$ in the case of
$\RR_o$.
 \endproof

Since the links of cusps are $\Z/2$-homology spheres, the variety
$Y$ is a $\Z/2$-homology manifold.
 Thus, its
 $\Z/2$-valued intersection form
  $\la x,y\ra$, $x,y\in
H_2(Y;\Z/2)$,
is well-defined and non-degenerate.
One can also twist the intersection form via the involution
 $c\:H_2(Y;\Z/2)\to H_2(Y;\Z/2)$ and define
$\la x,y\ra_c=\la x,c(y)\ra$.

We let $H_2^{c}(Y;\Z/2)=\{x\in H_2(Y;\Z/2)\,|\,c(x)=x\}$
and put, similarly to the notation in Section \ref{via_gluing},
$r_2(H_2(Y;\Z/2),c)=\rank H_2(Y;\Z/2)/H_2^{c}(Y;\Z/2).$

\lemma\label{dY=r2} If $c\:H_2(Y;\Z/2)\to H_2(Y;\Z/2)$ is induced by a complex
conjugation in $Y$ such that $Y(\R)\ne\oo$, then $d(Y)=r_2(H_2(Y;\Z/2),c)$.
\endlemma

\proof Since $H_1(Y;\Z/2)=H_3(Y;\Z/2)=0$ and
the Smith sequence
$$\aligned
\dots\to &H_{r+1}(Y/\conj,Y(\R);\Z/2)\to H_{r}(Y/\conj,Y(\R);\Z/2)\oplus
 H_{r}(Y(\R);\Z/2)@>\operatorname{tr^r}+\operatorname{in_r}>>\\
 \to&H_{r}(Y;\Z/2)@>\operatorname{pr_r}>> H_{r}(Y/\conj,Y(\R);\Z/2)\to
\dots
\endaligned
$$
is exact, it is sufficient to show that the
image of $\operatorname{tr}^2+\operatorname{in}_2$ is equal to
$H_2^{c}(Y;\Z/2)$.
Such an equality follows from the
relation $\operatorname{tr}^*\circ \operatorname{pr}_*=1+c$ and,
again, the exactness of the Smith sequence.
\endproof

\remark{Remark}
If $Y(\R)=\emptyset $ then $d(Y)=r_2(H_2(Y;\Z/2),c)-2$. Here, one can bypass the Smith theory and argue
a bit differently.
Namely, since the quotient $Y/\conj$
is a $\Z/2$-manifold, the image of $H_2(Y/\conj;\Z)/\operatorname{Tors}$ by the Gysin homomorphism
in $H_2(Y;\Z)/\operatorname{Tors}$ is a lattice of the form $L'(2)$ where $L'$ has an odd discriminant. Therefore, this lattice is primitively embedded in $H_2(Y;\Z)/\operatorname{Tors}$,
which implies that $r_2(H_2(Y;\Z), c)=\rank L'$. It remains to notice that Proposition \ref{even_twist}(2) implies $r_2(H_2(Y;\Z), c)=r_2(H_2(Y;\Z/2), c)$,
and that \ $\rank L'= \dim H_2(Y/\conj;\Q)=\frac12 \dim H_*(Y;\Q)-2 =\allowbreak
\frac12 \dim H_*(Y;\Z/2)-2$.
\boxedR\endremark

\lemma\label{transfer_toT}
The Gysin homomorphism
$\operatorname{\rho^!} \:H_2(Y;\Z/2)\to H_2(\til Y;\Z/2)$ induced by
the projection $\rho\:\til Y\to Y$
is a monomorphism and its image is
$T'\otimes\Z/2=T'/2T'\subset L/2L=L\otimes\Z/2=H_2(\til Y;\Z/2)$.

The isomorphism $H_2(Y;\Z/2)\to T'/2T'$ provided by
$\rho^!$ commutes with the involutions induced by the complex conjugation
in $Y$ and $\til Y$ and preserves the $\Z/2$-valued intersection forms.
It sends $H_2^c(Y;\Z/2)$ to $(T'_++T'_-)\otimes\Z/2$ and induces an isomorphism
$$H_2(Y;\Z/2)/H_2^c(Y;\Z/2)\to (T'/(T'_++T'_-))\otimes\Z/2=T'/(T'_++T'_-).$$
\endlemma

\proof It follows from $\rho_*\circ\rho^!=\id$ and absence of $2$-torsion in $H_2(Y)$ and $H_2(\til Y)$.
\endproof

\corollary\label{Mobius-defect}
For any real Zariski sextic $A$ and M\"obius involution in $\til Y$ we have
$r_2(T_+)=r_2(T',c)=d(Y)=d(A)$.
If $A(\R)\ne\oo$, and we choose a non-M\"obius involution in $\til
Y$, then $r_2(T_+)=r_2(T',c)=d(A)+1$.\qed
\endcorollary

\proof
From Proposition \ref{even_twist} it follows that $r_2(T'_+)=r_2(T',c)=r_2(T'_2,c)$.
Since $T'_+=T_+$ and $r_2(T'_2,c)=
r_2(H_2(Y;\Z/2),c)$, the remaining parts of the statement follow from Lemmas \ref{defects-are-the-same}, \ref{dY=r2}, and \ref{transfer_toT}.
\endproof

\subsection{Real varieties of type I}
The following is a version of the so-called Arnold lemma, which is valid for any
real algebraic surface $Y$
that is a $\Z/2$-homology manifold
and, in particular, in our case
of the double plane $Y\to P^2$ branched along Zariski sextic.

\lemma\label{Arnold_lemma} The fundamental class $[Y(\R)]\in
H_2(Y;\Z/2)$ is the characteristic class of the bilinear form $\la
x,y\ra_c$, that is, for any $x \in H_2(Y;\Z/2)$ we have
$\la x,x\ra_c=\la x,[Y(\R)]\ra_c$.
\endlemma

\proof
We pick a  $\Z/2$-cycle $C$ representing an element $x \in
H_2(Y;\Z/2)$, which is smooth at the non-singular points of $Y$
and generic with respect to $Y(\R)$, and then check that
any intersection point of $C$ with $Y(\R)$ gives an odd
contribution to both $\la x,c(x)\ra$ and $\la x,[Y(\R)]\ra$.
At the smooth points of $Y(\R)$ genericness means transversality and this contribution is $1$.
At the singular points, the intersection number can be
replaced by the linking number (here, genericness means that the links of $C$, $c(C)$, and $Y(\R)$ are disjoint).
Namely, in a $3$-link of each real singular point
we can fill the local link, $\operatorname{lk}(C)$,  of $C$  by a smooth
$2$-membrane, $M_C$, $\operatorname{lk}(C)= \partial
M_C$,  transversal to its conjugate, $\operatorname{lk}(c(C))=\partial c(M_C)$,
and notice that $M_C\cap c(M_C)$ is a $c$-equivariant $1$-manifold bounding the $0$-cycle $(C\cap c(M_C)) \cup (c(C)\cap M_C)$. Finally, it remains to observe that each non-closed connected component of this $1$-manifold fixed by
$c$ intersects $Y(\R)$ at one point, each closed connected component fixed by $c$ intersects $Y(\R)$ at two points, while the other (non fixed) components come in pairs and are disjoint from $Y(\R)$.
 This implies that the linking number of $\operatorname{lk}(C)$ with $\operatorname{lk}(Y(\R))$ has the same parity as with $\operatorname{lk}(c(C))$.
 \endproof

\corollary\label{Arnold_corollary}
The following conditions are equivalent. \roster\item The
fundamental class $[Y(\R)]\in H_2(Y;\Z/2)$ vanishes. \item The
c-twisted intersection form is even, that is $\la x,x\ra_c=0$ for
all $x\in H_2(Y;\Z/2)$. \qed
\endroster
\endcorollary

A real algebraic surface
$(Y,c)$ is said to be of type I if the conditions of
Corollary \ref{Arnold_corollary}
are satisfied, otherwise, we say that it is of type II.

\lemma\label{even-versus-type}
The following conditions are equivalent.
\roster\item
The pair $(Y,c)$ is of type I;
\item
The pair  $(T',c\vert_{T'})$ is of type I;
 \item
The discriminant form in $\discr_2(T'_+)=\discr_2(T_+)$
is even, or in the other words, $\delta_2(T_+)=0$.
\endroster
\endlemma

\proof
Since the Gysin homomorphism $\operatorname{\rho^!} \:H_2(Y;\Z/2)\to H_2(\til Y;\Z/2)$
commutes with $\conj$, preserves the $\Z/2$-intersection form and, by Lemma \ref{transfer_toT},  establishes
an isomorphism between $H_2(Y;\Z/2)$ and
$T'/2T'$, the equivalence between (1) and (2) follows from Corollary \ref{Arnold_corollary}.
The equivalence between (2) and (3) follows from Proposition \ref{even_twist}.
\endproof

 \lemma\label{Rokhlin_lemma}
 A Zariski sextic
 $A$ is of type I if and only if $Y$ is of type I with respect to the M\"obius real structure.
\endlemma

\proof
Due to Corollary \ref{Arnold_corollary} and
exactness of the Smith sequence written for the deck transformation of the double covering $Y\to P^2$, it is sufficient to check
that $[\Cal A_n]\in H_2(P^2,A;\Z/2)$ is zero if $A$ is of type I (note that it is non zero, otherwise). But if $A$ is of type I, one can apply Rokhlin's trick:
consider instead $[\Cal A_n]+[\R P^2]=[\Cal A_o]$, lift the latter up to an integer homology cycle $R=[\Cal A_o]+[A_\pm]\in H_2(P^2)$ (here, $A_\pm$
denote the two components of $A\sm A(\R)$)
and check
that $R=\frac 12(R-\conj_* R)=\frac12 [A]\in H_2(P^2)$, while $\frac12[A]$ reduced modulo $2$ is equal to $[\R P^2]\in H_2(P^2;\Z/2)$.
\endproof

\proposition\label{type-via-lattice} A Zariski sextic $A$ is of type I if and only if $\delta_2(T_+(\til Y))=0$,
where $\til Y$ is endowed with the M\"obius real structure.
\endproposition

\proof It is a straightforward consequence of Lemmas \ref{Rokhlin_lemma} and \ref{even-versus-type}.
\endproof

\subsection{Resolution decoration of $H_2(\til Y)$}
 Given a {\it K3-lattice} $L= 3\U+2\E_8$, we say that
 $(\D,h)$ is a  {\it decoration of $L$} and
that $L$ is {\it $(\D,h)$-decorated}, if
$\D\subset L$ is a set of twelve elements called the {\it
exceptional classes} forming a basis of $6\A_2$ (i.e.,
$\D=\{e_1',e_1'',\dots,e_6',e_6''\}$, $(e_i')^2=(e_i'')^2=-2$,
$e_i'e_i''=1$, and $e_i'e_j'=e_i'e_j''=e_i''e_j''=0$ for $1\le
i\neq j\le6$) and $h\in L$ is an element orthogonal to $\D$ such
that $h^2=2$.

In the case of the K3-surface $\til Y$ obtained
by resolution of singularities  of the double plane $Y\to P^2$ ramified along a
Zariski sextic $4p^3+27q^2=0$,  such kind of decoration appears naturally in $L=H_2(\til Y)\cong H^2(\til
Y)$: the exceptional curves lying above the six
$\A_2$-singularities of $Y$ provide the exceptional classes
$e_i',e_i''\in\D$, and $h\in L$ is the polarization class of the
projection $\til Y\to Y\to P^2$ (it can be geometrically thought
of as the pull-back of a generic line in $P^2$). We call
this decoration a {\it resolution decoration}.

This resolution, and hence the decoration, can be enriched by a
consideration of the pullback
of the conic $p=0$. Indeed, the
resolution of the cusps together with the pullback
of the conic
can be obtained in two steps: first, blowing up of $P^2$ at the
six cups, and, second, taking the double covering ramified along
the proper transform of the
Zariski sextic. As a result, if the conic
is nonsingular, we get as its pullback
the two
$(-2)$-curves that are proper transforms of the conic and
the twelve exceptional curves $E_i',
E_i''$. Let us denote these two
$(-2)$-curves by $Q', Q''$.  Under appropriate ordering, the
Dynkin-Coxeter graph of this collection of
curves looks as on Figure \ref{delta-graph}a).
\midinsert \hskip10mm\epsfbox{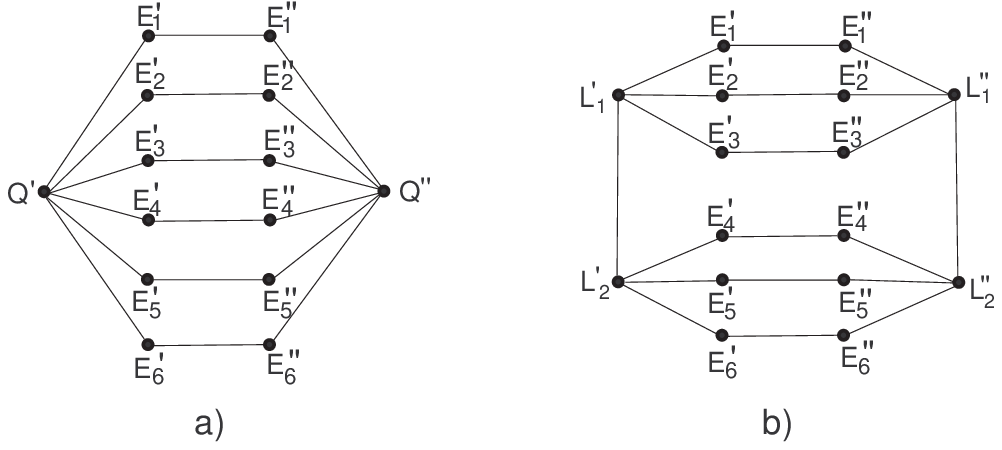} \figure
\label{delta-graph}
\endfigure
\endinsert
If the conic is singular, it splits in a union of two distinct
lines, each containing three of the six cusps. The pullback of
each line splits into two $(-2)$-curves. Let us denote these two
pairs by $L'_1,L''_1$ and  $L'_2, L''_2$. Under an appropriate
ordering, $L'_1$ intersects $L'_2$ and $L''_1$ intersects $L''_2$.
We fix such an ordering and put $Q'=L'_1+L'_2$, $Q''=L''_1+L''_2$.
The Dynkin-Coxeter graph of the collection  $L'_1,L'_2, L''_1,
L''_2,
E'_1,\dots  E''_6$ is as on Figure \ref{delta-graph}b).

\lemma\label{Q-formula} For any Zariski sextic, the following relations hold in  $L=H_2(\til Y)$:
$$
[Q']=h-\sum\frac{e''_i+2e'_i}3, \quad [Q'']=h-\sum\frac{e'_i+2e''_i}3.
$$
\endlemma

 \proof Since the both sides in each of the identities to be proved
belong to sublattice $S$,
the result follows just from a straightforward verification that the both sides
have the same intersection indices with $h, e'_1,e''_1,\dots,e'_6,e''_6$ generating $S$, which in turn follows from
the incidence relations shown on the above Dynkin-Coxeter graphs.
 \endproof

\subsection{The Galois $S_3$-coverings}\label{construction_Galois}
 One of the consequences of Lemma \ref{Q-formula} is that $6\A_2$ is not a
 primitive sublattice of $L$.  A more conceptual explanation of this non-primitiveness consists in appealing to
 the Galois covering $Z\to P^2$ with the Galois group $S_3$ (symmetric group on a set of three elements)
induced by the central projection $\pi_X :X\to P^2$, that is the Galois covering with the Galois group $S_3$ branched along
the Zariski curve $A$ of $\pi_X :X\to P^2$ and whose unramified part is formed by the fibers $Z_s, s\in P^2\setminus A$,
consisting of the six orderings
of the 3-elements sets $\pi_X^{-1}(s)$ (in algebraic terms, it is the ramified Galois covering whose Galois group is the
monodromy group of $\pi_X$). This covering has the following remarkable properties:
\roster\item
$Z$ is a nonsingular $K3$-surface with a holomorphic $S_3$-action;
 \item
the subgroup $C_3\subset S_3$ formed by even permutations acts on $Z$
symplectically (i.e., trivially on the holomorphic differential $2$-forms);
 \item
the quotient $Z/C_3$  is canonically identified with $Y$;
 \item
the $C_3$-action has precisely six fixed points
and they are the pullback of the six cusps of $Y$.
\endroster

In particular, we obtain the following diagram of projections
$$\CD
Z @>>> Y\\
@VVV @VV\pi_YV \\
X @>\pi_X>> P^2.
\endCD$$

Just an existence of such a nontrivial cyclic order three covering $Z\to Y$
implies that $6\A_2$ is not primitive in $L$.
Namely, we obtain the following well-known property.

 \lemma\label{primitive-closure}
 The primitive closure $S^0$ of $6\A_2$ in $L$ is an extension of index $3$.
 \endlemma

 \proof Since the discriminant order $|\discr6\A_2|=3^6$ is a power of $3$,
the index $[S^0:6A_2]$ is also a power of $3$
 and, therefore, coincides with the order of  $H^1(Y\sm {\operatorname{Sing}}\, Y; \Z/3)$.
 The latter group is nontrivial, since
a connected $C_3$-covering of $Y\sm {\operatorname{Sing}}\,Y$ is provided by $Z\to Y$. It is isomorphic to $\Z/3$, since $Z$
is smooth at the six points $z_1,\dots, z_6$ that form the pull back of
${\operatorname{Sing}\,} Y$ and by this reason $Z\sm \{z_1,\dots,z_6\}$ does not admit any further nontrivial $C_3$-covering.
 \endproof

\proposition\label{S}
For any resolution decoration of $L=H_2(\til Y)$,
 \roster\item
the order in
each pair
of $\A_2$-generators $e_i',e_i''\in\D$ can
be chosen in such a way that the
element $\s=
(e_1'-e_1'')+\dots+(e_6'-e_6'')$ is divisible by 3 in $L$;
 \item
 the primitive closure $S^0$ is spanned by $\D\cup\{\frac{\s}3\}$;
 \item
the sublattice $S$ spanned by $\D\cup\{\frac{\s}3,h\}$ in $L$ is
primitive and splits into a direct sum $\la2\ra+S^0$, where
$\la2\ra$ is generated by $h$.
\endroster
\endproposition

\proof Statements (1) and (2) are immediate consequences of Lemmas  \ref{Q-formula} and \ref{primitive-closure}. Statement (3) follows from absence of 2-torsion
in the discriminant of $S^0$.
\endproof

\corollary\label{discr-S} $(1)$ $\discr(S^0)=\dip+3\din$.
\newline
$(2)$
$S$ has the following direct sum decompositions:
$S=
\la2\ra+S^0=
\la2\ra+\E_6+3\A_2
=\U+\A_5+3\A_2.
$
\endcorollary

\proof
 By Proposition \ref{S}(2) and Lemma \ref{discr-extesion}, we have
$\discr S^0=[\frac{\sigma}3]^\perp/[\frac{\sigma}3]$, where
$[\frac{\sigma}3]\in\discr6\A_2=6\la-\frac23\ra$ is the diagonal
element. It follows now from Lemma \ref{3-classification} that
$\discr S^0=p\dip+q\din$, $p+q=4$, where
$2(p-q)=\Br(S^0)=\Br(6\A_2)=4\mod8$, which gives (1).

Part (2) follows from Theorem \ref{stability-criterion}, since
$\discr(S^0)= \la -\frac23\ra + 3 \la \frac23\ra $, $\discr \E_6=\la \frac23\ra  $, $\discr \A_2= \la -\frac23\ra$, and
$\discr \A_5=\la -\frac56\ra = \la \frac12\ra +\la\frac23\ra $.
\endproof

\proposition\label{Galois-conj}
The Galois covering $Z\to P^2$
inherits a real structure, $c_Z:Z\to Z$, from the real structure
of the cubic surface $X$.
This real structure,
$c_Z$, commutes with the Galois action of
$C_3\subset S_3$ on $Z$ and descends to real
structures $c_Y:Y\to Y$, $c_{\til Y}:\til Y\to \til Y$
(as $Y$ is identified with $Z/C_3$)
commuting with the deck transformation of \, $Y\to P^2$ and the
resolution of singularities $\til Y\to Y$;
the real parts $Y(\R)$ and $\til Y(\R)$ are projected to
$\RR_-\subset P^2_{\R}$.

Conversely, if $Y(\R)$ is nonempty and a
real structure $c:\til Y\to \til Y$ commutes with the
involution $\til Y\to\til Y$ induced by the
deck transformation $Y\to Y$, then $c$ lifts to a real structure $c_Z\:Z\to Z$. Moreover:
\roster
 \item
if $c_*(\sigma)=\sigma$, then the above $c_Z$ commutes with the Galois action of $C_3$;
\item
if $c_*(\sigma)=-\sigma$, then
$c_Z$ together with the $C_3$-action form an action of
a semi-direct product, $C_3\rtimes \Z/2=S_3$.
\endroster
\endproposition

\proof The direct statements is a straightforward consequence of
the construction of Galois coverings induced by a given
projection. The assumption $Y(\R)\ne\emptyset$ in the converse
statement is only for ensuring
 that the lift of $c$ is an
involution. Finally, the commutator relation
$c_Z\circ\theta=\theta^{o(c)}\circ c_Z$, where $o(c)$ is defined by
$(c_Z)_*(\sigma)=o(c)\sigma$, follows from the Poincar\'e-Lefschetz
duality (that transforms $\frac13\sigma$ into a characteristic
element of the Galois covering over
$Y\sm{\operatorname{Sing}}\,Y$, cf., proof of Lemma
\ref{primitive-closure}) and the construction of cyclic coverings.
\endproof

\subsection{Abstract K3-lattice conical $(\D,h)$-decorations}\label{conical_decoration}
We say that a $(\D,h)$-decoration of a K3-lattice  $L$
is {\it conical} if it satisfies the properties (1)--(2) (and thus, also (3)\,)
of Proposition \ref{S}. According to this Proposition, each resolution decoration is conical.

The properties (1) and (2) in Proposition \ref{S} imply that, for a given conical $(\D,h)$-decoration of $L$, the choice of $\sigma$ is unique up to sign. We call $\s$ the {\it $\D$-master element}.
Thus, the choice of a $\D$-master element determines the order
$(e_i',e_i'')$ of generators
in each of the $\A_2$-components of $\Delta$, but not the order
of the indices $1\le i\le6$. The choice of the opposite
$\D$-master element alternates the order of $\A_2$-generators.

Given a conically $(\D,h)$-decorated K3-lattice $L$, we denote by
$T$ and $T'$ (in accordance with Section \ref{Covering_K3}) the orthogonal complements of sublattices $S$ and
$S^0$, respectively.

\lemma\label{T}
The following properties hold for any K3-lattice with a conical $(\D,h)$-decoration.
\roster
 \item $T'$ is an extension of $T+\la2\ra$ of index 2. Namely,
$T'$ is obtained by adding to $T$ an element $\frac12(h+x)$ where
 $x\in T$ is a primitive element and
 $x.T\subset2\Z$.
 \item There are following
isomorphisms
$$\aligned
&T=\U+\U(3)+2\A_2+\A_1=\la6\ra+\U+3\A_2=\A_2(-1)+3\A_2+\A_1,\\
&T'=2\U+\U(3)+2\A_2.\endaligned$$
\endroster
\endlemma

\proof Since $\discr(S^0)$ and $\discr(T')$ are isomorphic as groups, the first statement follows from the absence of 2-torsion in $\discr(S^0)$.
The second statement follows from Theorem \ref{stability-criterion}, $\discr(T')=-\discr(S^0)= \la \frac23\ra +3 \la -\frac23\ra $, $ \discr \U(3)= \la \frac23\ra
+ \la -\frac23\ra$,
$\discr \A_2= \la -\frac23\ra$, $\discr(T)=-\discr(S)= \la -\frac12\ra + \discr(T')$, $\discr\la 6\ra=\la \frac16\ra= \la -\frac12\ra + \la \frac23\ra $,
and Lemma \ref{3-classification}.
\endproof

 By an isomorphism between a pair of K3-lattices, $L_1$ and $L_2$,
 decorated with $(\D_1,h_1)$ and $(\D_2,h_2)$ respectively, we mean a
 lattice isometry $f\:L_1\to L_2$ such that $f(\D_1)=\D_2$ and
 $f(h_1)=h_2$.

 We denote by $\Aut(L,\D,h)$ the group of
 automorphisms of a $(\D,h)$-decorated K3-lattice $L$.
 We let also $\overline{\Aut}(L,\D,h)=\{f\,|\,-f\in\Aut(L,\D,h)\}$.

\lemma If a $(\D,h)$-decoration of a K3-lattice $L$ is conical, then
for a $\D$-master element $\s$ and any $f\in \Aut(L,\D,h)$ we have
$f(\s)=\pm\s$.
\endlemma

\proof
Straightforward consequence of the definitions (properties (1) and (2) in Proposition \ref{S}).
\endproof

For any $f$ from $\Aut(L,\D,h)$ or $\overline{\Aut}(L,\D,h)$ we
define $o(f)\in\{+,-\}$
by imposing the relation $f(\s)=o(f)\s$.

\corollary\label{cusps-to-halfs}
Let $c: L\to L$ be induced by the M\"obius involution in $\til Y$.
Then $c\in\overline{\Aut}(L,\D,h)$ and $o(c)=o(A)$, where $A$ is
the Zariski sextic of $\til Y$.
\endcorollary

\proof Due to Proposition \ref{Galois-conj}, if $f':L\to L$ is induced by $c_Z$, then $f'\in\overline{\Aut}(L,\D,h)$ and
$o(f')=+$. Thus, there remains to notice that $o(\RR)=+$ if and only if $p:X(\R) \to P^2(\R)$ is three-fold over the
nonorientable half of $P^2(\R)$, and that the deck transformation of $\til Y\to P^2$ reverses the master element
$\sigma$. The latter follows
from Lemma \ref{Q-formula}, since the deck transformation permutes $Q'$, $Q''$.
\endproof

\proposition\label{decorated-isomorphism} If $L_i$, $i=1,2$ are
conical $(\D_i,h_i)$-decorated K3-lattices, then they are
isomorphic as decorated K3-lattices.
\endproposition

\proof Using an isomorphism for $T$ in Lemma \ref{T}, we obtain a required isomorphism of K3-lattices applying
Lemma \ref{automorphisms-gluing}.
\endproof

\subsection{$S^0$-eigenlattices} As above, let us consider a real Zariski sextic $A$, the M\"obius involution $\conj$ on $\til Y$,
and the induced involution $c\:L\to L$. The latter obviously
preserves each of the
sublattices $S$, $T$, $S^0$, $T'$ invariant. Let us mark with indices $\pm$
(e.g., $S_\pm$, $T_\pm$, etc.), the corresponding $\pm1$-eigenlattices.

If in a lattice $M$ its element $v\in M$ is not divisible by
$d>1$, but $vx$ is divisible by $d$ for all $x\in M$, and $v^2$ is
divisible by $d^2$, then by adding element $\frac{v}d$ to $L$ we
obtain its extension (still an integral lattice) of index $d$,
which we denote by $[M]_{\frac{v}d}$.
In this notation,
$S^0=[6\A_2]_{\ss}$.

\proposition\label{S-eigenlattices} The
eigenlattices $S^0_\pm$
are determined up to isomorphism by the number $\n_i$ of imaginary
pairs of cusps of $A$ and the sign $o(A)\in\{+,-\}$, as is
indicated in the Table 2.
\midinsert\topcaption{Table 2}\endcaption
\centerline{$\matrix
{\matrix
\boxed{\text{The case of $o=-$}}\\
\boxed{\matrix \n_i  &S_+^0  &S_-^0\\
\text{---}&\text{---------}&\text{---------}\\
0& 0      & [6\A_2]_\ss \\
1& \A_2(2) & [4\A_2+\A_2(2)]_\ss \\
2& 2\A_2(2) & [2\A_2+2\A_2(2)]_\ss \\
3& 3\A_2(2) &[3\A_2(2)]_\ss=\E_6(2) \\
\endmatrix}\endmatrix
\matrix\boxed{\text{The case of $o=+$}}\\
\boxed{\matrix \n_i  &S_+^0  &S_-^0\\
\text{---}&\text{---------}&\text{---------}\\
0& [6\la-6\ra]_\ss& 6\A_1 \\
1&  [4\la-6\ra+\A_2(2)]_\ss&4\A_1+\A_2(2) \\
2& [2\la-6\ra+2\A_2(2)]_\ss&2\A_1+2\A_2(2)   \\
3& [3\A_2(2)]_\ss=\E_6(2)&3\A_2(2)  \\
\endmatrix}\endmatrix}\endmatrix
 $}
\endinsert
\endproposition

\proof
 As follows from Corollary \ref{cusps-to-halfs}, the master element $\s$ defining the extension
$S^0= [6\A_2]_\ss$ is preserved by $c$ if $o=+$, and reversed otherwise. If $o=-$, then each real cusp gives
$\A_2\subset S_-^0$. If $o=+$, then at each real cusp we have $c(e'_i)=-e_i''$, so that
$e'_i+e_i''\in S_-^0$ and $e'_i-e_i''\in S_+^0$.
A pair of imaginary cusps gives a copy of $\A_2(2)$ in $S^0_+$ and
another copy in $S^0_-$.
\endproof

\corollary\label{p-q-of-S}
For any real Zariski sextic and the associated eigenlattices $S^0_\pm$, we have
$$
\discr_3S_+^0=p\din+q\dip, \quad
\discr_3S_-^0=(1-p)\din+(3-q)\dip
$$
where
$p=1$, $0\le q\le3$ if $o(c)=+$ and $p=0$, $0\le q\le3$ if $o(c)=-$.

Under such a presentation, the number $\n_i$ of imaginary pairs of cusps is
$q$ if $p=1$, and $3-q$ if $p=0$.
\qed\endcorollary

A pair of lattices isomorphic to one of the eight pairs $(S^0_+,S^0_-)$ in Table 2 will be called an {\it S-pair}.
Each component, $S^0_\pm$, of such a pair will be called an {\it S-half}.

\subsection{Geometric involutions}\label{geometric-involutions}
We say that a lattice {\it hyperbolic} if its positive inertia
index equals $1$ (with the usual abuse of terminology if the
negative inertia index is $0$).
We say that an involution $c:L\to L$
in a conical $(\D,h)$-decorated
K3-lattice $L$ is {\it geometric} if
 the following conditions are satisfied:
 \roster
 \item
 $c\in \overline{\Aut}(L,\D,h)$;
 \item
 the eigenlattices $T_\pm(c)$ are hyperbolic.
 \endroster

\lemma\label{geometric-properties}
For any geometric involution $c:L\to L$ of a conical $(\D,h)$-decorated
 K3-lattice $L$, the pair of eigenlattices $T_\pm=T_\pm(c)$ has the following properties:
 \roster
 \item
$r(T_+)+r(T_-)=9$;
\item $T_\pm$ have no other discriminant factors than
2 and 3, and the latter ones are both elementary;
\item $|r_2(T_+)-r_2(T_-)|=1$,
\item the eigenlattice whose rank $r_2$ is greater has $\delta_2=1$;
\item $\discr_3T_++\discr_3T_-=\dip+3\din$.
 \endroster
\endlemma

\proof
 Statement (1) follows from $r(S)+r(T_+)+r(T_-)=r(L)$.
Statement (2)
follows from  Lemma \ref{evenness-criteria} applied to $h\in T'$ and Proposition \ref{gluing-involutions}. If $T'_-=\Z h + T_-$ then $\discr_2(T_+) =-\la\frac12\ra +\discr_2(T_-)$ and $r_2(T_+)=r_2(T_-)+1$, otherwise, $T'_-$
is obtained by a nontrivial gluing $\Z h +_\phi T_-$ and then $r_2(T_+)= r_2(T_-)-1$, so that in both cases we get (3) and (4).
Statement (5) follows from $\discr(S^0)= \la -\frac23\ra + 3 \la \frac23\ra $
(see Corollary \ref{discr-S}) and Proposition \ref{gluing-involutions}.
\endproof

Let us call a geometric involution {\it ascending} if $r_2(T_+)<r_2(T_-)$ and {\it descending} otherwise.
Similarly,
a complex conjugation on $\til Y$
is called {\it ascending involution} or {\it ascending real structure} if the induced by it geometric involution is ascending;
otherwise, we call it {\it descending involution} or {\it descending real structure}.

\lemma\label{ascending_descending} Assume that $A$ is a Zariski sextic,
and $c$ is the involution in $L=H_2(\til Y)$
 induced by one of the
two lifting of the complex conjugation from $P^2$ to $\til Y$.
Then:
\roster\item
$c$ is geometric;
\item
if $c$ is induced by the M\"obius real structure and
$A(\R)\ne\oo$, then $c$ is ascending;
\item if $c$ is induced by the non-M\"obius real structure and
$A(\R)\ne\oo$, then $c$ is descending.
\endroster
\endlemma

\proof The involution $c$ is geometric, since each lift of the complex conjugation from $P^2$ to $\til Y$
sends the exceptional divisors and the polarization to the exceptional divisors and the polarization,
reversing their orientation, and as any anti-holomorphic involution, permutes the Hodge summands, $H^{2,0}$ and $H^{0,2}$.
Since $r_2(T_+)=r_2(T'_+)=r_2(T',c)$, Corollary \ref{Mobius-defect} implies
(2) and (3).
\endproof

Let $C(L,\D,h)$ denote the set of geometric involutions, and
let $C^<(L,\D,h)$, $C^>(L,\D,h)$
denote the set of the ascending and descending ones respectively.
The group $\Aut(L,\D,h)$ acts on $C(L,\D,h)$
preserving
the
subsets $C^<(L,\D,h)$ and $C^>(L,\D,h)$ invariant.
 Let $C[L,\D,h]$, $C^<[L,\D,h]$, and $C^>[L,\D,h]$
denote the orbit spaces (i.e., the set of conjugacy classes) of
geometric involutions, ascending ones, and descending ones respectively.
We say that two geometric involutions have the same {\it homological type}, if they represent the same element
in $C[L,\D,h]$.

\subsection{T-pairs and T-halves}\label{T-pairs-T-halves}
We say that a pair of even hyperbolic lattices $(T_1,T_2)$ form a {\it
T-pair} if they satisfy the five properties stated in Lemma
\ref{geometric-properties} for $(T_+,T_-)$. A T-pair will be
called {\it ascending (descending)} if $r_2(T_1)<r_2(T_2)$
(respectively $r_2(T_1)>r_2(T_2)$).

\lemma\label{r2-estimate} For any ascending T-pair $(T_1,T_2)$ or
descending T-pair $(T_2,T_1)$
$$
\align
r_2(T_1)&\le\min(r(T_1),8-r(T_1)),\\
r_2(T_2)&\le\min(r(T_2),10-r(T_2)).\\
\endalign
$$
In particular, $r_2(T_1)\le4$, and $r_2(T_2)\le5$, where $r_2(T_2)=5$ implies that $r(T_1)=r_2(T_1)=4$
and $r(T_2)=5$.
\qed\endlemma

\proof It follows from
properties (1) and (3) combined with
the bounds $r_2(T_i)\le r(T_i)$, see Proposition \ref{p-divisible}.
\endproof

A lattice is called a {\it T-half} if it is a component, $T_1$ or
$T_2$, of some T-pair. The following is a
straightforward consequence
of Lemma
\ref{geometric-properties}.

\proposition\label{T-halves-properties}
If $M$ is a $T$-half, then:
 \roster
 \item $M$ is hyperbolic of rank $1\le r\le 8$;
 \item $M$ has no discriminant p-factors different from $p=2,3$, and the latter ones (if exist) are elementary;
 \item $M$ has either $r_2\le4$, or $r_2=r=5$, and in the latter case, it can be only the second component of an ascending T-pair;
 \item $\discr_3 M=p\dip+q\din$, where $0\le p\le1$ and $0\le q\le 3$.
\qed\endroster
\endproposition

\proposition\label{properties-of-involutions}
For any $c\in C(L,\D,h)$, the  $S^0$-eigenlattices
$(S_+^0,S_-^0)$, the $T$-eigen\-lattices $(T_+,T_-)$, and $T_-'$
must satisfy the following relations.
 \roster
 \item The pair $(S_+^0,S_-^0)$ is an S-pair.
 \item $T_\pm$ are  T-halves and
 either $(T_+,T_-)$ or $(T_+,T_-)$  is an ascending T-pair.
 \item $\discr_3 T_\pm$ is anti-isomorphic to $\discr_3 S_\pm^0$,
 and in particular $r_3(T_\pm)=r_3(S^0_\pm)$.
 \item $\discr_2 T_+$ is anti-isomorphic to $\discr_2
 T_-'$, and in particular $|r_2(T_+)-r_2(T_-)|=1$.
\endroster
 \endproposition
 \proof
Any $c\in C(L,\D,h)$ by definition preserves sublattice $6\A_2\subset S^0$ invariant.
Moreover, an automorphism of
 $6\A_2$ must permute its $\A_2$-components, and then, the arguments
similar to that of Proposition \ref{S-eigenlattices} can be
applied to obtain (1).
(2) follows immediately from definitions.
Since $\discr(L_\pm)$ for the
$\pm$-eigenlattices $L_\pm$ are elementary 2-groups as it follows from
Proposition \ref{gluing-involutions}, gluing of $L_\pm$ from $T_\pm$, and $S_\pm$
involves an anti-isomorphism of their $\discr_3$-components according to Proposition \ref{discr-gluing},
which implies (3).
 Since $\discr T'$ is a 3-group, the $\discr_2$-components of $T_+$
 and $T_-'$ are anti-isomorphic, again by Proposition \ref{discr-gluing}, which implies (4).
 \endproof

Property (3) in Proposition \ref{properties-of-involutions} shows that the pair $(p,q)$ describing
$\discr_3(T_+)$ (or equivalently, $\discr_3(T_-)$) determines and is determined by the pair $r_3(T_+)=p+q$ and
$o(c)$. Namely, it shows that Corollary \ref{p-q-of-S}
can be restated as follows.

\corollary\label{p-q-of-T}
The ranks $0\le p\le1$ and $0\le q\le 3$ in the decompositions
$\discr_3 T_+=p\dip+q\din$, $\discr_3 T_-=(1-p)\dip+(3-q)\din$
determine $o(c)$ and the number $\n_i$ as follows:
 \roster\item
if $p=1$, then $o(c)=+$ and $\n_i=q$;
\item
if $p=0$, then $o(c)=-$ and $\n_i=3-q$.
\qed\endroster
\endcorollary

\lemma\label{sumB2}
For any T-pair $(T_1,T_2)$, $\Br_2(T_1)+\Br_2(T_2)=-1$.
\endlemma

\proof
Applying Corollary \ref{vanderblij-refined} to
$T_i$, $i=1,2$, we
get $\Br_2(T_i)+\Br_3(T_i)=(2-r(T_i))$, where $\Br_3(T_i)=2(p_i-q_i)$ if $\discr_3 T_i=p_i\dip+q_i\din$.
 Using $r(T_1)+r(T_2)=9$ and $\discr_3 T_1+\discr_3 T_2=\dip+3\din$
 we conclude that $\Br_2(T_1)+\Br_2(T_2)=(4-9)-2(1-3)=-1$.
\endproof

\lemma\label{anti-isomorphism-criterion}
Assume that $(T_1,T_2)$ is an ascending  T-pair. Then $K_1=\discr_2 T_1$ is anti-isomorphic to a subgroup
$K_2\subset\discr_2 T_2$ if and only if
$K_2$ is the orthogonal complement $v^\perp$ of an element
$v\in\discr_2T_2$
satisfying the following conditions:
\roster\item
$\q_{T_2}(v)=-\frac12\in\Q/2\Z$;
\item
$v$ is a Wu element if and only if $\delta_2(T_1)=0$.
\endroster
\endlemma

\proof
If $K_2$ is anti-isomorphic to $K_1$, then
$K_2$ is a non-degenerate subgroup of $\discr_2 T_2$ of corank $r_2(T_2)-r_2(T_1)=1$
and, therefore, $K_2=v^\perp$ for some $v\in\discr_2 T_2$. Furthermore, $\Br_2(T_1)=\Br(K_1)=-\Br(K_2)$, and additivity of $\Br$ implies
$\Br_2(T_2)=\Br(K_2)+\Br\la v\ra$.
and $\Br\la v\ra=2\q_{T_2}(v)\mod4$ by Lemma \ref{Wu-even}(3).
 Thus, $\Br_2(T_1)+\Br_2(T_2)=2\q_{T_2}(v)\mod4$, which equals $-1\mod4$ by Lemma \ref{sumB2}
 and thus gives condition (1). Lemma \ref{Wu-even}(2) implies condition (2).

Conversely, if $v\in\discr_2 T_2$ satisfies condition (2), then $\delta_2(K_2)=\delta_2(K_1)$. If $v$ satisfies condition (1), then $\Br\la v\ra=-1$ (see Lemma \ref{Wu-even}(4)), and thus $\Br(K_2)=\Br_2(T_2)+1$.
Applying Lemma \ref{sumB2} we get
$\Br(K_1)=\Br_2(T_1)=-1-\Br_2(T_2)=-\Br(K_2)$
and by means of Theorem
\ref{2-classification} conclude that
$K_1$ is anti-isomorphic to $K_2$.
\endproof

\subsection{T-halves and pairs of M\"obius involutions}
Given a Zariski sextic $A$, let $c_1$ denote the M\"obius
involution in $L=H_2(\til Y)$, and $c_2$
be the non-M\"obius one.

\proposition\label{r2c1c2}\label{exceptional-T-relations} $T_\pm(c_2)=T_\mp(c_1)$ for any Zariski sextic $A$.
Moreover, in the case $A(\R)=\oo$,
\roster
 \item  $T_+(c_1)= T_-(c_2)$ have $r=r_2=5$ and $\delta_2=1$;
  \item  $T_-(c_1)=T_+(c_2)$ have $r=r_2=4$ and $\delta_2=0$.
  \endroster
\endproposition

\proof Note that $c_2\circ c_1=c_1\circ
c_2$ is induced by the deck transformation of the covering $\til Y\to \til P^2$, and the action of the deck transformation
in $H_2(\til Y)$ is equal to
the product of
pairwise commuting maps  $-\rho_h$  and
$\rho_{e_i'+e_i''}$, $i=1,\dots,6$, where $\rho_v$ stands for
the reflection $x\mapsto x-2\frac{xv}{v^2}v$. Therefore, $T_\pm(c_2)=T_\mp(c_1)$.

 In the case of $A(\R)=\oo$, Lemma \ref{codes-and-ranks}
gives values $r_+=r_-=d(A)=5$ for the involution $c_1$, and applying Lemma \ref{ranks-sing-nonsing} we obtain
the values of $r(T_\pm(c_1))=r(T_\pm(c_2))$ indicated in (1)--(2).

Corollary \ref{Mobius-defect} implies that $r_2(T_+(c_1))=5$.
The parity $\delta_2=0$ for $T_+(c_2)$,
and thus for $T_-(c_1)$, follows from Corollary \ref{Arnold_corollary}.
The imparity $\delta_2=1$ for $T_+(c_1)$, and thus for $T_-(c_2)$, follows from Proposition \ref{type-via-lattice},
since $A(\R)=\oo$ implies that $A$ has type II.
Finally, the values of $r_2(T_\pm(c_1))$, and thus of $r_2(T_\pm(c_2))$, indicated in (1)--(2)
follow from Lemmas \ref{defects-are-the-same} and Corollary \ref{Mobius-defect}.
\endproof

\section{Arithmetics of geometric involutions}\label{Chapter_Tpairs}

\subsection{Automorphisms of 3-elementary inner product
groups of small rank}\label{3-Aut}
Here, we analyze the automorphisms of $\discr(S^0)\cong
 \din+3\dip\cong\dip+3\din$ (see Corollary \ref{discr-S}) and of its non-degenerate subgroups.

Note that each permutation of coordinates followed by
an alternation of signs of some components, $(x_1,\dots,x_n)\mapsto(\pm
x_{\s(1)},\dots,\pm x_{\s(n)})$, provides an automorphism of
$n\di{}23$, as well as of $n\di-23$. These "coordinatewise"
automorphisms form a group that we denote by $\Aut_\co(n\di{\pm}23)$. It is
isomorphic to the group of symmetries of an $n$-cube and fits canonically in
an exact sequence $1\to(\Z/2)^n\to \Aut_\co(n\di{\pm}23) \to S_n \to1$.

\lemma\label{3elliptic} If $1\le n\le3$, then
$\Aut(n\di{\pm}23)=\Aut_\co(n\di{\pm}23)$.
\endlemma

\proof
Note that for $q=n\di{\pm}23$ the expression
$q(x_1,\dots,x_n)=\pm\frac23(x_1^2+\dots+x_n^2)$ (mod $2\Z$)
depends only on the
number of non-zero coordinates modulo 3. Thus, for $n\le3$ the
elements of the direct summands in $n\di{\pm}23$ are distinguished
from all the other elements if $n\le3$.
Hence, any automorphism in $\Aut(n\di{\pm}23)$ with $n\le 3$ is coordinatewise.
\endproof

Now, consider an enhanced group $(G,q)=6\din$ and
denote by $\delta=\sum_{i=1}^6a_i$
the {\it diagonal element}, that is the sum of the generators
$a_i$ of all the $\din$-components of $G$. Note that $\delta$ is
{\it isotropic}, that is $\delta^2=0$, and thus $\delta\in
G^\delta=\{x\in G\,|\,x\delta=0\}$.

We put $\bb{i_1\dots i_p}{j_1\dots j_q}=(a_{i_1}+\dots
a_{i_p})-(a_{j_1}+\dots a_{j_q})\in G$, and in the case $\delta
\bb{i_1\dots i_p}{j_1\dots j_q}=-\frac23(p-q)=0$,
we denote by
$[\bb{i_1\dots i_p}{j_1\dots j_q}]$ the coset in $G^\delta/(\delta)$.
Since $\delta$ is isotropic, the quotient group $G^\delta/(\delta)$ inherits from $G$
the inner product and quadratic form.

\lemma\label{elements-description}
 The non-trivial elements of the group
$G^\delta/(\delta)=\din+3\dip$
break up into 3 sets: \roster\item 30
 elements $[\bb{i}{j}]$, $1\le i,j\le6$, $i\ne j$, of square $\frac23$;
 \item 30 elements $[\bb{ij}{kl}]$
 of square $-\frac23$, where $1\le i,j,k,l\le 5$ are distinct;
 \item 20 isotropic (i.e., of square $0$) elements
 $[\bb{ijk}{}]$, $1\le i<j<k\le 6$.
\endroster
\endlemma

\proof Straightforward.
\endproof

\lemma\label{S-discr}
The elements $\bb{12}{34}$,$\bb12$,$\bb34$,$\bb56$
split the enhanced group $G^\delta/(\delta)$ in an orthogonal direct sum
$\din+3\dip$.
\endlemma

\proof Straightforward.
\endproof

\subsection{Reduction homomorphism}
Keeping the notation of the previous Subsection,
we consider the automorphism group  $\Aut G$ and  study
its subgroups:
$$\aligned\Aut(G,\delta)=&\{f\in\Aut(G)\,|\,f(\delta)=\pm\delta\},\\
\Aut_\co(G,\delta)=&\{f\in\Aut_\co(G)\,|\,f(\delta)=\pm\delta\},\\
\Aut_\co^+(G,\delta)=&\{f\in\Aut_\co(G,\delta)\,|\,f(\delta)=\delta\}.\endaligned$$

The induced homomorphism $\Aut(G,\delta)\to
\Aut(G^\delta/(\delta))=\Aut(\din+3\dip)$ will be called the {\it
reduction homomorphism}.

\proposition\label{3discr-isomorphism}
 \roster \item $\Aut_\co^+(G,\delta)=S_6$ and
$\Aut_\co(G,\delta)=S_6\times\Z/2$. \item The reduction
homomorphism restricted to $\Aut_\co(G,\delta)$ is an isomorphism
$S_6\times\Z/2=\Aut_\co(G,\delta)\to\Aut(G^\delta/(\delta))=\Aut(\din+3\dip)$.
\endroster
\endproposition

 \lemma\label{action-on-elements}
$\Aut_\co(G,\delta)$ acts effectively and transitively on the set
of $\frac23$-elements $[\bb{i}j]\in G^\delta/(\delta)$.
\endlemma

\proof
It follows from the description of
$\frac23$-elements in Lemma \ref{elements-description}.
\endproof

 \demo{Proof of Proposition \ref{3discr-isomorphism}}
The claim (1) is straightforward, since $\delta$ is
preserved by a coordinatewise
automorphism if and only if it is defined by a
simple permutation (without any sign reversion).

Lemma \ref{action-on-elements}
implies that the reduction homomorphism is injective on
$\Aut_\co(G,\delta)$. To show that it is isomorphic, it is
sufficient
to check that the order of the group $\Aut(\din+3\dip)$ coincides
with that of $\Aut_\co(G,\delta)$.

The order of $\Aut_\co(G,\delta)$
is equal to $|S_6\times\Z/2|=2\cdot 6!$.
On the other hand, by Lemma \ref{action-on-elements} the group
$\Aut(\din+3\dip)$ acts transitively on the 30 elements
$[\bb{i}j]$. The stabilizer of an element $[\bb{i}j]$ is
isomorphic to $\Aut(3\din)$, since the orthogonal complement of
$[\bb{i}j]$ in $G^\delta/(\delta)=\dis23+3\di-23$ is isomorphic to
$3\din$ (see Lemma \ref{3-classification}).
This implies that the order of $\Aut(\din+3\dip)$ is
$30|\Aut(3\dip)|$, where by Lemma \ref{3elliptic}
$|\Aut(3\dip)|=8|S_3|=48$. \qed
\enddemo

\rk{Remark} In fact, it is not difficult to prove that
$\Aut(G,\delta)=\Aut_\co(G,\delta)$, and, therefore, Proposition
\ref{3discr-isomorphism}(2) means that the reduction homomorphism
is an isomorphism. \boxedR\endrk

Now, assume that $c_G\in S_6\times \Z/2=\Aut_{\co}(G,\delta)$ is an
involution acting on $G$. It induces an involution $c_G^\delta$ in
$\din+3\dip=G^\delta/(\delta)$.
 Let $\Aut_{\co}(G,\delta,c_G)=\{f\in\Aut_{\co}(G,\delta)\,|\,f\circ c_G =c_G\circ f\}$
and $\Aut(\din+3\dip,c_G^\delta)=\{f\in\Aut(\din+3\dip)\,|\,f\circ c_G^\delta=c_G^\delta\circ f\}$.
 As an immediate corollary, we obtain the following equivariant version of Proposition
 \ref{3discr-isomorphism}(2).

\corollary\label{3equi-isomorphism} The isomorphism in Proposition
\ref{3discr-isomorphism}(2) restricts to an isomorphism
$$\Aut_{\co}(G,\delta,c_G)\to\Aut(\din+3\dip,c_G^\delta).\qed$$
\endcorollary

\subsection{Equivariant epistability of $S^0$} The proof of the following Lemma is straightforward.

\lemma\label{discrS0} There is a canonical isomorphism
$\discr(S^0)=G^\delta/(\delta)$.
\qed \endlemma

Given a subset $D\subset L$ in a lattice $L$, we let $\Aut(L,D)=\{f\in\Aut L\,|\,f(D)=D\}$.
We say that $L$ is {\it $D$-relatively epistable} if the induced homomorphism $\Aut(L,D)\to\Aut(\discr L)$ is
surjective.

\proposition\label{Aut-lifting} The induced projection $\Aut(S^0,\D)\to\Aut(\discr
S^0)$ is an isomorphism.
In particular, $S^0$ is $\D$-relatively epistable.
\endproposition

\proof
Each automorphism $f\in\Aut(S^0,\D)$ preserves the sublattice
$6\A_2\subset S^0$, permutes its $\A_2$-components, and preserves or
reverses the $\Delta$-master element $\s\in6\A_2$
that is responsible for the extension
of $6\A_2$ to $S^0$. Thus, the induced automorphism in
$G=\discr(6\A_2)$ preserves or reverses $\delta$. This yields a
homomorphism $\Aut(S^0,\D)\to\Aut_\co(G,\delta)=S_6\times\Z/2$,
which is an isomorphism.
 Now, it is left to apply Lemma \ref{discrS0} and Proposition \ref{3discr-isomorphism}(2).
\endproof

 Let $\overline{\Aut}(S^0,\D)=\{f\in\Aut(S^0)\,|\,-f\in\Aut(S^0,\D)\}$.

\corollary\label{bar-Aut-lifting}
The projection $\overline{\Aut}(S^0,\D)\to\Aut(\discr
S^0)$ is bijective. Thus, any involution in $\discr S^0$ can be lifted to an involution of
$S^0$ which sends $\D$ to $-\D$.
\endcorollary

\proof
Given $\phi\in\Aut(\discr S^0)$, let $f\in\Aut(S^0,\D)$ be
 a lifting of $-\psi$ existing by Proposition \ref{Aut-lifting}. Then $-f\in\overline{\Aut}(S^0,\D)$
 is a lifting of $\phi$ that we need.
\endproof

 Now, consider an involution $c\in C(L,\D,h)$ and the induced involutions, $c_G$ on $G$ and
$c_G^\delta$ on $\discr S^0=G^\delta/(\delta)$.
 Let $\Aut(S^0,\D,c)=\{f\in\Aut(S^0,\D)\,|\,fc=cf\}$, and
$\Aut(\discr S^0,c_G^\delta)=\{\phi\in\Aut(\discr S^0)\,|\,\phi
c_G^\delta=c_G^\delta\phi\}$.

\corollary\label{S-equi-epistability} A restriction of the isomorphism of Proposition
\ref{Aut-lifting} yields an isomorphism
$\Aut(S^0,\D,c)\cong\Aut(\discr S^0,c_G^\delta)$.
\endcorollary

\proof
By definition, $-c$ belongs to $\Aut(S^0,\D)$ and
 $\Aut(S^0,\D,c)$ coincides with the centralizer of $-c$, whereas $\Aut(\discr S^0,c_G^\delta)$ is the centralizer of
its image $-c_G^\delta$ under the isomorphism $\Aut(S^0,\D)\to\Aut(\discr
S^0)$.
\endproof

The property of $S^0$
indicated in Corollary \ref{S-equi-epistability} will be referred to as
 {\it $c$-equivariant $\D$-relative epistability of $S^0$}.
 In a bit more general setting, the definition looks as follows.
Given a lattice $L$ with a subset  $D\subset L$ and an involution $c\:L\to L$ inducing
an involution $c^{\discr}$ in $\discr L$,
we say that $L$ is {\it $c$-equivariant $D$-relative epistable} if
the induced homomorphism from $\Aut(L,D,c)=\{f\in\Aut(L,D)\,|\,f\circ c=c\circ f\}$ to
$\Aut(\discr L,c^{\discr{}})=\{\phi\in\Aut(\discr L)\,|\,\phi\circ c^{\discr{}}=c^{\discr{}}\circ f\}$
is surjective.

\subsection{Gluing of involutions}
If two lattices $L_1$, $L_2$ are glued along different anti-isomorphisms,
$\phi$ and $\phi'$ between the same subgroups $K_i\subset\discr L_i$, $i=1,2$,
Lemma \ref{automorphisms-gluing} gives the following
criterion of existence
of an isomorphism $f\:L\to L'$ with
$f(L_i)=L_i$ for $i=1,2$.

\lemma\label{isomorphic-gluing} Assume that $L=L_1+_\phi L_2$ and
$L'=L_1+_{\phi'} L_2$, where $\phi$ and $\phi'$ are
anti-isomorphisms $K_1\to K_2$ between the same subgroups
$K_i\subset\discr L_i$, $i=1,2$. Then, a pair of automorphisms
$f_i\:L_i\to L_i$, $i=1,2$ can be extended to an isomorphism
$f\:L\to L'$ if and only if $f_i$ are $(\phi,\phi')$-compatible,
that is if and only if the following two conditions are satisfied:
 \roster\item
for each $i=1,2$, the automorphism $(f_i)_*$ induced by $f_i$ on $\discr L_i$ preserves $K_i$;
 \item
 $(f_2)_*\vert_{K_2}\circ\phi=\phi'\circ(f_1)_*\vert_{K_1}$.
 \qed
\endroster
\endlemma

If we let $L'=L$ using one of such isomorphisms for identification,
then all the others can be characterized as follows.

\corollary\label{automorphism-gluing}
 Automorphisms $f\:L\to L$ of $L=L_1+_\phi L_2$, $\phi\:K_1\to K_2$,
  such that $f(L_i)=L_i$, $i=1,2$,
are in one-to-one correspondence with the pairs
$(f_1,f_2)$ of $(\phi,\phi)$-compatible automorphisms of $L_1$ and $L_2$
 (here, the compatibility means that the induced
 homomorphisms $(f_1)_* : \discr L_1\to \discr L_1,(f_2)_*  : \discr L_2\to \discr L_2$
 preserve  $K_1, K_2$ and, being restricted to  $K_1,K_2$, commute with $\phi$).
 \qed
\endcorollary

Next, we give a criterion for existence of an extension of an
automorphism $f_1\:L_1\to L_1$ from $L_1$ to $L=L_1+_\phi L_2$.

\proposition\label{automorphism-extension}
Let $L=L_1+_\phi L_2$, $\phi\:K_1\to K_2$, be a gluing such that
$K_2\subset\discr L_2$ is a direct summand, and lattice $L_2$ is epistable. Assume that
$f_1\:L_1\to L_1$ is an automorphism of $L_1$ and
$(f_1)_*(K_1)=K_1$.
\roster\item
 Then $f_1$ can be extended to an automorphism $f\:L\to L$.
 \item
 Furthermore, if
$L'=L_1+_{\phi'}L_2$, $\phi'\:K_1\to K_2$, then $f_1$ can be
extended to an isomorphism $f\:L\to L'$.
\endroster
\endproposition

\proof Here (1) is a special case of (2), so we shall just prove
the latter. Let us define
$\psi_2|_{K_2}=\phi'\circ(f_1)_*|_{K_1}\circ\phi^{-1}$ and extend
it to an automorphism $\psi_2\:\discr L_2\to\discr L_2$ using the assumption
 that $K_2$ splits off as a
direct summand. Due to epistability of $L_2$, we can find an automorphism
$f_2\:L_2\to L_2$ such that $\psi_2=(f_2)_*$. Now, the assumptions of
Lemma \ref{isomorphic-gluing}
are satisfied and we obtain $f\:L\to L'$ by gluing $f_1$ and $f_2$.
\endproof

Let $L_i$ be equipped with involutions $c_i$ and
$c_i^{\discr{}}$ denote the induced involutions on $\discr L_i$.
An anti-isomorphism $\phi\:K_1\to K_2$ between
$c_i$-invariant subgroups $K_i\subset\discr L_i$ is said to be {\it
equivariant} if it commutes with $c_i^{\discr{}}|_{K_i}$.
The following lemma is also an immediate consequence of
Lemma \ref{automorphisms-gluing}

\lemma\label{gluing-involutions2}
There is an involution
$c\:L\to L$,  $L=L_1+_\phi L_2$,
$c|_{L_i}=c_i$, $i=1,2$, if
 and only if the anti-isomorphism $\phi$ is equivariant. In this case,
 such an involution $c$ is unique.
\qed\endlemma

\proposition\label{gluing-equi-automorphisms}
Assume that $L=L_1+_\phi L_2$, $\phi\:K_1\to K_2$, where
 $K_2$ splits off
as a direct summand, $\discr L_2=K_2+G_2$.
 Fix an involution $c\:L\to L$
 with $c_i=c|_{L_i}$,
and assume that
$L_2$ is $c_2$-equivariant $D$-relative epistable for some $D\subset L_2$.
 Let $f_1\:L_1\to L_1$ be a $c_1$-equivariant automorphism
 inducing $(f_1)_*\:\discr L_1\to \discr L_1$ such that
 $(f_1)_*(K_1)=K_1$. Then:
\roster\item $f_1$ can be extended to a $c$-equivariant
automorphism $f\in\Aut(L,D)$.
 \item For any given
automorphism $\phi^G_2\:G_2\to G_2$,
the automorphism
$f$ in (1) can be chosen so that its
restriction, $f_2\:L_2\to L_2$ induces $\phi^G_2$ on $G_2$.
 \item
 Given another equivariant anti-isomorphism $\phi'\:K_1\to K_2$,
 with the induced from $c_1$ and $c_2$ involution $c'$ on $L'=L_1+_{\phi'}L_2$,
 there exists an extension of $f_1$ to an isomorphism $f\:L\to
 L'$, which commutes with $c$ and $c'$,
whose restriction, $f_2\:L_2\to L_2$ induces on $G_2$ automorphism
$\phi^G_2$, and for which $f(D)=D$.
 \endroster\endproposition

\proof In the parts (1) and (2), we consider an automorphism $(f_2)_*$ on
$\discr L_2$ which is $\phi\circ(f_1)_*\circ\phi^{-1}$ on $K_2$ and $\phi_2^G$
on $G_2$. Equivariant $D$-relative epistability of $L_2$
 implies that $(f_2)_*$ can be lifted to an automorphism $f_2\:L_2\to
 L_2$ such that $f_2(D)=D$.
 The anti-isomorphism $\phi$ is equivariant, so, by Lemma
 \ref{gluing-involutions2} the required $f$ exists.

The part (3) is similar, except that we need to start with an
involution $\phi'\circ(f_1)_*\circ\phi^{-1}$ on $K_2$.
\endproof

\subsection{Realizability of  T-pairs by geometric involutions}
In this subsection we prove existence of a geometric involution $c\in C(L,\D,h)$ on a $(\D,h)$-decorated K3-lattice $L$ whose pair of eigenlattices $(T_+(c),T_-(c))$ is isomorphic to a given  T-pair $(T_1,T_2)$.
Throughout the subsection this pair
is supposed to be ascending, although the case of descending  T-pairs is analogous.
First, we show that $T_1$ and $T_2$ after an appropriate gluing give a lattice isomorphic to $T\subset L$.

\proposition\label{possibility-to-glue}
For any ascending T-pair $(T_1,T_2)$,
there exists a subgroup $K_2\subset\discr_2 T_2$
 anti-isomorphic to $\discr_2T_1$. Any such subgroup is the orthogonal complement $K_2=v^\perp$
 of some $v\in\discr_2(T_2)$, such that $\q_{T_2}(v)=-\frac12$.
\endproposition

\lemma\label{Wu-element}
Assume that $(T_1,T_2)$ is an ascending  T-pair such that $\delta_2(T_1)=0$, and $v\in\discr_2T_2$ is
the characteristic element. Then $\q_{T_2}(v)=-\frac12$.
\endlemma

\proof
Lemma \ref{sumB2} gives $\Br_2(T_1)+\Br_2(T_2)=-1$, and by additivity $\Br_2(T_2)=\Br\la v\ra+\Br(v^\perp)$,
where $\la v\ra\subset\discr_2T_2$ is spanned by $v$,
while $v^\perp$ is its orthogonal complement.
Since the quadratic forms in $\discr_2T_1$ and $v^\perp$ are even, their Brown invariants are divisible by 4,
and thus $\Br\la v\ra=-1\mod4$. This implies that $\q_{T_2}(v)=-\frac12$.
\endproof

\lemma\label{non-Wu-element} Assume that $(G,\q)$
 is an elementary enhanced 2-group,
which is not isomorphic to $n\la\frac12\ra$ for $n\le3$, and has $\delta_2(G)=1$
(i.e., $\q$ is odd). Then there exists an element $v\in G$, such that $\q(v)=-\frac12$.
 If in addition $G\ne\la-\frac12\ra$, then such $v$ can be chosen
non-characteristic.
\endlemma

\proof Existence of $v$ such that $\q(v)=-\frac12$
follows from decomposability
$G=p\la\frac12\ra+q\la-\frac12\ra$ (see Theorem \ref{2-classification})
and isomorphism $4\la\frac12\ra\cong 4\la-\frac12\ra$. If $p+q>1$, then a generator of a summand
$\la-\frac12\ra$ in such a decomposition is non-characteristic.
\endproof

\demo{Proof of Proposition \ref{possibility-to-glue} }
By Lemma \ref{anti-isomorphism-criterion}, it is sufficient to prove existence of $v\in\discr_2(T_2)$
satisfying two conditions: first,  $\q_{T_2}(v)=-\frac12$, and second, $v$  is
characteristic if $\delta_2(T_1)=0$, and non-characteristic otherwise.
 Lemma \ref{Wu-element} proves it in the case $\delta_2(T_1)=0$. If $\delta_2(T_1)=1$, then
we can use Lemma \ref{non-Wu-element} after showing that $\discr_2T_2$ cannot be isomorphic to
 $n\la\frac12\ra$, $1\le n\le3$, or to $\la-\frac12\ra$,
if $\delta_2(T_1)=1$.

Indeed, if we suppose that $\discr_2T_2=n\la\frac12\ra$, then we have $r_2(T_2)=\Br_2(T_2)=n\mod8$, and thus,
$r_2(T_1)=n-1$ and due to Lemma \ref{sumB2}, $\Br_2(T_1)=-n-1\mod8$. On the other hand, if $\delta_2(T_1)=1$, then
$\discr_2T_1$ is isomorphic to $p\la\frac12\ra+q\la-\frac12\ra$, $p,q\ge0$, and $p+q=n-1$, $\Br_2(T_1)=p-q=-n-1\mod8$,
which is impossible for $1\le n\le 3$. If we suppose that
$\discr_2T_2=\la-\frac12\ra$, then we have $r_2(T_1)=r_2(T_2)-1=0$, which contradicts to the assumption that $\delta_2(T_1)=1$.
\qed\enddemo

\proposition\label{T-eigenlattices&T-pairs}
Any  ascending T-pair $(T_1,T_2)$ is isomorphic to the pair of eigenlattices
$(T_+(c_T),T_-(c_T))$ of some involution $c_T$ in the lattice $T$.
\endproposition

\proof
By Proposition \ref{possibility-to-glue} an anti-isomorphism
$K_1=\discr_2 T_1@>\phi>> K_2\subset\discr_2 T_2$ does exist, which gives due
to Proposition \ref{gluing-involutions}, an involution, $c$, in $T_1+_\phi T_2$
with the eigenlattices $T_1$, $T_2$.
 So, it is left to verify that
 $T_1+_\phi T_2$ is isomorphic to $T=\U+\U(3)+2\A_2+\A_1\subset L$ (see Lemma \ref{T}).
 This follows from
Nikulin's stability criterion \ref{stability-criterion} applied to
$T_1+_\phi T_2$ and $T$, because both lattices have
the same inertia indices and isomorphic discriminants. For the latter, note that
$\discr_3(T_1+_\phi T_2)=\discr_3 T_1+\discr_3 T_2=\dip+3\din=\discr T$, and that $\discr_2(T_1+_\phi T_2)=\la-\frac12\ra=\discr_2 T$.
\endproof

\proposition\label{extension-T2L}
Any involution $c_T\:T\to T$ can be extended to a geometric involution $c\in C(L,\D,h)$.
\endproposition

\proof
First, let us extend a given involution $c_T\:T\to T$ to
the lattice $T'=T+_{\phi_h}\la2\ra$ glued along the anti-isomorphism $\phi_h$
between $\discr_2T=\la-\frac12\ra$ and $\la\frac12\ra=\discr\la2\ra$, so that the generator
$h\in\la2\ra$ is sent to $-h$ (see Lemma \ref{gluing-involutions2}).
 Next, we extend the involution $c_{T'}\:T'\to T'$ that we obtain to $L$ as follows.

 Let $\phi\:\discr S^0\to\discr T'$ be the anti-isomorphism defined by the sublattices $S^0,T'\subset L$,
 $L=S^0+_\phi T'$. Consider the pull-back $c_S^{\discr}\:\discr S^0\to\discr S^0$ of the involution induced by
 $c_{T'}$ in $\discr T'$ via $\phi$.
 According to Corollary \ref{bar-Aut-lifting}, involution $c_S^{\discr}$ can be lifted to a involution
 $c_S\in \overline{\Aut}(S^0,\D)$. Then, $c_{T'}$ and $c_S$ are compatible and thus yield an involution, $c$,
 in $L$ (see Lemma \ref{gluing-involutions2}) which is geometric.
\endproof

Finally, we obtain the following theorem as a direct corollary of Propositions \ref{T-eigenlattices&T-pairs} and \ref{extension-T2L}.

\theorem\label{eigenlattices&T-pairs}
Any ascending
T-pair $(T_1,T_2)$ is geometric.
\qed\endtheorem

\rk{Remark}
As was mentioned in the beginning of previous Subsection, the case of
descending T-pairs is analogous, so such pairs are also geometric.
\boxedR\endrk

\section{Deformation classes via periods}\label{Chapter_periods}

Through all this section we fix a $K3$-lattice $\L$ and
its conical $(\bold{\Delta},\bold{h})$-decoration. Note that any other
conical $(\D,h)$-decoration can be identified with the fixed decoration via
an isomorphism, which exists due to Proposition
\ref{decorated-isomorphism}.

 \subsection{The complex period map}\label{C-periods}
Recall that, according to Proposition \ref{S}, the resolution decoration of the K3-surface
$\til Y$ associated with a Zariski sextic is conical. Turning this property into a definition we extend it
to any six-cuspidal K3-surface $Z$ which can be equipped with a lattice isomorphism $\phi : H^2(\til Z)\to\L$
such that the exceptional divisors of the six cusps are mapped to $\bold{\D}$ and the preimage of $h$ belongs to the closure
of the cone generated by ample divisors (note that not any six-cuspidal K3-surface has these properties,
cf., \cite{Zar29} and \cite{Degt-Oka}).
By an {\it $\LDh$-premarked K3-surface} we mean any six-cuspidal K3-surface $Z$ equipped with a lattice isomorphism having the above properties,
the latter is called an {\it $\LDh$-premarking}.
Thus, for K3-surfaces
$\til Y$ associated with Zariski sextics, an
$\LDh$-premarking of $\til Y$ is nothing but a lattice isomorphism $f\:H_2(\til Y)\to \L$
identifying the resolution decoration with the reference conical decoration of $\L$.

Assume that $\phi\:H^2(\til Z)\cong H_2(\til Z) \to \L$ is an $\LDh$-premarking of
a six-cuspidal  K3-surface $\til Z$. Then the holomorphic $2$-forms in $\til Z$
form a line $\phi(H^{2,0}(\til Z))\subset
 \L_\C=\L\otimes\C$ called the {\it period line of $\til Z$}.
 This line is orthogonal to the sublattice $S(\bold{\D})\subset \L$ (that is the primitive closure of the sublattice generated by $\bold{\D}$ and $h$) and thus
lies in $T_\C=T\otimes\C$.
Taking the projectivization $P(T_\C)\subset P(\L_\C)$ of $T_\C\subset \L_\C$ we obtain
 a point $\om=\om(\til Z)\in P(T_\C)$
 called the
{\it period point (or simply the period) of $\til Z$}.

More generally, one can define an $\LDh$-premarking of a
holomorphic, or continuous, family of six-cuspidal $K3$-surfaces
(for example, associated with a holomorphic family of Zariski sextics)
as a locally trivial family of premarkings, and given such a family of premarkings
one obtains a well defined holomorphic, or continuous,
family of period points
in $P(T_\C)$.

According to Hodge-Riemann bilinear relations,
the period point $\om$ belongs to
the quadric $Q= \{w^2=0\}\subset P(T_\C)$ and, more precisely, to
its open subset $\widehat{\Cal D}=\{w\in Q\,|\,w\overline w>0\}$.
 This subset has two connected components,
which are exchanged by the complex conjugation (this reflects also
switching from the given complex structure on $\til Z$ to the complex-conjugate one and multiplying the marking by $-1$).
 Writing $w\in\om\subset T_\C$ as $w=u+iv$, $u,v\in
T_\R=T\otimes\R$, we can reformulate the conditions $w^2=0, w\overline w>0$ as
$u^2=v^2>0$, $uv=0$, which implies that the real plane $\la
u,v\ra\subset T_\R$ spanned by $u$ and $v$ is positive definite
and bears a natural orientation $u\wedge v$ given by $u=\Re w, v=\Im w$. Note
that the orientation determined similarly by the conjugate complex
line $\overline{\om}\subset T_\C$ is the opposite one.

The orthogonal projection of a positive definite real plane in
$T_\R$ onto another one is non-degenerate. Thus, to select one of
the two connected components of $\widehat{\Cal D}$ we fix an
orientation of positive definite real planes in $T_\R$ so that the
orthogonal projection preserves it. We call it the {\it prescribed
orientation}
and define an {\it $\LDh$-marking} as an $\LDh$-premarking for
which the orientation $u\wedge v$  of $\phi(H^{2,0}(\til Z))$ defined by the
pairs $u=\Re w, v=\Im w$ for $w\in\phi(H^{2,0}(\til Z))$ is the
prescribed one. We denote this component by $\Cal D$ and call it
the {\it period domain}. By the {\it period mapping} we understand the mapping from the set of
$\LDh$-marked K3-surfaces to $\Cal D$ respecting the above conventions.

By $\Aut^+\LDh$ we denote the group
of those automorphisms of the triple $\LDh$ that preserve the prescribed
orientation (and thus preserve $\Cal D$). The complementary coset
 $\Aut^-\LDh=\Aut\LDh\sm\Aut^+\LDh$
consists of automorphisms exchanging the connected components of
$\widehat{\Cal D}$.

Let us call an $\LDh$-marked K3-surface $Z$ {\it regular}, if
there is no $v\in \L\sm S(\D)$ such that $v^2=-2$ and $vh=v\om(\til Z)=0$.
The following statement is well known and follows, for example, from \cite{SDonat}.

\proposition\label{SD-simple-model} If $Z$ is a regular $\LDh$-marked K3-surface, then the linear system defined by
$h$ induces a degree $2$ map $Z\to P^2$ and this map is a double covering of $P^2$ branched over a Zariski sextic.
Moreover, this correspondence establishes an isomorphism between
the space of regular $\LDh$-marked K3-surfaces $Z$
and the space of $\LDh$-marked K3-surfaces $\til Y$ associated with Zariski sextics.
\qed
\endproposition

\remark{Remark} By Riemann-Roch inequality, the homology
classes $h-\sum\frac{e''_i+2e'_i}3$ and $h-\sum\frac{e'_i+2e''_i}3$ (cf., \ref{Q-formula})
are realized in $Z$ by effective divisors, and the conic through the cusps
making the ramification sextic to  be a Zariski sextic is nothing but the image of each of these divisors.
\boxedR\endremark

Part (1) of the following theorem is a standard consequence of
surjectivity of the period mapping (see \cite{KulSurj}) and
Proposition \ref{SD-simple-model}; part (2) follows from
the strong Torelli theorem (see \cite{BR})
(a systematic study of plane sextics with arbitrary homological conditions on collections
of simple singularities via such an approach is undertaken
in \cite{Degt-def}).

\theorem\label{surjectivity}
Assume that $\om\subset \L_\C$ is a line orthogonal to $h$ and $\D$. Then:
 \roster
 \item
 A point $\om\in\Cal D$ is the period of a surface $\til Y$ associated with some Zariski sextic $A\subset P^2$
 and some $\LDh$-marking $\phi\:H^2(\til Y)\to \L$
 if and only if there is no $v\in \L\sm S(\D)$ such that $v^2=-2$ and $vh=v\om=0$.
 \item
If we are given another Zariski sextic $A'$ with an
$\LDh$-marking $\phi'\:H^2(\til Y')\to \L$ having the same period $\om$,
 then $A$ with $A'$ are projectively equivalent, and in particular,
the projections of $Y$ and $Y'$ to $P^2$ are projectively equivalent.
\quad\qed
\endroster
 \endtheorem

By this theorem, the image of the period mapping is the complement in $\Cal D$
of a certain arrangement $\Cal H$ of hyperplane sections.
As we check below, this arrangement splits naturally in three
sub-arrangements corresponding to three different kinds of codimension one degenerations:
appearance of a node, gluing of two $\A_2$-singularities to an $\A_5$-singularity,
 and degeneration to a double ruled surface.
In arithmetical terms, such a splitting of the
arrangement $\Cal H$ is given by distinguishing three sets of $(-2)$-roots:
$V_2^n=\{v\in T\,|\, v^2=-2, v\ne h\mod 2\L\}$,
$V_2^h=\{v\in T \,|\, v^2=-2, v=h\mod 2\L\}$, and
$V_2^g=\{v\in \L\,|\, v^2=-2, vh=0, vS^0\ne 0, v\ne h\mod 2\L\}$.
We denote by $H_v$ the hyperplane section of $\Cal D$ defined by $xv=0$
and, finally, define
$\Cal H$ to be the union of three arrangements, $\Cal H=\Cal H^n\cup \Cal H^h\cup \Cal H^g$,  where $\Cal H^n=\bigcup_{v\in V_2^n}H_v$, $\Cal
H^h=\bigcup_{v\in V_2^h}H_v$, and $\Cal H^g=\bigcup_{v\in V_2^g}H_v$.

\rk{Remark}
Note that the hyperplane $xv=0$ does not intersect $\widehat{\Cal D}$ if $v\in \L, v^2=-2, vh=0, vS^0\ne 0,$ and $ v=h\mod 2\L$. Indeed,
in such a case the orthogonal complement to the plane generated by $v$ and $e\in S^0$ with $e^2=-2, ev\ge 2$
has the positive inertia index
$<3$, while if the orthogonal complement contains a point $\omega\in \widehat{\Cal D}$ then it contains a positive definite
$3$-plane generated by $h$, $\Re\omega$, $\Im\omega$.
By this reason we do not need to consider the set $\{v\in \L\,|\, v^2=-2, vh=0, vS^0\ne 0, v=h\mod 2\L\}$,
and even could replace $V_2^g$ by $\{v\in \L\,|\, v^2=-2, vh=0, vS^0\ne 0\}$.
\boxedR\endrk

In this notation, Theorem \ref{surjectivity} can be rephrased as
follows: {\it the set of periods of the K3-surfaces of Zariski
sextics is the complement $\Cal D\sm \Cal H$.}
Furthermore, together with the strong Torelli theorem
it implies the following
statement, which is also well known.

\theorem\label{fine-moduli}
The space $\Cal D\sm\Cal H$ is a fine moduli space of regular
$\LDh$-marked K3-surfaces. Its quotient by $\Aut^+\LDh$ is
naturally identified with the space of projective classes of Zariski sextics.
\qed\endtheorem

{\it Codimension 1 degenerations} of regular $\LDh$-marked K3-surfaces are represented by
non-singular points of $\Cal H$. Speaking more formally, such a degeneration is represented by
a holomorphic disc, $f\:D^2\to\Cal D$, intersecting $\Cal H$ transversally at a single point, say, $f(0)\in\Cal H$.
This gives rise to a holomorphic family $\til{Y}_t$, $t\in D^2$, of K3-surfaces presented as double planes
$\til{Y}_t\to P^2$ ramified along Zariski sextics $A_t$, experiencing a certain degeneration at $t=0$.

\proposition\label{three kinds}
Consider a codimension 1 degeneration of regular $\LDh$-marked K3-surfaces $\til{Y}_t$, $t\in D^2$,
with the periods $f(t)\in\Cal D$ degenerating to $f(0)\in H_v$, for some $(-2)$-root $v$.
Then the type of the degeneration of $\til{Y}_t$ and
of the corresponding Zariski sextics, $A_t$, depends on $v$ as follows.
\roster\item
If $v\in V_2^n$, then $\til Y_t$ as well as
$A_t$ experiences a nodal degeneration.
And any nodal degeneration of $A_t$ and $\til Y_t$ can be represented via such a family.
\item If $v\in V_2^g$, then a pair of cusps of $A_t$ as well as a pair of cusps of $\til Y_t$, experiences
a degeneration to $\A_5$-singularity.
And any codimension one degeneration in which at least one of the $6\A_2$ experiences a
degeneration to a deeper singularity can be
 represented via such a family.
\item If $v\in V_2^h$, then $A_t$
degenerates into a triple conic and $\til{Y_t}$ degenerates into an
elliptic K3. Moreover,
this elliptic K3 is a double of the ruled surface $\Sigma_4=P_{P^1}(\Cal O_{P^1}(4)\oplus \Cal O_{P^1})$
branched along a union of the $(-4)$-section and a 6-cuspidal trigonal curve
in $\Sigma_4$ where all the six cusps are coplanar, that is belong to a section of $\Sigma_4$ disjoint from the $(-4)$-section.
And any degeneration of regular $\LDh$-marked K3-surfaces
to such a double $\Sigma_4$ and of Zariski sextics to a triple conic
can be represented via such a family.
 \endroster
 \endproposition

\proof
According to Saint-Donat's results on the
 projective models of K3-surfaces, see \cite{SDonat}, those surfaces that carry
a numerically effective divisor of degree $2$ generating the
linear system without fixed components are the double covers of the plane branched in sextics with simple singularities, and
those that carry a numerically effective divisor of degree $2$ generating the linear system with a fixed component are the double covers of $\Sigma_4$ branched along a curve with simple singularities, where the curve is linear equivalent to $2c_1(\Sigma_4)$
and splits in a union of the $(-4)$-section
$E_0$ and a trigonal curve $3E_0+12F$ disjoint from $E_0$ (here, $F$ stands for the fiber of $\Sigma_4$; note that
$c_1(\Sigma_4)= 2E_0+6F$). The coplanar condition in
statement (3) reflects (and is equivalent to) the existence of a
6-cuspidal divisor sublattice and a master element in its extension; indeed, the proof of the sufficiency of coplanarity is literally the same as the proof of Lemma \ref{Q-formula}, while the necessity follows then from transitivity of $\Aut^+\LDh$ action on
the set of $V_2^h$-hyperplanes.
Thus, there remain to examine codimension one degenerations of Zariski sextic to a sextic with simple singularities.
As is well known and, for example,  can be easily deduced from \cite{Pho},
the only codimension one
degenerations of Zariski sextics to a sextic with simple singularities are nodal degenerations
and $\A_5$-degenerations of a pair of cusps, $2\A_2\subset \A_5$.
\endproof

Arithmetically, each $\A_5$-degeneration of a pair of cusps
results in a primitive embedding of $2\A_2$ into $\A_5$.
When $2\A_2$ is primitively embedded into $\A_{5}$, then, since $\vert \A_2\vert=3$ and $\vert \A_5\vert=6$,
the orthogonal complement of $2\A_2$ in $\A_5$ is $\la- 6\ra$. Thus, if $v\in V_2^g$ then a natural $(-6)$-vector appears:
it is nothing
but the orthogonal projection of $3v$ on $T$. We denote the set of all these $(-6)$-vectors by $V_6$.

\proposition\label{V6} The set $V_6$ coincides with $\{u\in T\,|\,u^2=-6,uT\subset 3\Z\}$.
\endproposition

\proof Necessity follows immediately from the definition of $V_6$.
So, let us assume that $u\in T, u^2=-6$, and $uT\subset 3\Z$.
Then, $[\frac{u}3]$ is an element of square $-\frac23$ in the 3-primary discriminant component
 $\discr_3(T)$. Thus, $[\frac{u}3]$ is glued to an element of square $\frac23$ in $\discr_3(S_0)$.
By Lemma \ref{elements-description},
there exist 30 such elements $[b_j^i]=[\frac{e_i'+2e_i''}3+\frac{2e_j'+e_j''}3]$, $1\le i, j\le 6$, $i\ne j$,
$[b_j^i]=-[b_i^j]$.
Gluing of $u$ with, say, $[b_i^j]$ produces
a  $(-2)$-root $e=\frac13(u-(e_i'+2e_i''+2e_j'+e_j''))\in L$
that generates $\A_5\subset L$ together with $e_i',e_i'',e_j',e_j''$; namely, $e$ represents
the ``middle'' root of $\A_5$, or equivalently, $u=e_i'+2e_i''+3e+2e_j'+e_j''$.
Note, that $eS^0\ne 0\mod 2\Z$ and therefore $e$ belongs to
$V_2^g$.
\endproof

To each $v\in V_2^n\cup V_2^h$ we associate the $(-2)$-reflection $\rho^T_v\:T\to T, x\mapsto x+\langle x,v\rangle v$,
and to each $u\in V_6$ the $(-6)$-reflection $\rho^T_u\:T\to T, x\mapsto x +\frac13 \langle x,u\rangle u$.
The first one is evidently the restriction to $T$ of the $(-2)$-reflection $\rho_v\: \L\to \L$, while the second
one is the restriction to $T$ of the reflection $\rho^\L_u\: \L\to \L$
in the $3$-dimensional mirror generated by the triple $e, e_i'+e_j'', e_j'+e_i''$
like in the proof of Proposition \ref{V6}.
Note also that for each $u\in V_6$, the hyperplane section $H_u$
of $\Cal D$ defined by $xu=0$ coincides with $H_e$ we used in our definition of $\Cal H^g$.

\subsection{The real period map}\label{real-period}
Let us fix a geometric involution $c\in C\LDh$ (see Section \ref{geometric-involutions}), extend
$c$ to a complex linear involution on $\L \otimes \C$ and denote
also by $c$ the induced involutions on $T_\C$, $P=P(T_\C)$, and
$\widehat{\Cal D}$.
An $\LDh$-marking $\phi$ of a real K3-surface $(Z,\conj)$ is
called {\it $c$-real} if $c\circ\phi=\phi\circ\conj$.
If a real K3-surface
$(Z,\conj)$ can be equipped with  a $c$-real $\LDh$-marking, we say that this surface is of {\it homological type $c$}.
A real Zariski sextic is said to be of homological type $c$, if the K3-surface $\til Y$ associated to
it is equipped with a real structure $\conj$
lifted from the real structure of $P^2$ and $(\til Y,\conj)$ is equipped with a $c$-real $\LDh$-marking.

The following Lemma shows that this definition of homological type concords with the one
given in Section \ref{geometric-involutions}.

\lemma\label{homological types} If $(Z,\conj)$ is of homological type $c$, and $c'$ represents the same element
in $C[\L,\D,h]$ as $c$, then $(Z,\conj)$ is of homological type $c'$ as well.
\endlemma

\proof If $c=fc'f^{-1}$ with $f\in\Aut^-\LDh$, then $c=gc'g^{-1}$ where $g=cf\in\Aut^+\LDh$.
\endproof

The involution $c$ permutes the two components,
$\Cal D$ and $\overline{\Cal D}$, of $\widehat{\Cal D}$,
and thus $\overline c(\Cal D)=\Cal D$, where $\overline c\:T_\C \to T_\C$ is
the composition of $c$ with the complex conjugation in $T_\C$.
Let $\widehat{\Cal D}_\R^{c}$ and $\Cal D_\R^{c}$ denote the fixed
point sets of $\overline c$ restricted to $\widehat{\Cal D}$ and
$\Cal D$. The latter fixed point set consists of the lines generated by
$w=u_++iu_-$ such that $u_\pm\in T_\pm\otimes\R$,
$u_+^2=u_-^2>0$, and the orientation $u_+\wedge u_-$ is the prescribed
one. Since $c$ is geometric, both ${\Cal D}_\R^{c}$ and its
(trivial) double covering $\widehat{\Cal D_\R^{c}}$ are nonempty.

As it follows from definitions, the period of a c-real
$\LDh$-marked K3-surface
belongs to $\Cal D_\R^{c}=\{x\in\Cal D\,|\,c(x)=\overline x\}$. Therefore, we call $\Cal D_\R^{c}$ the
{\it real period domain}. It splits in a direct product,
$$
\Cal D_\R^{c}=\Cal D(T_+)\times \Cal D(T_-),
$$
where $\Cal D(T_+)=\Cal D\cap P(T_+\otimes\R)$ and $\Cal D(T_-)=\Cal D\cap P(T_-\otimes\R)$ are the real hyperbolic (Lobachesvki) spaces associated with the (hyperbolic) lattices $T_\pm$.
In other words, one can fix a half of the cone $u_+^2>0$ in $T_+\otimes \R$ and a half of the cone $u_-^2>0$ in
$T_-\otimes \R$ in a way that the prescribed halves respect the prescribed orientation (that is to make the orientation
$u_+\wedge u_-$ to be the prescribed orientation for any $u_+,u_-$ from the fixed half-cones) and then $\Cal D(T_\pm)$
become the spaces of the vector half-lines in the chosen two half-cones.

\theorem\label{real-surjectivity} The periods of $c$-real  regular $\LDh$-marked
$K3$-surfaces form in $\Cal D_\R^{c}$
the complement of $\Cal D_\R^{c}\cap \Cal H$.
The space $\Cal D_\R^{c}\sm\Cal H$
is a fine moduli space of $c$-real regular $\LDh$-marked
K3-surfaces.
\endtheorem

\proof This is a straightforward consequence of Theorems \ref{surjectivity} and \ref{fine-moduli}.
\endproof

Let us put $V_2^n(T_\pm)=V_2^n\cap T_\pm, V_2^h(T_\pm)=V_2^h\cap
T_\pm$, $V_6(T_\pm)=V_6\cap T_\pm$ and define
$\Cal H^n(T_\pm)=\bigcup_{v\in V_2^n(T_\pm)}H_v\cap \Cal D(T_\pm)$, $\Cal
H^h(T_\pm)=\bigcup_{v\in V_2^h(T_\pm)} H_v\cap \Cal D(T_\pm)$,
$\Cal H^g(T_\pm)=\bigcup_{u\in V_6(T_\pm)}H_u\cap \Cal D(T_\pm)$.
We denote by $\Ch_\pm$ the set of connected components of the complement
 $\Cal D(T_\pm)\sm\Cal H_\pm$ of the hyperplane
arrangement
 $\Cal H_\pm=\Cal H^n(T_\pm)\cup \Cal
H^h(T_\pm)\cup \Cal H^g(T_\pm)$.
Every of these (infinite in number) components
is obviously a convex polyhedron.

\lemma\label{up-to-codimension-two}
Every connected component of $\Cal D_\R^{c}\sm\Cal H$ is obtained from a product $P_+\times P_-$ of two
polyhedra $P_\pm\in\Ch_\pm$ by removing a codimension two subset.
In particular, the inclusion map identifies the set of
connected components of $\Cal D_\R^{c}\sm\Cal H$
with $\Ch_+\times\Ch_-$.
\endlemma

\proof The difference is due to $(-2)$- and $(-6)$-vectors $v\in V_2^n\cup V_2^h\cup V_6$,
$v\notin T_\pm$, that
have nonzero components
$v_\pm\in T_\pm\otimes\R$
and hence $H_v\cap\Cal D_\R^{c}=v_+^\perp\times v_-^\perp$ (orthogonal complement $v_\pm^\perp$ is in $\Cal D(T_\pm)$)
is of codimension $2$.
\endproof

By $\Aut^+\LDhc$ we denote the subgroup of $\Aut^+\LDh$ consisting of those elements of the latter
that commute with $c$. The important
examples of such elements are $\rho_v$ with $v\in V_2^n(T_\epsilon)\cup V_2^g(T_\epsilon)$ and
$\rho^\L_u$ with $u\in V_6(T_\epsilon), \epsilon=\pm$. Note that all these elements act as a reflection in $T_\epsilon$ and as the identity in $T_{-\epsilon}$.

Each element of $\Aut^+\LDhc$ preserves the eigenspaces $T_\pm$ and hence, group $\Aut^+\LDhc$
naturally acts on $\Ch_+, \Ch_-$, and $\Ch_+\times\Ch_-$.

\theorem\label{K3-Zariski-deformation} The set of deformation classes of real Zariski sextics of homological type $c$
is in a natural bijection with the set of orbits of the action of $\Aut^+\LDhc$ on $\Ch_+\times\Ch_-$.
\endtheorem

\proof
Given a real Zariski sextic of homological type $c$, we replace it by the associated K3-surface $\til Y$ and equip  $\til Y$ with a real structure
$\conj$ (lifted from the real structure of $P^2$) and with a $c$-real $\LDh$-marking. Due to Proposition \ref{SD-simple-model},
the sextic can be reconstructed back, up to projective transformation, as the branch curve of the linear system given by $\h$. Furthermore,
if we start from a regular $c$-real $\LDh$-marked K3-surface $Z$ and apply this reconstruction procedure, the real structure we get
on $P^2$ is the standard one, since all real structures on $P^2$ are isomorphic to each other.
Now, it remains
 to notice that the $c$-real $\LDh$-markings of $\til Y$ form an orbit of the action of $\Aut^+\LDhc$
as it follows from
Theorem \ref{real-surjectivity}, Lemma \ref{up-to-codimension-two}, and
the strong Torelli theorem.
\endproof

\subsection{Deformation classification via geometric involutions}

\theorem\label{deformation-classification} For any homological type $c\in C\LDh$ there is
one and only one deformation class of real Zariski sextics of homological type $c$.
In particular, the deformation classes of real Zariski sextics are in one-to-one correspondence with the set of conjugacy classes of ascending
geometric involutions.
\endtheorem

Let us fix a geometric involution $c\in C\LDh$.

\proposition\label{transitivity} The action of $\Aut^+\LDhc$
on $\Ch_+\times\Ch_-$ is transitive.
\endproposition

Let us analyze how the reflections $\rho^T_v$ for $v\in V(T_\pm)$ can be extended from $T_\pm$
to the whole $\L$. The extension in case of $v\in V_2(T_\pm)$ is obvious.

\lemma\label{2-reflections} For any $v\in V_2(T_\pm)$, the
reflection $\rho_v(x)=x-\frac{xv}{v^2}v$ is well-defined on the
whole lattice $\L$, it acts
as the identity in $S$ as well as in $T_\mp$ (opposite to $T_\pm$), and
$\rho_v\in\Aut^+\LDhc$. \qed\endlemma

The case of $v\in V_6(T_\pm)$ is a bit more subtle.

\lemma\label{6-reflections} For any $v\in
V_6(T_\pm)$, the reflection
$\rho^T_v\vert_{T_\pm} \:T_\pm\to T_\pm$ can be
extended to an automorphism $f\in\Aut^+\LDhc$ that
acts on $T_\mp$ as
$\id_{T_\mp}$.
\endlemma

\proof
Note that $\rho^T_v\vert_{T_\pm}$ induces the identity map in $\discr_2 T_\pm$ and, thus,
 by Proposition \ref{gluing-involutions}
can be extended to the whole $T$ by gluing with the identity automorphism on $T_\mp$.
 Applying Proposition \ref{gluing-involutions} once more, we extend
the result of gluing to an involution $\rho_v^{T'}\:T'\to T'$ that maps $\h$ to $\h$.
Finally, we extend $\rho_v^{T'}$ to $f\:\L\to \L$ applying Proposition \ref{gluing-equi-automorphisms},
where we use the $c$-equivariant $\D$-relative epistability of $S^0$
(see Corollary \ref{S-equi-epistability}).
To check that $f\in\Aut^+\LDhc$ it is left to notice that
$\rho^T_v\vert_{T_\pm}$ preserves the halves of the
 cones $\{v\in T_\pm\,|\,v^2>0\}$, and
 therefore $f$ preserves the prescribed orientation of the
positive definite real planes in $T_\R$.
\endproof

\demo{Proof of Proposition \ref{transitivity}}
The reflection group generated by $\rho^T_v$,
$v\in V_2^n(T_\pm)\cup V_2^g(T_\pm)\cup V_6(T_\pm)$,
acts transitively on $\Ch_\pm$.
Thus, Lemmas \ref{2-reflections} and \ref{6-reflections} imply Proposition \ref{transitivity}. \qed
\enddemo

\demo{Proof of Theorem \ref{deformation-classification}} The "only if" part is trivial. The other part follows from Theorem \ref{K3-Zariski-deformation}
and Proposition \ref{transitivity}.
\enddemo

\subsection{From a trigonal curve in $\Sigma_4$ to Zariski sextics in $P^2$}\label{explicit-partners}
With a given geometric ascending involution $c\in C^<\LDh$ and an element $v\in V_2^h(T_-)$
we associate its {\it partner involution $c^v$}
that acts in  $T$ as a composition of $c$ with reflection against the plane generated by $\h$ and $v$, and in $S^0$
as a composition of $c$ with the involution permuting the elements of $\D$ so that
the generators in each of the six $\A_2\subset S^0$ are transposed.
 Then $c^v\in C^<\LDh$ as well, since by the definition $T_+(c^v)=T_-(c)\cap v^\perp$
and $T_-(c^v)\cap v^\perp=T_+(c)$.

We can peak a c-real $\LDh$-marked K3-surface $Z$, whose period point is a generic point of $H_v$.
Then according to Proposition \ref{three kinds}(3), it is a double covering of $\Sigma_4$
branched along the union of the $(-4)$-section and a $6$-cuspidal trigonal curve with coplanar cusps.
 In addition to the real structure in $Z$ inducing involution $c$ in its lattice, there is another one
inducing $c^v$; the two real structures in $Z$ are the two liftings of a real structure in $\Sigma_4$ and differ by
the deck transformation of the double covering $Z\to\Sigma_4$.
 By real perturbations of $(Z,c)$ and $(Z,c^v)$ we obtain regular $\LDh$-marked
  K3-surfaces and the deformation types of
 the corresponding real Zariski sextics do not depend on the choice of perturbations. Hence, we may speak on the
 well defined partner $Z$-deformation
classes of real Zariski sextics. Below, we make some explicit choice of perturbations, which allows
 us to compare the partner real Zariski sextics obtained in this way.
Namely, we prove that such a pair of real Zariski sextics is placed in $\R P^2$ in a trigonal
reverse position (see Section \ref{reversion} for definitions
of reverse position, trigonality, and reversion partners).

\proposition\label{partner_topology} Real Zariski curves in a pair of partner
$Z$-deformation classes are reversion partners.
\endproposition

\proof
Consider the weighted projective space $P(2,1,1,1)$ with coordinates $q,x,y,z$ of weights $2,1,1,1$ and the conic $B$ defined by equations
$q=0,\, Q=0$ where $Q$ denotes the polynomial $xy-z^2$. Note that the cone over $B$ with vertex at $(1,0,0,0)$ is the ruled surface $\Sigma'_4$ (the $(-4)$ section of $\Sigma_4$
is contracted here to the vertex of the cone). Note also that due to the isomorphism
$P^1\to B$ given by $x=u^2, y=v^2,z=uv$ any homogeneous degree $2n$ polynomial $g_{2n}(u,v)$ in variables $u,v$ can be rewritten
as a homogeneous degree $n$ polynomial $f_n(x,y,z)$ in $x,y,z$ (such a presentation is unique modulo the ideal generated by $Q$).

Now, pick a 6-cuspidal trigonal curve $C$ in $\Sigma'_4$ so that all the six cusps belong to the plane $q=0$
and write a defining polynomial of $C$
as $q^3+g_4(u,v)q^2+g_8(u,v)q+g_{12}(u,v)$.
Then, express the polynomials
$g_4$, $g_8$, $g_{12}$ in variables $u$, $v$ as polynomials $f_2$, $f_4$, $f_6$
in variables $x,y,z$. This allows us to define $C\subset\Sigma'_4$ by equations
$$
q^3+f_2(x,y,z)q^2+f_4(x,y,z)q+f_{6}(x,y,z)=0,\, Q=0,
$$
and then include it in a family
$$
q^3+f_2(x,y,z)q^2+f_4(x,y,z)q+f_{6}(x,y,z)=0,\, Q=tq.
$$
Due to the special position of the cusps $q^3+g_4(u,v)q^2+g_8(u,v)q+g_{12}(u,v)=q^3+(g_2(u,v)q+g_6(u,v))^2$,
 so that $q^3+f_2(x,y,z)q^2+f_4(x,y,z)q+f_{6}(x,y,z)=q^3+(f_1(x,y,z)q+f_3(x,y,z))^2$.

As a straightforward calculation shows, the associated family of plane curves defined by
$$
Q^3+t(f_1(x,y,z) Q+ t f_3(x,y,z))^2=0
$$
is a family of Zariski sextics (with 6 cusps on $B$) degenerating to
$Q^3=0$. The double coverings of $P^2$ branched along these sextics are naturally embedded into
the weighted projective space $P(3,2,1,1,1)$ with coordinates $r,q,x,y,z$ of weights $3,2,1,1,1$
where they become defined by equations
$$
r^2+(q^3+(f_1(x,y,z)q+f_3(x,y,z))^2)=0,\, Q=tq,
$$
and thus form a family converging to the double covering of $\Sigma_4$ branched along $C$:
$$
r^2+(q^3+(f_1(x,y,z)q+f_3(x,y,z))^2)=0,\, Q=0.
$$
To get a pair of real Zariski curves in a pair of $Z$-deformation classes we need to select a pair of M\"obius real structures converging
to $(Z,c)$ and $(Z,c^v)$. In other words to select appropriate real forms of the double coverings, namely, to choose an appropriate sign of $t$ and an appropriate sign of $r^2$. Such a pair of real forms is given, with respect to the standard real structure on $P(3,2,1,1,1)$, by
$$
r^2+(q^3+(f_1(x,y,z)q+f_3(x,y,z))^2)=0,\, Q=tq,
$$
for $t>0$ and by
$$
r^2-(q^3+(f_1(x,y,z)q+f_3(x,y,z))^2)=0,\, Q=tq,
$$
for $t<0$. The real Zariski curves $C'_t, t>0$, of the first family converge to $C$, while the real Zariski curves $C''_t, t<0$, of the second one
converge to the image $C'$ of $C$ under the reflection $q\mapsto -q$ in the cylindrical part of $\Sigma'_4$, and thus
the trigonality and reversion position property from the definition of reversion partners follows from
the possibility to trivialize the real part of the surface family over the real locus of the cylindrical part of $\Sigma'_4$.
The signs $o(C)$ and $o(C')$ are opposite due to  Corollary \ref{cusps-to-halfs}, while Proposition \ref{type-via-lattice}
implies that  $C$ and $C'$ are both of the same type.
\endproof

\section{Arithmetics of the  T-pairs}\label{mainProof}

\subsection{Geography of the  ascending T-pairs}
 We start analysis of the  T-halves, which we denote here by $M$,
 with reviewing their possible
 numerical invariants. On the first step, we ignore its 3-primary component
 and look only at the rank $r(M)$, the discriminant 2-rank $r_2(M)$, and the discriminant parity $\delta_2(M)$.
 Namely, we prove the following

\lemma\label{pairs-of-2ranks} Assume that $(M,M')$ is an
 ascending T-pair. Then
the combination of the invariants $r$ and
$r_2$ of $M$ is one of the fifteen ones listed in
the Table 3A, and
the combination of the invariants $r'=9-r$, $r_2'=r_2+1$
of $M'$ appears at the corresponding
position in Table 3B.
 The value of $\delta_2$ for $M$  indicated in each row of Table 3A
is the only possible value for all the pairs $(r,r_2)$ in this row
 (for the last row, the both values, $0$ and $1$, of $\delta_2$ are allowable).
The value $\delta_2'$ for $M'$ is $1$ for each pair $(r',r_2')$.
\endlemma
\midinsert
\hskip15mm$\matrix\text{Table 3A. T-halves $M$}\\ \\
\boxed{\matrix
 (r,r_2) & \delta_2\\
 \text{---------------------------}&\text{---}\\
 (2,0), (4,0), (6,0), (8,0) &0\\
 (1,1), (3,1),(5,1),(7,1) &1\\
 (3,3), (5,3) &1\\
 (2,2), (4,2), (6,2), (4,4) &0,1\\
  \endmatrix}\endmatrix$
 \hskip-3mm
$\matrix\text{Table 3B. T-halves $M'$}\\ \\
\boxed{\matrix (r',r_2') & \delta_2'\\
 \text{---------------------------}&\text{---}\\
 (7,1), (5,1),(3,1),(1,1)& 1\\
 (8,2), (6,2), (4,2), (2,2) &1\\
 (6,4), (4,4)
   &1\\
 (7,3), (5,3), (3,3), (5,5) &1\\
 \endmatrix}\endmatrix$
\endinsert
\proof
 For  T-pairs the relations $1\le r,r'\le8$, $r+r'=9$ hold by definition, the ascending condition means $r_2'=r_2+1$,
 and Lemma \ref{r2-estimate} yields $r_2\le\min(r,r'-1)$. Together with the congruence
 $r_2=r\mod2$ (see Lemma \ref{r2=r-mod2}), these restrictions forbid all the cases except the fourteen ones listed in the Tables 3A-B.

 The value $\delta_2'=1$ follows from the definition of  ascending T-pairs in \ref{T-pairs-T-halves} (property (4) in Lemma
 \ref{geometric-properties}). For the values of $\delta_2$ in Table 3A note that for $r_2=0$ we have $\discr_2M=0$, and thus $\delta_2=0$. On the other hand, for odd values of $r_2$, we  have $\delta_2=1$,
 see Theorem \ref{2-classification}(1).
 \endproof

Now, for each combination of values $r,r_2,\delta_2$ from Tables 3A-B we will finalize the enumeration of the IDs by giving
possible values of $(p,q)$ describing the discriminant 3-component, $\discr_3 M=p\la\frac23\ra+q\la-\frac23\ra$.
 Namely, we prove the following.
\proposition\label{preliminary-list}
Assume that $(M,M')$ is an ascending  $T$-pair.
Then the combination of the invariants $r$,
$r_2$, $\delta_2$, $p$, $q$ for $M$ and the corresponding
invariants
 $r'=9-r$, $r_2'=r_2+1$, $p'=1-p$, $q'=3-q$ for $M'$ (recall that $\delta_2'=1$)
is contained in one or the rows of Table 4 (the matching values of $(p,q)$ and
 $(p',q')$ appear in the corresponding positions of that row).
\endproposition
{\eightpoint
\midinsert
\centerline{
$\matrix\text{\tenrm Table 4.}
\\ \\
\boxed{\matrix
\delta_2&(r,r_2)&(p,q)&(r',r_2')&(p',q')\\
\text{-----}&\text{-----}&\text{---------------}&\text{------}&\text{---------------}\\
0&(2,0)&(0,0),(1,1)&(7,1)&(1,3),(0,2)\\
0&(4,0)&(0,1),(1,2)&(5,1)&(1,2),(0,1)\\
0&(6,0)&(0,2),(1,3)&(3,1)&(1,1),(0,0)\\
0&(8,0)&(0,3)&(1,1)&(1,0)\\
1&(1,1)&(0,0),(1,0)&(8,2)&(1,3),(0,3)\\
1&(3,1)&(0,0),(0,1),(1,1),(1,2)&(6,2)&(1,3),(1,2),(0,2),(0,1)\\
1&(5,1)&(0,1),(0,2),(1,2),(1,3)&(4,2)&(1,2),(1,1),(0,1),(0,0)\\
1&(7,1)&(0,2),(0,3),(1,3)&(2,2)&(1,1),(1,0),(0,0)\\
0,1&(2,2)&(0,0),(1,1)&(7,3)&(1,3),(0,2)\\
1&(2,2)&(0,1),(1,0)&(7,3)&(1,2),(0,3)\\
0&(4,2)&(0,3),(1,0) &(5,3)&(1,0),(0,3)\\
0,1&(4,2)&(0,1),(1,2) &(5,3)&(1,2),(0,1) \\
1&(4,2)&(0,0),(0,2),(1,1),(1,3) &(5,3)&(1,3),(1,1),(0,2),(0,0)  \\
0&(6,2)&(1,1)&(3,3)&(0,2)\\
0,1&(6,2)&(0,2),(1,3)&(3,3)&(1,1),(0,0)\\
1&(6,2)&(0,1),(0,3),(1,2)&(3,3)&(1,2),(1,0),(0,1)\\
1&(3,3)&(0,0),(0,1),(0,2),(1,0),(1,1),(1,2)&(6,4)&(1,3),(1,2),(1,1),(0,3),(0,2),(0,1)\\
1&(5,3)&\text{\rm any } p\le1,q\le3&(4,4)&p'=1-p,q'=3-q\\
0,1&(4,4)&(0,3),(1,0)&(5,5)&(1,0),(0,3)\\
1&(4,4)&(0,0),(0,1),(0,2),(1,1),(1,2),(1,3)&(5,5)&(1,3),(1,2),(1,1),(0,2),(0,1),(0,0)\\
\endmatrix}
\endmatrix$}
\endinsert}
The proof is based on the following relations.

\proposition\label{generalized-restrictions}
Assume that lattice $L$ is even, the discriminant $p$-ranks $r_p$ vanish for all $p\ne2,3$, and the primary components $\discr_p(L)$ for $p=2$ and $p=3$ are elementary. Let $\sigma\in\Z$ denote the signature of $L$, $\delta_2\in\{0,1\}$ the parity of
$\discr_2L$, and $\Br_3$ the Brown invariant of $\discr_3L$.
Then
 \roster
 \item $r_2=0$ implies that $\Br_3=\sigma\mod8$;
\item $r_2=1$ implies that $|\Br_3-\sigma|=1\mod8$;
\item  $r_2=2$ together with
$|\Br_3-\sigma|=4\mod8$ imply that $\delta_2=0$;
\item $\delta_2=0$ implies that
$r_2$ is even and $\Br_3=\sigma\mod4$.
\item
$r_2=r$ together with $\delta_2=0$
imply that $\Br_3=-\sigma\mod8$;
\item
$r_3=r$ implies that $\Br_3=2\sigma\mod8$.
\endroster
\endproposition

\proof
 The Brown invariant $\Br_2$ of $\discr_2L$ can be expressed as the modulo 8 residue of
 $\sigma-\Br_3$ by the van der Blij theorem (see
Corollary \ref{vanderblij-refined}). Since $\Br_2=0$ for $r_2=0$, this implies (1).
If $r_2=1$, then $\Br_2=\pm1$ (see Example \ref{2-rank-1}), which implies (2).
The case of $r_2=2$ with $\delta_2=0$ is considered in the Example \ref{2-rank-2}, it yields (3).
Theorem \ref{2-classification}(1) yields (4).

If $r_2=r$, then $L$ is divisible by $2$ (see Proposition \ref{p-divisible}), moreover
the condition $\delta_2=0$ implies that division by $2$ yields an even lattice
$L'=L(\frac12)$.
 Moreover, parts (2) and (3) of Proposition \ref{p-divisible} imply that
  $\discr L'$ is anti-isomorphic to $\discr_3 L$, and thus, $\Br(L')=-\Br_3(L)$.
  Applying the van der Blij theorem to $L'$, we obtain $\Br(L')=\s\mod8$, and thus
$\Br_3(L)=-\s$, and applying it to $L$, we get $\Br(L)=\Br_2(L)+\Br_3(L)=\s\mod8$, and thus,
$\Br_2(L)=2\s\mod8$, which yields (5).

If $r_3=r$, then $L$ is divisible by $3$ and $L'=L(\frac13)$ is an
even lattice whose discriminant $\discr L'=\discr_2L'$
(again, due to  Proposition \ref{p-divisible})
is anti-isomorphic to $\discr_2L$. Like in the previous case, applying
van der Blij's theorem to $L'$ and $L$
we obtain $\Br(L')=-\Br_2(L)=\s\mod8$, and $\Br(L)=\Br_2(L)+\Br_3(L)=\s\mod8$, which yields (6).
\endproof

\corollary\label{T-half-restrictions} Assume that $M$ is an
 $T$-half of rank $r$, $2$-rank $r_2$, and
$\discr_3M=p\dip+q\din$.
Then \roster
\item if $r_2=0$, then
$q-p=\frac{r}2-1\mod4$; \item if $r_2=1$, then
$q-p=\frac{r\pm1}2-1\mod4$; \item if $r_2=2$ and
$q-p=\frac{r}2+1\mod4$, then $\delta_2=0$; \item if
$\delta_2=0$, then $r=r_2=0\mod2$ and $q-p=\frac{r}2-1\mod2$;
\item if $r_2=r$ and $\delta_2=0$, then $p-q=\frac{r}2-1\mod4$.
\item if $r_3=r$, then $q-p=r-2\mod4$;
\item $r\ge r_3=p+q$.
\endroster
\endcorollary

\proof (1)--(6) follows from Proposition \ref{generalized-restrictions},
since for $T$-halves $\Br_3=2(p-q)\mod8$ and
$\sigma=2-r$. (7) follows from Proposition \ref{p-divisible}.
\qed\endproof

\demo{Proof of Proposition \ref{preliminary-list}}
It follows from the definition of  T-halves $M$ (see Lemma \ref{geometric-properties}(5)) and classification of the elementary 3-group (see Lemma \ref{3-classification})
that $\discr_3M=p\la\frac23\ra+q\la-\frac23\ra$, where $0\le p\le1$ and $0\le q\le3$.
 For each combination of $(r,r_2)$ and $\delta_2$ from Table 3A, we included in Table 4 those
 pairs $(p,q)$ which are not forbidden by one of
 the restrictions of Corollary \ref{T-half-restrictions} applied to M and M'.
\qed\enddemo

\remark{Remark}
 In a few cases Corollary \ref{T-half-restrictions} does not forbid a combination of $\delta_2$, $r$, $r_2$, $p$, and $q$, but forbids the complementary combination (i.e., for the other T-half of an  T-pair).
 For instance, the case $(r,r_2)=(7,1)$, $(p,q)=(1,0)$ is not forbidden, and is in fact realizable for the lattice
 $\la2\ra+\E_6$. However, the complementary lattice would have $r'=2$ and $(p',q')=(0,3)$, which is forbidden, since
 $r'<r_3'=3$. Thus, the combination $(r,r_2)=(7,1)$, $(p,q)=(1,0)$ is excluded from the list.
\boxedR\endremark

\subsection{The list of T-halves}\label{subsection-with-tables}
Below we provide an example of a T-half for every combination
of the numerical invariants from Table 4 (verification is straightforward, {\it cf.,} \cite{Nik1}, Theorem 1.10.1).

\proposition\label{T-halves-list}
For every combination of invariants
 $r$, $r_2$, $\delta_2$, $p$, $q$, or  $r'$, $r_2'$, $\delta_2'=1$, $p'$, $q'$ listed in the Table 4, there exists
 a T-half with such invariants. Namely, such a lattice can be found in
Tables 5A-J.
\qed\endproposition

Tables 5A-J contain only  T-halves. Sign ``-'' signifies that a lattice with given invariants is forbidden
by Corollary \ref{T-half-restrictions}. Sign ``*'' stands if a lattice actually exists, but it is not a T-half, because the complementary lattice is forbidden.
{\eightpoint
$$\matrix\text{\tenrm Table 5A. $r_2=0$ and $p=0$}
\\
\boxed{\matrix
(r,r_2)&   q=0&q=1&q=2&q=3\\
(2,0)&\U&-&-&-\\
(4,0)&-&\U+\A_2&-&-\\
(6,0)&-&-&\U+2\A_2&-\\
(8,0)&-&-&-&\U+3\A_2\\
\endmatrix}\endmatrix
\hskip-1mm \matrix\text{\tenrm Table 5B. $r_2=0$ and $p=1$}\\
\boxed{\matrix
q=0&q=1&q=2&q=3\\
-&\U(3)&-&-\\
-&-&\U(3)+\A_2&-\\
-&-&-&\U(3)+2\A_2\\
\U+\E_6&-&-&-\\
\endmatrix}\endmatrix$$
$$\matrix\text{\tenrm Table 5C. $r_2=1$ and p=0}\\
\boxed{\matrix
(r,r_2)& q=0&q=1&q=2&q=3\\
(1,1)&\la2\ra&-&-&-\\
(3,1)&\U+\A_1&\la2\ra+\A_2&-&-\\
(5,1)&-&\U+\A_2+\A_1&\la2\ra+2\A_2&-\\
(7,1)&-&-&\U+2\A_2+\A_1&\la2\ra+3\A_2\\
\endmatrix}\endmatrix
\hskip-1.2mm \matrix\text{\tenrm Table 5D. $r_2=1$ and
p=1}\\
\boxed{\matrix
q=0&q=1&q=2&q=3\\
\la6\ra&-&-&-\\
-&\la6\ra+\A_2&\U(3)+\la-6\ra&-\\
-&-&\la6\ra+2\A_2&\U(3)+\A_2+\la-6\ra\\
*&-&-&\la6\ra+3\A_2\\
\endmatrix}\endmatrix$$
$$\matrix\text{\tenrm Table 5E. $r_2=2$ and p=0}\\
\boxed{\matrix
\delta_2&(r,r_2)& q=0&q=1&q=2&q=3\\
0&(2,2)& \U(2)&-&-&-\\
1&(2,2)&\la2\ra+\A_1&\la2\ra+\la-6\ra&-&-\\
0&(4,2)&-& \U(2)+\A_2&-&\U(3)+\A_2(2)\\
1&(4,2)&\U+2\A_1&\U+\A_1+\la-6\ra&\U+2\la-6\ra&-\\
0&(6,2)&*&-&\U(2)+2\A_2&-\\
1&(6,2)&-&\U+\A_2+2\A_1&\U+\A_2+\la-6\ra+\A_1&\U+\A_2+2\la-6\ra\\
1&(8,2)&-&-&*&\la2\ra+3\A_2+\A_1\\
\endmatrix}\endmatrix$$
$$\matrix\text{\tenrm Table 5F. $r_2=2$ and p=1}\\
\boxed{\matrix
\delta_2&(r,r_2)& q=0&q=1&q=2&q=3\\
0&(2,2)&-&\U(6)&-&-\\
1&(2,2)&\la6\ra+\A_1&\la6\ra+\la-6\ra&-&-\\
0&(4,2)&\U+\A_2(2)&-&\U(6)+\A_2&-\\
1&(4,2)&-&\U(3)+2\A_1&\U(3)+\A_1+\la-6\ra&\U(3)+2\la-6\ra\\\
0&(6,2)&-&\U(3)+\DD_4&-&\U(6)+2\A_2\\
1&(6,2)&*&-&\U(3)+\A_2+2\A_1&\U(3)+\A_2+\A_1+\la-6\ra\\
1&(8,2)&*&*&-&\U(3)+2\A_2+2\A_1\\
\endmatrix}\endmatrix$$
$$\matrix\text{\tenrm Table 5G. $r_2=3$, or $r_2=5$ and p=0}\\
\boxed{\matrix
(r,r_2)& q=0&q=1&q=2&q=3\\
(3,3)&\la2\ra+2\A_1&\la2\ra+\A_1+\la-6\ra&\la2\ra+2\la-6\ra&-\\
(5,3)&\U+3\A_1&\U+2\A_1+\la-6\ra&\U+\A_1+2\la-6\ra&\U+3\la-6\ra\\
(5,5)&\U(2)+3\A_1&\U(2)+2\A_1+\la-6\ra&\U(2)+\A_1+2\la-6\ra&\U(2)+3\la-6\ra\\
(7,3)&*&*&\U+\A_2+\la-6\ra+2\A_1&\U+\A_2+2\la-6\ra+\A_1\\
\endmatrix}\endmatrix$$
$$\matrix\text{\tenrm Table 5H. $r_2=3$, or $r_2=5$ and p=1}\\
\boxed{\matrix
(r,r_2)& q=0&q=1&q=2&q=3\\
(3,3)&\la6\ra+2\A_1&\la6\ra+\A_1+\la-6\ra&\la6\ra+2\la-6\ra&-\\
(5,3)&\la6\ra+\DD_4&\U(3)+3\A_1&\U(3)+2\A_1+\la-6\ra&\U(3)+\A_1+2\la-6\ra\\
(5,5)&\la6\ra+4\A_1&\U(6)+3\A_1&\U(6)+2\A_1+\la-6\ra&\U(6)+\A_1+2\la-6\ra\\
(7,3)&*&*&\U(3)+\A_2+3\A_1&\U(3)+\A_2+\la-6\ra+2\A_1\\
\endmatrix}
\endmatrix$$
$$\matrix\text{\tenrm Table 5I. $r_2=4$ and p=0}\\
\boxed{\matrix
\delta_2&(r,r_2)& q=0&q=1&q=2&q=3\\
0&(4,4)&-&-&-&\U(6)+\A_2(2)\\
1&(4,4)&\la2\ra+3\A_1&\la2\ra+2\A_1+\la-6\ra&\la2\ra+\A_1+2\la-6\ra&\la2\ra+3\la-6\ra\\
1&(6,4)&*&\U+\la-6\ra+3\A_1&\U+2\la-6\ra+\A_1&\U+3\la-6\ra+\A_1\\
\endmatrix}\endmatrix$$
$$\matrix\text{\tenrm Table 5J. $r_2=4$ and p=1}\\
\boxed{\matrix
\delta_2&(r,r_2)& q=0&q=1&q=2&q=3\\
0&(4,4)&\U(2)+\A_2(2)&-&-&-\\
1&(4,4)&\la6\ra+3\A_1&\la6\ra+\la-6\ra+2\A_1&\la6\ra+2\la-6\ra+\A_1&\la6\ra+3\la-6\ra\\
1&(6,4)&*&\U(3)+4\A_1&\U(3)+\la-6\ra+3\A_1&\U(3)+2\la-6\ra+2\A_1\\
\endmatrix}
\endmatrix$$}
\rk{Remark}
Note that there is some ambiguity in the direct sum presentations
of the lattices in the above tables, due to the following
isomorphisms:
$\la6\ra+\A_2=\U(3)+\A_1$, $\la2\ra+\A_2=\U+\la-6\ra$,
$\la6\ra+\A_2(2)=\la2\ra+2\la-6\ra$,
$\la2\ra+\A_2(2)=\la6\ra+2\A_1$.

Besides, for any lattice $L$ with an odd $\discr_2L$ $($for
instance, $L=\A_1$, or $L=\la-6\ra$~$)$ we have
$\U(2)+L=\la2\ra+\A_1+L$, and $\U(6)+L=\la6\ra+\la-6\ra+L$.
\hskip2.2in \boxedR\endrk

\subsection{Stability of the  T-halves}
Our aim in the next four subsections is to prove stability of the  T-halves, and thus,
to show that their list in Tables 5A-J is complete.

It is easy to check which of these lattices do not satisfy
Nikulin's stability criterion (see \ref{stability-criterion}) using Table 4.
We see from it that Nikulin's condition for $r_2$, namely, $r>r_2$, or $r=r_2>2$,
fails only for the T-halves listed in Table 6A below.
 Table 6B contains the remaining cases, in which Nikulin's condition $r_3\le r-2$ fails.

\lemma
For any  T-pair $(T_1,T_2)$,
Nikulin's stability criterion \ref{stability-criterion} is satisfied for $T_i$, $i\in\{1,2\}$,
unless the combination of the invariants $r(T_i)$, $r_2(T_i)$, $\delta_2(T_i)$, $p(T_i)$ and $q(T_i)$
is among the ones listed in Tables 6A-B.
\qed\endlemma
\midinsert
$\phantom{AAAAAAA}\matrix\text{Table 6A.}\ \
r_2=r\le2\\
\matrix \\ 1\\ 2\\ 3\\ 4\\ 5\\ 6\\ 7\\ 8 \endmatrix
\boxed{\matrix
\delta_2&(r,r_2)&(p,q)&T_i\\
1&(1,1)&(0,0)&\la2\ra\\
1&(1,1)&(1,0)&\la6\ra\\
0&(2,2)&(0,0)&\U(2)\\
0&(2,2)&(1,1)&\U(6)\\
1&(2,2)&(0,0)&\la2\ra+\A_1\\
1&(2,2)&(0,1)&\la2\ra+\la-6\ra\\
1&(2,2)&(1,0)&\la6\ra+\A_1\\
1&(2,2)&(1,1)&\la6\ra+\la-6\ra\\
\endmatrix}
\endmatrix$
\endinsert

\midinsert
$$\matrix\text{Table 6B.} \  r_3\ge r-1\phantom{AAAAAAAAAA}\\
\matrix \\ 1\\ 2\\ 3\\ 4\\ 5\\ 6\\ 7\\ 8\\ 9\\ 10\\ 11\\ 12\\ 13\\ 14\\ 15\\ 16\\ 17\\ 18\endmatrix
\boxed{\matrix
\delta_2&(r,r_2)&(p,q)&T_i\\
0&(2,0)&(1,1)&\U(3)\\
0&(4,0)&(1,2)&\U(3)+\A_2\\
1&(3,1)&(1,1)&\U(3)+\A_1\\
1&(3,1)&(1,2)&\U(3)+\la-6\ra\\
1&(5,1)&(1,3)&\U(3)+\A_2+\la-6\ra\\
0&(4,2)&(0,3)&\U(3)+\A_2(2)\\
0&(4,2)&(1,2)&\U(6)+\A_2\\
1&(4,2)&(1,2)&\U(3)+\A_1+\la-6\ra\\
1&(4,2)&(1,3)&\U(3)+2\la-6\ra \\
1&(3,3)&(0,2)&\la2\ra+2\la-6\ra\\
1&(3,3)&(1,1)&\la6\ra+\A_1+\la-6\ra\\
1&(3,3)&(1,2)&\la6\ra+2\la-6\ra\\
1&(5,3)&(1,3)&\U(3)+\A_1+2\la-6\ra\\
0&(4,4)&(0,3)&\U(6)+\A_2(2)\\
1&(4,4)&(0,3)&\la2\ra+3\la-6\ra\\
1&(4,4)&(1,2)&\la6\ra+\A_1+2\la-6\ra\\
1&(4,4)&(1,3)&\la6\ra+3\la-6\ra\\
1&(5,5)&(1,3)&\U(6)+\A_1+2\la-6\ra \\
\endmatrix}
\matrix
\\ r_3=r\\ r_3=r-1\\  r_3=r-1\\ r_3=r\\ r_3=r-1\\  r_3=r-1\\ r_3=r-1\\ r_3=r-1\\  r_3=r\\ r_3=r-1\\ r_3=r-1\\  r_3=r\\
 r_3=r-1\\ r_3=r-1\\  r_3=r-1\\ r_3=r-1\\ r_3=r\\  r_3=r-1
\endmatrix \hskip3mm
\matrix
\\ r_2=0 \\ r_2=0 \\ r_2<r \\ r_2<r \\ r_2<r \\ r_2<r \\ r_2<r \\ r_2<r \\ r_2<r \\ r_2=r \\ r_2=r \\
 r_2=r \\  r_2<r \\ r_2=r \\  r_2=r \\ r_2=r \\  r_2=r \\ r_2=r \\
\endmatrix
\endmatrix$$
\endinsert

\corollary\label{either-pm-stability}
For any  T-pair $(T_1,T_2)$ either $T_1$ or $T_2$ satisfies Nikulin's stability condition, and so
is stable and epistable.
\endcorollary

\proof
Is can be easily seen from Table 4, that there is no  pairs $(T_1,T_2)$ for which
both $T_1$ and $T_2$ are represented in the above list of exceptions.
\endproof

Our next aim is to verify stability of the lattices in Tables 6A-B.
The Miranda-Morrison criterion (see Proposition \ref{MM-23elementary}) gives the following statement.

\proposition\label{MM-corollary}
Assume that $M$
is a T-half of rank $r>2$ and that one of the discriminant $p$-ranks $r_p$, $p=2,3$
of $M$  is less than $r$. Then $M$ is stable and epistable.
\qed
\endproposition

 The following fact is well known.

\lemma\label{unimodular-hyperbolic} The only unimodular hyperbolic
lattices of rank $\le8$ are $\U$ and $\la1\ra+n\la-1\ra$, $n\le7$.
\qed\endlemma

\lemma\label{simplest-lemma}
Assume that $M$ is a T-half.
Let $r$ be its rank and $r_p$, $p=2,3$ its discriminant $p$-ranks. Then:
\roster
\item If $r_2=0$ and $r_3=r$, then $M=\U(3)$. 
\item If
$r_2=r$ and $r_3=0$, then either $\delta_2=0$ and $M=\U(2)$, or
$\delta_2=1$ and $M=\la2\ra+n\la-2\ra$, $n\le7$. \item If
$r_2=r_3=r$, then either $\delta_2=0$ and $M=\U(6)$, or
$\delta_2=1$ and $M=\la6\ra+n\la-6\ra$, $n\le7$.
\endroster\endlemma

\proof It is enough to apply Proposition \ref{T-halves-properties}
and Lemma \ref{unimodular-hyperbolic}
to $M(\frac13)$, $M(\frac12)$, and $M(\frac16)$ in the cases (1),
(2), and (3), respectively (the estimate $n\le7$ is due to Proposition \ref{T-halves-properties}(1)).
\endproof

\lemma\label{r2=r=2}  Under the assumption of Lemma \ref{simplest-lemma}, if $r_2=r=2$ and $r_3=1$, then $L$ is either
$\la2\ra+\la-6\ra$, or $\la6\ra+\la-2\ra$.
\endlemma

\proof The lattice $L'=L(\frac12)$ has discriminant order $|\discr L'|=3$.
 Fixing any basis of $L$, we can present it by a
matrix, $\left(\matrix a&b\\b&c\endmatrix\right)$, $ac-b^2=-3$,
(the sign of the discriminant is due
 to that $L$ is indefinite). The gaussian theory of normal
 forms says that in a suitable basis this matrix has $0\le
 b\le\sqrt 3$, i.e., $b=0$, or $b=1$.
 If $b=1$, then $ac=b^2-3=-2$, and the only solutions
(up to reordering the basis) are
$\left(\matrix\pm1&1\\1&\mp2\endmatrix\right)$
 and it is easy to check that the both matrices are
 diagonalizable. Diagonalization yields
 $\left(\matrix\pm3&0\\0&\mp1\endmatrix\right)$, that is
 $\la1\ra+\la-3\ra$, or $\la3\ra+\la-1\ra$.
\endproof

\corollary\label{simplest-corollary}
The lattices of Table 6A and the lattices of Table 6B having $r_3=r$,
namely, the ones in rows 1,4,9,12, and 17, are all stable.
\endcorollary

\proof
Lattices of Table 6A satisfy either assumptions (2)--(3) of
Lemma \ref{simplest-lemma}, or the ones of Lemma \ref{r2=r=2}.
 For the indicated lattices of Table 6B we apply Lemma \ref{simplest-lemma}.
\endproof

\theorem\label{Stability-of-T-halves}
All the lattices in Tables 5A-J are stable.
\endtheorem

\proof The cases of lattices with $r\le 2$ for which  Nikulin's and Miranda-Morrisson's criteria do not
guarantee stability are covered by Lemmas \ref{simplest-lemma} and
\ref{r2=r=2}.
The case (3) of Lemma \ref{simplest-lemma} covers also $r_2=r_3=r\ge3$.
In the remaining cases Proposition \ref{MM-corollary} can be applied.
\endproof

\corollary\label{completeness-T-halves}
Tables 5A-J give a complete list of  T-halves.
\endcorollary

\proof
Proposition \ref{preliminary-list} gives a list of the invariants $(r,r_2,\delta_2,p,q)$ which can be
potentially realized for  T-halves. Proposition \ref{T-halves-list} shows that this list is actually
realized by lattices in Tables 5A-J. Theorem \ref{Stability-of-T-halves} shows that there is no other
T-halves with the same invariants, and so, Tables 5A-J give a complete list.
\endproof

\corollary\label{T-half-via-invariants}
A T-half $M$ is determined by its rank $r$, 2-rank
$r_2$, parity $\delta_2$ of $\discr_2M$, and the characteristics $0\le p\le1$, $0\le q\le3$
of $\discr_3M=p\dip+q\din$. In particular, $M$ is determined by $r$ and $\discr M$.
\endcorollary

If we combine  T-halves in Tables 5A-J into ascending  T-pairs in accordance  with the Table 4,
then we obtain the following result.

\midinsert
$$
\matrix \\ \\
1\\2\\3\\4\\5\\6\\7\\8\\9\\10\\11\\12\\13\\14\\15\\16\\17\\18\\19\\20\\21\\22\\23\\24\\25\\26\\27\\28\\29\\
30\\31\\32\\33\\34\\35\\
\endmatrix
\matrix\text{Table 7A. Ascending  T-pairs $T_\pm$ in the case of $p=0$}\\
\boxed{\matrix
\delta_2&(r,r_2)&q&T_+ &             T_-\\
 0&(2,0)&0&\U         &\U(3)+2\A_2+\A_1\\
 0&(4,0)&1&\U+\A_2       &\U(3)+\A_2+\A_1& \\
 0&(6,0)&2&\U+2\A_2     &\U(3)+\A_1& \\
 0&(8,0)&3&\U+3\A_2         &\la6\ra& \\
 1&(1,1)&0&\la 2\ra              &\U(3)+2\A_2+2\A_1\\
 1&(3,1)&0&\U+\A_1       & \la6\ra+\la -6\ra+2\A_2 \\
 1&(3,1)&1&\U+\la -6\ra       &\U(3)+\A_2+2\A_1\\
 1&(5,1)&1&   \U+\A_2+\A_1    & \la6\ra+\la -6\ra+\A_2\\
 1&(5,1)&2&\la2\ra+2\A_2      &\U(3)+2\A_1& \\
 1&(7,1)&2&\U+2\A_2+\A_1   & \la6\ra +\la-6\ra& \\
 1&(7,1)&3   &\la 2\ra+3\A_2   &\la 6\ra+\A_1& \\
 0&(2,2)&0&\U(2)     & \la6\ra+\la -6\ra+2\A_2+\A_1\\
 1&(2,2)&0&\la 2\ra+\A_1     & \la6\ra+\la -6\ra+2\A_2+\A_1\\
 1&(2,2)&1&\la 2\ra+\la -6\ra         &\la 6\ra+2\A_2+2\A_1& \\
 0&(4,2)&3&\U(3)+\A_2(2)      &\U+\A_2(2)+\A_1 & \\
 1&(4,2)&0  &\U+2\A_1  & \la6\ra+2\la -6\ra+\A_2& \\
 1&(4,2)&1&  \la 2\ra+\A_2+\A_1         & \la6\ra+\la -6\ra+\A_2+\A_1 \\
 1&(4,2)&2&\U+2\la -6\ra      &\U(3)+3\A_1\\
 0&(4,2)&1& \U(2)+\A_2         &  \la6\ra+\la -6\ra+\A_2+\A_1\\
 0&(6,2)&2&  \U(2)+2\A_2           &\la6\ra+\la -6\ra+\A_1\\
 1&(6,2)&2   &\U+\A_2+\la -6\ra+\A_1       & \la6\ra+\la -6\ra+\A_1 \\
 1&(6,2)&3   &\la 2\ra+2\A_2+\la-6\ra  &\la6\ra+2\A_1 \\
 1&(6,2)&1&   \U+\A_2+2\A_1  & \la6\ra+2\la -6\ra&\\
 1&(3,3)&0&\la 2\ra+2\A_1       & \la6\ra+2\la -6\ra+\A_2+\A_1 \\
 1&(3,3)&1&\la 2\ra+\A_1+\la -6\ra  &\la6\ra+\la -6\ra+\A_2+2\A_1 \\
 1&(3,3)&2&\la 2\ra+2\la -6\ra     &\U(3)+4\A_1& \\
 1&(5,3)&3&\U+3\la -6\ra       &\la 6\ra+3\A_1&\\
 1&(5,3)&2& \la 2\ra+\A_2+\A_1+\la -6\ra & \la6\ra+\la -6\ra+2\A_1 & \\
 1&(5,3)&1&  \la 2\ra+\A_2+2\A_1       & \la6\ra+2\la -6\ra+\A_1 \\
 1&(5,3)&0  &\U+3\A_1 & \la6\ra+3\la -6\ra \\
 1&(4,4)&0& \la 2\ra+3\A_1  & \la6\ra+3\la -6\ra+\A_1& \\
 1&(4,4)&1& \la2\ra+2\A_1+\la-6\ra       &\la6\ra+2\la-6\ra+\A_1&\\
 1&(4,4)&2&\la 2\ra+\A_1+2\la -6\ra &\la6\ra+\la -6\ra+3\A_1 \\
 1&(4,4)&3&\la 2\ra+3\la -6\ra      &\la 6\ra+4\A_1& \\
 0&(4,4)&3&\U(6)+\A_2(2)  &\la6\ra+4\A_1\\
\endmatrix}\endmatrix$$
\endinsert

\midinsert
$$
\matrix \\ \\
1\\2\\3\\4\\5\\6\\7\\8\\9\\10\\11\\12\\13\\14\\15\\16\\17\\18\\19\\20\\21\\22\\23\\24\\25\\26\\27\\28\\29\\30\\31\\32
\endmatrix
\matrix\text{Table 7B. Ascending  T-pairs $T_\pm$ in the case of $p=1$}\\
\boxed{\matrix
   \delta_2&(r_+,r_2)&q&T_+ &         T_-\\
 0&(2,0)&1& \U(3)          &\U+3\A_2 \\
 0&(4,0)&2&\U(3)+\A_2        &\U+\A_2+\A_1& \\
 0&(6,0)&3&   \U(3)+2\A_2   & \U+\A_1\\
 1&(1,1)&0  &\la 6\ra           &\U+2\A_2+\A_1+\la-6\ra \\
 1&(3,1)&1   &\la 6\ra+\A_2            &\la 2\ra+2\A_2+\A_1 \\
 1&(3,1)&2&   \U(3)+\la-6\ra    &\U+\A_2+2\A_1\\
 1&(5,1)&2&   \U(3)+\A_2+\A_1  &\U+\la-6\ra+\A_1\\
 1&(5,1)&3&   \U(3)+\A_2+\la-6\ra   &\U+2\A_1\\
 1&(7,1)&3&\la 6\ra+3\A_2  &\la 2\ra+\A_1&\\
 0&(2,2)&1& \U(6)    &\U(2)+2\A_2+\A_1 \\
 1&(2,2)&1&  \la6\ra+\la-6\ra      &\U+\A_2+2\A_1+\la-6\ra \\
 1&(2,2)&0   &\la 6\ra+\A_1   &\la2\ra+2\A_2+\A_1+\la-6\ra \\
 0&(4,2)&0&\U+\A_2(2)    &\U(3)+\A_2(2)+\A_1 & \\
 0&(4,2)&2&\U(6)+\A_2       &\U(2)+\A_2+\A_1& \\
 1&(4,2)&2&  \la6\ra+\A_2+\la-6\ra      &\la2\ra+\A_2+2\A_1 \\
 1&(4,2)&1   &\la 6\ra+\A_2+\A_1         &\la2\ra+\la-6\ra+\A_2+\A_1\\
 1&(4,4)&3&  \la 6\ra+3\la-6\ra          &\la 2\ra+4\A_1 \\
 1&(4,2)&3&   \U(3)+2\la-6\ra    &\U+3\A_1\\
 0&(6,2)&3& \U(6)+2\A_2       &\U(2)+\A_1& \\
 0&(6,2)&1&   \U(3)+\DD_4       &\la6\ra+\A_2(2)&  \\
 1&(6,2)&2&   \la 6\ra+2\A_2+\A_1 &\la 2\ra+\A_1+\la -6\ra \\
 1&(6,2)&3&  \U(3)+\A_2+\A_1+\la-6\ra     &\la2\ra+2\A_1 \\
 1&(3,3)&0   &\la 6\ra+2\A_1        & \la2\ra+2\la-6\ra+\A_2+\A_1 \\
 1&(3,3)&1& \la 6\ra+\A_1+\la-6\ra      &\la 2\ra+\A_2+2\A_1+\la -6\ra \\
 1&(3,3)&2& \la6\ra+2\la-6\ra         &\U+\la-6\ra+3\A_1 \\
 1&(5,3)&0   &\la 6\ra+\DD_4           &\la6\ra+\A_2(2)+\la -6\ra\\
 1&(5,3)&1&\la 6\ra+\A_2+2\A_1  &\la 2\ra+\A_1+2\la -6\ra &  \\
 1&(5,3)&3&\la 6\ra+\A_2+2\la-6\ra         &\la 2\ra+3\A_1 \\
 1&(5,3)&2&\la 6\ra+\A_2+\A_1+\la-6\ra          &\la 2\ra+2\A_1+\la -6\ra \\
 1&(4,4)&0&\la 6\ra+3\A_1      &\la 2\ra+\A_1+3\la -6\ra  & \\
 1&(4,4)&1   &\la 6\ra+\la -6\ra+2\A_1          &\la2\ra+2\A_1+2\la -6\ra \\
 1&(4,4)&2 & \la6\ra+\A_1+2\la-6\ra        &\la2\ra+3\A_1+\la-6\ra \\
\endmatrix}\endmatrix$$
\endinsert

 \theorem\label{eigenlattices-enumeration}
 Tables 7A-B give a complete list of the ascending T-pairs
 $(T_+,T_-)$.\!\qed \endtheorem

\subsection{Property of ($-\frac12$)-transitivity for the eigenlattices $T_-$}\label{-1/2-transitivity}
Here we establish a certain property of T-halves $T_-$ to be used in the next section.
We say that a lattice $M$ is {\it $(-\frac12)$-transitive} if the
automorphism group $\Aut(M)$ induces a transitive action on the
non-characteristic elements $x\in\discr_2(M)\subset\discr M$ with
$\q_M(x)=-\frac12$. Note that for the characteristic element,
$v\in\discr_2(M)$ the value $\q_M(v)$ is determined by the Brown
invariant of $\discr_2M$. For
instance, if $\discr_2M=p\la\frac12\ra+q\la-\frac12\ra$, then
$\q_M(v)=\frac{p-q}2\mod2\Z$, so, $\q_M(v)=-\frac12$ if and only if
$\Br_2(M)=3\mod4$.

\proposition\label{prop-1/2-transitivity} All the eigenlattices $T_-$ in the Tables 5 and 6A-B
are $(-\frac12)$-transitive.
\endproposition

\proof First, note that the automorphisms of $\discr_2T_-$ act
transitively on the non-characte\-ristic ($-\frac12$)-elements.

\lemma For any elementary enhanced 2-group $G$ the isometry group
$\Aut(G)$ acts transitively on the non-characteristic elements
$x\in G$ having $\q(x)=-\frac12$.
\endlemma

\proof If $V$ contains $v\in G$ such that $\q(v)=-\frac12$, then
$q$ is odd and $G=p\la\frac12\ra+q\la-\frac12\ra$, where
$p-q=\Br(\q)$. The orthogonal complement, $G^v$, of $v$ is odd if
and only if $v$ is not characteristic, and in the latter case,
$G^v=\la\frac12\ra+(q-1)\la-\frac12\ra$. Thus, for any other
non-characteristic element $w\in G$ with $\q(w)=-\frac12$, we
obtain an isomorphism $G^v\to G^w$ which is extended to $G$ so
that $v$ is sent to $w$.
\endproof

It follows that 2-epistability implies ($-\frac12$)-transitivity.
So, we should analyze only the lattices $T_-$ for which Nikulin's
stability criterion fails. Among these cases there are trivial
ones, with $\discr_2 T_-$ containing only one or no
non-characteristic $(-\frac12)$-elements. This happens if
$r_2(T_-)\le1$, or if
 $\discr_2(T_-)$ is isomorphic to
$2\la\frac12\ra$, $\la\frac12\ra+\la-\frac12\ra$,
$3\la\frac12\ra$, or $2\la\frac12\ra+\la-\frac12\ra$.

 So, it is sufficient
to consider only those lattices $T_-$ non
 satisfying Nikulin's criterion for which
 $\discr_2(T_-)$ is $2\la-\frac12\ra$,
$\la\frac12\ra+2\la-\frac12\ra$, $3\la-\frac12\ra$,
 or $r_2(T_-)\ge4$.
Analyzing the lattices $T_-$ in Tables 5 and 6A-B, in Section
\ref{subsection-with-tables}, we find that there remain only the
following cases to consider:
 \roster \item
$\la6\ra+3\la-6\ra$,
 \item
$\la2\ra+3\la-6\ra$,
 \item
 $\la6\ra+k\la-6\ra+\A_1$, $0\le k\le3$.
\endroster

If $T_-=\la6\ra+3\la-6\ra$, then $T_-(\frac13)$ is epistable by Nikulin's criterion.
In the cases (2) and (3) with $k>0$, Proposition \ref{MM-23elementary} is applicable.

In the remaining case,
$T_-=\la6\ra+\A_1$, the discriminant component $\discr_2T_-=2\la-\frac12\ra$ has
 only one non-identity automorphism.
This automorphism interchanges the summands of $\discr_2T_-$, and therefore it is induced by the reflection $\rho_v$ in $T_-$
with respect to $v=(1,1)\in T_-$. Thus, $T_-$ is 2-epistable and,
in particular, $(-\frac12)$-transitive.
\endproof

\section{Back to Zariski curves}\label{Chapter_back}

\subsection{Classification of geometric involutions}\label{classification-involutions}
\theorem\label{eigenlattices&involutions}
Geometric involutions $c,c'\in C(L,\D,h)$
have the same homological type $($i.e., $[c]=[c']\in C[L,\D,h])$
if and only if the corresponding  T-pairs $(T_+(c),T_-(c))$ and $(T_+(c'),T_-(c'))$ are isomorphic
$($i.e., $T_\pm(c)\cong T_\pm(c'))$.
\endtheorem

\proof
The ``only if'' part is trivial. Assume that there are isomorphisms $T_\pm(c)\cong T_\pm(c')$.
 Without loss of generality we may assume that $r_2(T_+)<r_2(T_-)$, i.e., that the involutions $c$ and $c'$ are ascending. Using Proposition \ref{conjugate-prop}
 we can observe
that the restrictions $c|_T$ and $c'|_T$ to the sublattice $T\subset L$ are conjugate via some automorphism
$f_T\:T\to T$.
 Namely, Corollary \ref{either-pm-stability} shows that either $T_-$ or $T_+$ is epistable, and in the first case,
the condition (1) of Proposition \ref{conjugate-prop} is satisfied.
 In the case  $T_+$ is epistable, the condition (2) of Proposition \ref{conjugate-prop} includes,
 in addition, a certain transitivity assumption;
 this assumption is satisfied according to
Lemma \ref{anti-isomorphism-criterion}
(implying that transitivity
of the action of $\Aut(T_-)$ on the subgroups of $\discr_2(T_-)$ anti-isomorphic to $\discr_2(T_+)$ is equivalent to
$(-\frac12)$-transitivity of this group action on $T_-$),
and Proposition \ref{prop-1/2-transitivity} (which establishes the latter transitivity).

 Next, we extend $f_T$ to an automorphism $f_{T'}\:T'\to T'$ by letting $h\mapsto h$
(see Section \ref{2-extensions}). Since $c(h)=c'(h)=-h$, the restrictions $c|_{T'}$ and $c'|_{T'}$
 are conjugate via $f_{T'}$.

Finally, using Proposition \ref{gluing-equi-automorphisms} and Corollary \ref{S-equi-epistability} we conclude that $f_{T'}$ can be extended to
$f\in\Aut(L,\D,h)$ that conjugates $c$ with $c'$.
\endproof

\subsection{Reversion roots}
 Given an ascending  T-pair $(T_1,T_2)$,
we say that an element $v\in T_2$,
is a {\it reversion root for $(T_1,T_2)$} if, first, it is an even $(-2)$-element
(see the definitions in Section \ref{2-extensions}),
and, second, $[\frac{v}2]\in\discr_2T_2$
is a characteristic element in the case $\delta_2(T_1)=0$, and
non-characteristic in the case $\delta_2(T_1)=1$.
Due to Lemma \ref{Wu-and-even}) the latter condition implies
the following property of the orthogonal complement $T_2^v\subset T_2$ of $v$.

\corollary\label{reversion-parity}
$\delta_2(T_1)=\delta_2(T_2^v)$.
\qed\endcorollary

\subsection{Reversion partners of  T-pairs}
Given an ascending  T-pair $(T_1,T_2)$ and a reversion root $v\in T_2$,
we can write $T_2=\Z v+T_2^v$. Interchanging $T_1$ and $T_2^v$ we obtain another pair
 $(T_2^v,\Z v+T_1)$, which will be called
{\it the reversion partner of $(T_1,T_2)$ with respect to $v$}.
If a reversion root does exist, we say that {\it $(T_1,T_2)$ is a reversible pair}.

 \lemma\label{partnership-lemma}
 The reversion partner $(T_1',T_2')=(T_2^v,T_v)$ is also an ascending  T-pair. Moreover,
 \roster\item
 $r(T_i')=8-r(T_i)$, $i=1,2$;
 \item
 $r_2(T_i')=r_2(T_i)$, and $\delta_2(T_i')=\delta_2(T_i)$;
 \item
$\discr_3(T_1')=\discr_3(T_2)$, $\discr_3(T_2')=\discr_3(T_1)$,
and, in particular, if $\discr_3(T_i)=p\dip+q\din$ then
$\discr_3(T_i')=(1-p)\dip+(3-q)\din$.
 \endroster
 Furthermore, $v\in T_2'$ is a reversion root of $(T_1',T_2')$ and the reversion partner of $(T_1',T_2')$ with respect to $v$ is $(T_1,T_2)$.
 \endlemma

\proof
 From the construction of $(T_1',T_2')$ it follows that $r(T_1')=r(T_2)-1=(9-r(T_1))-1$, and $r(T_2')-1=r(T_1)=9-r(T_2)$, which gives (1). Next, $r_2(T_1')=r_2(T_2)-1=r_2(T_1)$, and
Corollary \ref{reversion-parity} yields (2).
Lemma \ref{evenness-criteria} implies that $\discr_3(v^\perp)=\discr_3(T_1)$, which implies (3).
 The last claim is obvious from the construction.
\endproof

\proposition For any ascending  pair $(T_1,T_2)$, its
 reversion partner $(T_1',T_2')$ if exists is independent up to isomorphism of the
choice of a reversion root $v$.
\endproposition

\proof
By Lemma \ref{partnership-lemma} the invariants $r$, $r_2$,  $\delta_2$, $p$, and $q$ of $T_i'$ are determined by those of $T_i$.
These invariants determine an ascending  T-pair uniquely up to isomorphism, see Corollary \ref{T-half-via-invariants}.
\endproof

\theorem
A reversion partner
 exists for all the ascending  T-pairs except six pairs $(T_1,T_2)$ listed in Table 8A.
The partnership of the remaining 62 ascending  T-pairs is shown in Tables 8B-C, where partners
are placed in the same rows.
\endtheorem
{\eightpoint
\midinsert
\hskip20mm$
\matrix\\ \\ \\
1\\2\\3\\4\\5\\6
\endmatrix
\matrix\text{\tenrm Table 8A. Irreversible  T-pairs}
\\
\boxed{\matrix
\delta_2&(r,r_2)&(p,q)&(r',r_2')&(p',q')&T_1&T_2\\
\text{------}&\text{------}&\text{------}&\text{------}&\text{------}&\text{---------}&\text{---------}\\
0&(8,0)&(0,3)&(1,1)&(1,0)&\U+3\A_2&\la6\ra\\
1&(7,1)&(0,2)&(2,2)&(1,1)&\U+2\A_2+\A_1&\la6\ra+\la-6\ra\\
0&(6,2)&(1,1)&(3,3)&(0,2)&\U(3)+\DD_4&\la6\ra+\A_2(2)\\
1&(6,2)&(0,1)&(3,3)&(1,2)&\U+\A_2+2\A_1&\la6\ra+2\la-6\ra\\
1&(5,3)&(0,0)&(4,4)&(1,3)&\U+3\A_1&\la6\ra+3\la-6\ra\\
1&(5,3)&(1,0)&(4,4)&(0,3)&\la6\ra+\DD_4&\la2\ra+3\la-6\ra\\
\endmatrix}
\endmatrix$
\endinsert
}

\proof
 For the reversion partners the invariants  $r$, $r_2$,  $\delta$, $p$, and $q$ should match as is indicated in Lemma
\ref{partnership-lemma}.
For the six exceptional T-pairs in Table 8A there is no candidates to be a partner with the matching invariants.
 For the other T-pairs such a candidate is unique, since these invariants uniquely determine an ascending  T-pair by Corollary \ref{T-half-via-invariants}. It is straightforward to check that
 the T-pairs placed in the same rows in Tables 8B-C have
matching above invariants.
 Now, it is sufficient to check existence of a reversion root $v\in T_2$ for every T-pair $(T_1,T_2)$ in Table 8B
(this implies existence of a reversion root in the partner T-pair in Table 8C). This is also straightforward and easy.
\endproof

{\eightpoint
\midinsert
\hskip-12mm
$
\matrix\\ \\
1\\2\\3\\ \\ 4\\5\\6\\7\\8\\9\\ \\ 10\\11\\12\\13\\14\\15\\16\\17\\18\\19\\20\\ \\
21\\22\\23\\24\\25\\26\\ \\ 27\\28\\29\\30\\ \\ 31
\endmatrix
\matrix\text{\tenrm Table 8B. The case of $o=-$ ($p=0$)}\\
\boxed{\matrix
(r,r_2)&q&\delta_2&T_1&T_2\\
(2,0)&0&0&\U                 &\U(3)+2\A_2+\A_1\\
(4,0)&1&0&\U       +\A_2    &\U(3)+\A_2+\A_1& \\
(6,0)&2&0&\U       +2\A_2    &\U(3)+\A_1& \\
\\
(1,1)&0&1&\la 2\ra       &\U(3)+2\A_2+2\A_1\\
(3,1)&0&1&\U       +\A_1           & \la6\ra+\la -6\ra+2\A_2 \\
(3,1)&1&1&\U      +\la -6\ra      &\U(3)+\A_2+2\A_1\\
(5,1)&1&1&  \U    +\A_2+\A_1   & \la6\ra+\la -6\ra+\A_2\\
(5,1)&2&1& \la2\ra+2\A_2    &\U(3)+2\A_1& \\
(7,1)&3&1& \la 2\ra+3\A_2  &\la 6\ra+\A_1& \\
\\
(2,2)&0&0&\U(2)          &\la6\ra\!+\!2\A_2\!+\!\la-6\ra\!+\!\A_1\\
(2,2)&0&1&\la 2\ra+\A_1  &\la6\ra\!+\!2\A_2\!+\!\la-6\ra\!+\!\A_1\\
(2,2)&1&1&\la 2\ra +\la -6\ra     &\la 6\ra+2\A_2+2\A_1 \\
(4,2)&2&1&\U        +2\A_1 & \la6\ra+2\la -6\ra+\A_2 \\
(4,2)&1&0&\U(2)  +\A_2         &  \la6\ra\!+\!\la -6\ra\!+\!\A_2\!+\!\A_1\\
(4,2)&1&1& \la 2\ra +\A_2+\A_1    & \la6\ra\!+\!\la -6\ra\!+\!\A_2\!+\!\A_1 \\
(4,2)&2&1&\U      +2\la -6\ra       &\U(3)+3\A_1\\
(4,2)&3&0& \U(3) +\A_2(2)     &\U+\A_2(2)+\A_1 & \\
(6,2)&2&0&  \U(2) +2\A_2             &\la6\ra+\la -6\ra+\A_1\\
(6,2)&2&1& \U \!+\!\A_2\!+\!\la -6\ra\!+\!\A_1    &\la6\ra+\la -6\ra+\A_1 \\
(6,2)&3&1& \la 2\ra +2\A_2+\la-6\ra  &\la6\ra+2\A_1 \\
\\
(3,3)&2&1&\la 2\ra +2\A_1     & \la6\ra\!+\!2\la -6\ra\!+\!\A_2\!+\!\A_1 \\
(3,3)&1&1&\la 2\ra +\la -6\ra+\A_1  &\la6\ra\!+\!\la -6\ra\!+\!\A_2\!+\!2\A_1 \\
(3,3)&2&1&\la 2\ra +2\la -6\ra   &\U(3)+4\A_1& \\
(5,3)&1&1&\la 2\ra +\A_2 +2\A_1         & \la6\ra+2\la -6\ra+\A_1 \\
(5,3)&2&1& \la 2\ra\!+\!\A_2\!+\!\la -6\ra\! +\!\A_1 & \la6\ra+\la -6\ra+2\A_1 & \\
(5,3)&3&1& \U     +3\la -6\ra      &\la 6\ra+3\A_1&\\
\\
(4,4)&0&1& \la 2\ra +3\A_1  & \la6\ra+3\la -6\ra+\A_1& \\
(4,4)&1&1& \la2\ra +\la-6\ra +2\A_1  &\la6\ra+2\la-6\ra+2\A_1&\\
(4,4)&2&1& \la 2\ra +2\la -6\ra+\A_1&\la6\ra+\la -6\ra+3\A_1 \\
(4,4)&3&1&\la 2\ra +3\la -6\ra   &\la 6\ra+4\A_1& \\
\\
(4,4)&3&0& \U(6)   +\A_2(2) &\la 6\ra+4\A_1& \\
\endmatrix}\endmatrix$
\hskip-2mm
$\matrix\text{\tenrm Table 8C. The case of $o=+$ ($p=1$)}\\
\boxed{\matrix
(r,r_2)&T_1&T_2\\
(6,0)&   \U(3)+2\A_2  & \U+\A_1\\
(4,0)&  \U(3)+\A_2     &\U+\A_2+\A_1& \\
(2,0)&  \U(3)       &\U+3\A_2 \\
\\
(7,1)& \la 6\ra+3\A_2  &\la 2\ra+\A_1&\\
(5,1)&   \U(3)+\A_2+\la-6\ra &\U+2\A_1\\
(5,1)&  \U(3)+\A_2+\A_1 &\U+\la-6\ra+\A_1\\
(3,1)&  \U(3)+\la-6\ra  &\U+\A_2+2\A_1\\
(3,1)& \la 6\ra+\A_2   &\la 2\ra+2\A_2+\A_1 \\
(1,1)&\la 6\ra   &\U\!+\!2\A_2\!+\!\A_1\!+\!\la-6\ra \\
\\
(6,2)&  \U(6)+2\A_2              &\la2\ra+2\A_1& \\
(6,2)&  \U(3)\!+\!\A_2\!+\!\A_1\!+\!\la-6\ra  &\la2\ra+2\A_1 \\
(6,2)&    \la 6\ra+2\A_2+\A_1 &\la 2\ra+\A_1+\la -6\ra \\
(4,2)&     \U(3)+2\la-6\ra   &\U+3\A_1\\
(4,2)& \U(6)+\A_2   &\la2\ra+\A_2+2\A_1& \\
(4,2)&  \la6\ra+\A_2+\la-6\ra    &\la2\ra+\A_2+2\A_1 \\
(4,2)& \la 6\ra+\A_2+\A_1       &\la2\ra\!+\!\la-6\ra\!+\!\A_2\!+\!\A_1\\
(4,2)& \U+\A_2(2)  &\U(3)+\A_2(2)+\A_1 & \\
(2,2)&  \U(6)  &\U\!+\!\A_2\!+\!2\A_1\!+\!\la-6\ra \\
(2,2)&  \la6\ra+\la-6\ra     &\U\!+\!\A_2\!+\!2\A_1\!+\!\la-6\ra \\
(2,2)& \la 6\ra+\A_1  &\la2\ra\!+\!2\A_2\!+\!\A_1\!+\!\la-6\ra \\
\\
(5,3)&  \la 6\ra+\A_2+2\la-6\ra   &\la 2\ra+3\A_1 \\
(5,3)&  \la 6\ra\!+\!\A_2\!+\!\A_1\!+\!\la-6\ra    &\la 2\ra+2\A_1+\la -6\ra \\
(5,3)&  \la 6\ra+\A_2+2\A_1 &\la 2\ra+\A_1+2\la -6\ra &  \\
(3,3)&  \la6\ra+2\la-6\ra  &\U+\la-6\ra+3\A_1 \\
(3,3)& \la 6\ra+\A_1+\la-6\ra   &\la 2\ra+\A_2+2\A_1+\la -6\ra \\
(3,3)& \la 6\ra+2\A_1  & \la2\ra+2\la-6\ra+\A_2+\A_1 \\
\\
(4,4)&  \la 6\ra+3\la-6\ra     &\la 2\ra+4\A_1 \\
(4,4)&  \la6\ra+\A_1+2\la-6\ra   &\la2\ra+3\A_1+\la-6\ra \\
(4,4)& \la 6\ra+\la -6\ra+2\A_1    &\la2\ra+2\A_1+2\la -6\ra \\
(4,4)&\la 6\ra+3\A_1    &\la 2\ra+\A_1+3\la -6\ra  & \\
\\
(4,4)& \U(2)+\A_2(2)    &\la 2\ra+\A_1+3\la -6\ra  & \\
\endmatrix}\endmatrix$
\endinsert
}

\proposition\label{reversion-to-reversion}
If ascending T-pairs $(T_+(c),T_-(c))$ and $(T_+(c'),T_-(c'))$
associated with Zariski curves, $A$ and $A'$, are reversion partners,
then $A$ and $A'$ are reversion partners themselves.
\endproposition

\proof Lemma \ref{anti-isomorphism-criterion}, which reduces the transitivity
of the action of $\Aut(T_-)$ on the subgroups of $\discr_2(T_-)$ anti-isomorphic to $\discr_2(T_+)$ to
the transitivity on the elements of $\discr_2 T_-$ whose square is $-\frac12$, and Proposition \ref{prop-1/2-transitivity},
which establishes the latter transitivity, imply that
$h$ is glued with a reversion root $v$.
Thus, $h$ and $v$ generate an $\U$-summand in $T_-'$. Therefore, $c'=c^\vee$ in the sense of Section \ref{explicit-partners}. Applying Proposition \ref{partner_topology} we find  an example of
reversion partners, $B$ and $B'$,
that have the same homological types as $A$ and $A'$, respectively. By
Theorem \ref{eigenlattices&involutions} and Theorem \ref{deformation-classification}, $B$ is deformation equivalent to $A$ and $B'$ to $A'$.
\endproof

\subsection{Classification of real Zariski sextics}
Given a real Zariski sextic $A$, we associate to it a geometric involution $c\in C(L,\D,h)$ induced
on a K3-lattice $L=H_2(\til Y)$ with a conical $(\D,h)$-decoration by the
ascending real structure, $\conj_{\til Y}\:\til Y\to \til Y$ (see \ref{Covering_K3}, \ref{geometric-involutions}).
The homology type $[c]\in C[L,\D,h]$ depends only on the deformation class $[A]$ of $A$, and the mapping
$[A]\mapsto[c]$ gives by Theorem \ref{deformation-classification} a one-to-one correspondence between the set of deformation classes of real Zariski sextics and the set $C^<[L,\D,h]=\{[c]\in C[L,\D,h]\,|\,c \text{ is ascending}\}$.
 Then we obtain a mapping that sends $[c]$ to the isomorphism type of $(T_+(c),T_-(c))$.

\theorem\label{classification-via-T-pairs}
Associating with a real Zariski sextic $A$ the eigenlattices $(T_+(c),T_-(c))$ of $c\in C^<(L,\D,h)$ defined by the ascending real structure in $\til Y$, we obtain a one-to-one correspondence between the deformation classes
of real Zariski sextics and the ascending  T-pairs listed in the Tables 8A, 8B and 8C.
\endtheorem

\proof
 Theorem \ref{eigenlattices&involutions} shows injectivity of the map from $C^<[L,\D,h]$ to the set of
 isomorphism types of ascending  T-pairs $(T_+(c),T_-(c))$,
 and Theorem \ref{eigenlattices&T-pairs} shows its surjectivity. Theorem
\ref{eigenlattices-enumeration} enumerates the ascending  T-pairs
(they are listed first in Tables 7A-B, and then presented in a different order in Tables 8A-B to fit to our final description of the IDs of real Zariski sextics in Theorem \ref{main-sextics}).
\endproof

\subsection{The IDs of real Zariski sextics}
Our aim now is to obtain from Theorem \ref{classification-via-T-pairs}
 an enumeration of the deformation classes of real Zariski sextics $A$ in terms of their IDs.
The following lemma describes how the characteristics $\ell(A)$, $\chi(\Cal \A_-)$,
$\nu_r(A)$, $o(A)$, and the type of $A$
are expressed in terms of the invariants $r$, $r_2$, $\delta_2$, $p$ and $q$ of $T_+(c)$

\lemma\label{ID-via-lattices} The ID of a real Zariski sextic $A$ determines the values $r$, $r_2$, $\delta_2$, $p$ and $q$ of the eigenlattices $T_\pm\subset L=H_2(\til Y)$ associated with the ascending real structure in $\til Y$.
Conversely, the values  $r$, $r_2$, $\delta_2$, $p$ and $q$ determine $\ell(A)$, $\chi(\Cal \A_n)$, $\nu_r$, the type of $A$ (I or II), and the sign $o(A)$. Namely, if $A(\R)\ne\oo$, then the following relations hold.
\roster\item
$
r_2(T_+)=r_2(T_-)-1=5-\ell(A),\ \ \
r(T_+)=9-r(T_-)=4+\chi(\Cal \A_n),
$\newline
in particular, in the case of code $\a+1\la\b\ra$, we have
$r_2(T_+)=4-(\a+\b)$, $r(T_+)=4+(\b-\a)$, or equivalently
$\a=4-\frac{r+r_2}2$ and $\b=\frac{r-r_2}2$.
\item
The type of $A$ determines $\delta_2(T_+)$ and is determined by it, namely,
 $\delta_2(T_+)=0$ if $A$ has type I, and $\delta_2(T_+)=1$ if type II.
\item
The sign $o(A)$ and the number of real cusps, $2\nu_r(A)$, determine $0\le p\le1$ and $0\le q\le3$, and conversely,
$p$ and $q$ determine $o(A)$ and $2\nu_r$ as follows.
$$
(p,q)=\cases
(0,3-\nu_r(A) ), &\text{\eightrm if \ } o(A)=-,\\
(1,\nu_r(A) ), &\text{\eightrm if \ } o(A)=+,\\
\endcases
\hskip4mm
(o(A),\nu_r)=\cases
(-,3-q), &\text{\eightrm if \ } p=0,\\
(+,q), &\text{\eightrm if \ } p=1.\\
\endcases$$
\endroster

If $A(\R)=\oo$, then $r(T_+)=r_2(T_+)=4$, $\delta_2(T_+)=0$, and $(p,q)$ is $(1,0)$ if $o(A)=-$ and
$(0,3)$ if $o(A)=+$.
\endlemma

\proof
 In the case $A(\R)\ne\oo$, the ascending involution is the M\"obius one, see Lemma \ref{ascending_descending}.
 So, the relations for $r_2(T_\pm)$ in (1) follow from Lemma \ref{defects-are-the-same} and Corollary
 \ref{Mobius-defect}.
The relation between
 $r(T_+)$ and $\chi(\Cal \A_n)$ follows from Lemmas \ref{codes-and-ranks} and \ref{ranks-sing-nonsing}.
Item (2) follows from Proposition \ref{type-via-lattice}, and (3) from Corollary \ref{p-q-of-T}.

In the case of $A(\R)=\oo$, the values of $r$, $r_2$ and $\delta_2$ for $T_+$ are found in
Proposition \ref{exceptional-T-relations}, and the relation for $(p,q)$ are obtained from (3) by alternation of the sign
$o(A)$, because the ascending real structure is non-M\"obius in the case of $A(\R)=\oo$.
\endproof

\lemma\label{empty-K3-lattices}
Assume that
$A$ is a real Zariski sextic and the covering desingularized K3-surface $\til Y$
is endowed with the ascending real structure.
Then the following conditions are equivalent:
\roster\item $A(\R)=\oo$;
\item $\til Y(\R)=\oo$;
\item
the eigenlattice $T_+(c)$ of the involution $c$ induced in $T\subset H_2(\til Y)$ by the real structure
is either $\U(2)+\A_2(2)$, or $\U(6)+\A_2(2)$;
\item
the ascending  T-pair $(T_+(c),T_-(c))$ is either the one in the last row of Table 8B, or the one
in the last row of Table 8C.
\endroster
\endlemma

\proof
For $A(\R)\ne\oo$, we have $\til Y(\R)\ne\oo$ by definition, and for $A(\R)=\oo$ the ascending real structure in
$\til Y$ is the non-M\"obius one (see Lemma \ref{ascending_descending}),
with $\til Y(\R)=\oo$, which shows equivalence of (1) and (2).

Proposition \ref{exceptional-T-relations} says that (1) implies $r(T_+)=r_2(T_+)=4$, and $\delta_2(T_+)=0$.
Let $\discr_3 T_+=p\dip+q\din$, then
$\nu_r=0$ implies that $(p,q)$ is either $(1,0)$, or $(0,3)$ depending on $o(A)$ (see Lemma \ref{ID-via-lattices}).
These characteristics detect precisely two ascending  T-pairs from Tables 7A-B, which are characterized by $T_+(c)$ in (3).
Existence of at least two deformation classes of ``empty'' Zariski sextics $A$ (with $o(A)=\pm$) shows that such T-pairs
really correspond to the case of $\til Y(\R)=\oo$. The aforementioned T-pairs are placed in
the last rows of Tables 8B and 8C.
\endproof

\rk{Remark}
Recall that in the case $A(\R)=\oo$, the central projection of the cubic surface onto $P^2(\R)$ is one-to-one
if $o(A)=-$ and three-to-one if $o(A)=+$.
\boxedR\endrk

\theorem\label{correspondence}
The ID of a real Zariski sextic whose equisingular deformation class is determined by the associated eigenlattices
$(T_+,T_-)$ in one of the Tables 8A, 8B, and 8C, is listed in the same row of the Table 1A, 1B, and 1C respectively.
\endtheorem

\proof
It was already shown in \ref{empty-K3-lattices} that for a real Zariski sextic $A$ with $A(\R)=\oo$ the associated pair $(T_+,T_-)$ is like in the last
rows of Tables 8B and 8C (depending on $o(A)\in\{+,-\}$). So, in what follows we assume that $A(\R)\ne\oo$ and
exclude the last row of Tables 8B-C from consideration.

G.Mikhalkin \cite{M} found 49 complete codes of real Zariski sextics $A$ indicating the topological type of the cubics
$X(\R)$ whose apparent contours they do represent. As follows from Lemma \ref{cubic-via-code}, $\chi(X(\R))$ together with the complete code of $A$ detect $o(A)$. Thus,
using Lemma \ref{ID-via-lattices}, one can determine the invariants $r$, $r_2$, $p$ and $q$ for the corresponding
lattices $T_\pm$, while $\delta_2(T_+)$ remains unknown, since the type of $A$ was not determined by Mikhalkin.
 Reviewing Tables 7A-B we can observe that the above invariants $r$, $r_2$, $p$ and $q$ distinguish all the
rows there except six pairs of rows. In each of these pairs, lattices $T_+$ are distinguished by the value of $\delta_2$.

The ascending  T-pairs $(T_+,T_-)$ whose invariants  $r$, $r_2$, $p$ and $q$ do not match the real Zariski sextics found by Mikhalkin, appear
in the rows 4 and 6 of Table 8A, rows 2, 3, 7, 8 of Table 8B, and rows 1, 2, 5, 6, 7, 13 and 20 of Table 8C.

If we know the ID of a real Zariski sextic $A$ associated with $(T_+,T_-)$, then we may determine the ID of
 the one associated with the partner of $(T_+,T_-)$ using Proposition \ref{reversion-to-reversion},
 if there is a non-empty oval in $A(\R)$, see Corollary \ref{reversion-cor}.
Tables 8B and 8C are arranged to place the partner pairs $(T_+,T_-)$ in the same rows.
 The sextics $A$ represented by $(T_+,T_-)$ in rows 1, 5, 6, 13, 20 of Table 8B and rows 3, 8 of Table 8C turn out to have a non-empty oval, so it remains to determine the IDs of $A$ represented by $(T_+,T_-)$ from Table 8A, rows 4 and 6 (since they have no partners), and from Tables 8B and 8C, rows 2 and 7 (since Mikhalkin's examples are missing for the both partners).

In the remaining part of the proof we may suppose that the simple code $A$ looks like $\a\+1\la\b\ra$, since
the case of null-code was considered in the beginning and the two cases of 3-nest codes (with different values of
$o(A)$) are represented in the row 17 of Tables 8B and 8C, which requires no further analysis
(it follows from Lemmas \ref{nested-no-cusps} and \ref{type-lemma} that these are the only IDs with the 3-nest code).

 Lemma \ref{ID-via-lattices} shows how to reconstruct from $(T_+,T_-)$ the type of $A$, $\nu_r(A)$, and $o(A)$. So, it remains only to determine the distribution of the real cusps on the ovals of $A(\R)$ to reconstruct
 its ID.

\lemma\label{missing-cases}
The real Zariski sextics $A$ characterized by the eigenlattices $(T_+,T_-)$ in rows 4 and 6 of Table 8A and rows 2 and 7 in each of Tables 8B--8C have IDs indicated in the corresponding rows of Tables 1A and 1B--1C respectively.
\endlemma

\proof
In all the cases indicated $A$ has type II, since $\delta_2(T_+)=1$.
 In the case of row 6 of Table 8A, Lemma \ref{ID-via-lattices} says that $A$ has no real cusps and thus,
the complete code coincides with the simple one. This Lemma says also that this code is $1\la1\ra$, and determines
the sign $o(A)=+$ (since $p=1$), giving the ID indicated in row 6 of Table 1A.

For row 4 of Table 8A, we get $\nu_r=1$, the simple code $1\la2\ra$, and $o(A)=-$. Lemma \ref{cusp-parity} implies that the both real cusps should lie on the same oval, by Lemma \ref{cusp-direction}, they must be outward cusps if lying
on the ambient oval and inward cusps in lying on an internal one. But the latter case is forbidden by
Lemma \ref{polycuspidal-ovals}, so the complete code of $A$ must be $1_1\la2\ra$.

In the case of row 2 of Table 8B, we similarly find $\nu_r=2$, the sign $o(A)=-$, and the simple code $2\+1\la2\ra$.
Applying Lemmas \ref{cusp-direction} and \ref{polycuspidal-ovals}, we conclude that
the two pairs of outward cusps must be distributed among the external ovals and the ambient oval.
Lemmas \ref{few-smooth-ovals} and \ref{one-smooth-oval} exclude a possibility that one of the external ovals
is smooth, thus both of them have a pair of cusps and the complete code is $2_1\+1\la2\ra$.

In the case of row 7 of Table 8B, we find $\nu_r=2$, $o(A)=-$, and the simple code $1\+1\la2\ra$.
 As above, we conclude that there should be two pairs of outward cusps distributed among the ambient and the external ovals. Lemma \ref{one-smooth-oval} exclude possibility that all the cusps lie on the ambient oval, while Lemma
\ref{polycuspidal-ovals} exclude possibility that all the cusps lie on the exterior one. Thus, the complete code must be $1_1\+1_1\la2\ra$.

The rows 2 and 7 of Table 8C represent the reverse partners for the same rows of Table 8B, and thus, the corresponding
curves $A(\R)$ must be in reverse positions by Proposition \ref{reversion-to-reversion}.
Since not all the ovals are empty, the complete code of the reverse partners can be found by the rule in Corollary \ref{reversion-cor}, namely, reversion of  $2_1\+1\la1\ra$ gives
$1\+1\la2_1\ra$
and reversion of $1_1\+1_1\la2\ra$ gives $2\+1_{-1}\la1_1\ra$, like indicated in Tables 1C.
\endproof

 Finally, it remains to analyze the six pairs of rows which differ only by the value of $\delta_2(T_+)$, namely, rows 10 with 11, rows 14 with 15, and rows 18 with 19 in Table 8B, and the same pairs of rows in Table 8C.

\lemma
Each of the three pairs of rows, 10 with 11, 14 with 15, 18 with 19 in Table 8B,
represent a pair of deformation classes of real Zariski sextics $A$ which have the same complete codes.
The same is true for the same pairs of rows in Table 8C. The corresponding complete codes are like indicated in
the corresponding rows of Tables 1B and 1C.
\endlemma

 \proof
Using Lemma \ref{ID-via-lattices} like in the proof of Lemma \ref{missing-cases}, we conclude that
the complete codes of $A$ for the indicated pairs of rows may differ only by the distribution of cusps
on the ovals, since the values of $r$, $r_2$, $p$ and $q$ for $T_+$ in each of the pairs of rows is the same.

Rows 10 and 11 of Table 8B give both the simple code $3$, with $\nu_r=3$ and $o(A)=-$, which implies that
there are three pairs of outward cusps.
None of these three ovals can be smooth by Corollary \ref{empty-ovals-cor},
so, both rows should give complete code
$3_1$.

Rows 14 and 15 give both the simple code $1\+1\la1\ra$, with $\nu_r=2$, and $o(A)=-$.
The sign $o(A)$ together with  Lemmas \ref{cusp-direction} and \ref{polycuspidal-ovals} imply that the
cusps are outward (cf. the proof of Lemma \ref{missing-cases}). An external oval can be neither smooth nor 4-cuspidal by Lemmas \ref{one-smooth-oval} and \ref{polycuspidal-ovals} respectively. So, the complete code for the both rows is
 $1_1\+1_1\la1\ra$.

 Rows 18 and 19 give the simple code $1\la2\ra$ with $\nu_r=1$ and $o(A)=-$. This implies that
there is one pair of cusps which must lie on the ambient oval, since
 Lemmas \ref{cusp-direction} and \ref{polycuspidal-ovals} exclude other possibilities.

The same rows in Table 8C represent the reversion partners for the ones considered above, and thus, their complete code is obtained by reversion, according to the rules in Corollary \ref{reversion-cor}.

 As the result, we can deduce that in each of these six pairs, the cusps must be distributed on the ovals in the same way, and thus, the corresponding real Zariski sextics $A$ must have the IDs like indicated
in the 1B and 1C.  \qed  \qed

\subsection{Proof of Theorem \ref{main-sextics}}\label{proof_main}
Theorem \ref{classification-via-T-pairs} together with the correspondence established in Theorem
\ref{correspondence} imply Theorem \ref{main-sextics}.
\qed

\section{Concluding Remarks}\label{concluding}
\subsection{Purely real statements}
Although the ID of a real Zariski sextic refers to
its complex point set, it is sufficient to look at the Tables 1A-C to conclude that
in the majority of cases the deformation classes are determined only by
the complete codes, which are purely real invariants.
There are however a few exceptions. The first group of exceptions is given by real
Zariski sextics with the complete codes $\emptyset$, $1$, $1\langle 1\rangle$, and
$1\langle 1\langle 1\rangle\rangle$. In each of these cases there are 2 deformation classes.
But they can be distinguished by enhancing the code
with the invariant $o(A)=\pm$, or in other words, marking the domain where the projection
of the cubic surface is three-to-one and thus expressing
the classification in terms of ``purely real'' data.

The other group of exceptions contains
6 complete codes, each one again representing 2
deformation classes, but this time
the classes in each pair
differ by the types, I or II, of the sextics. The only remedy we can
suggest to distinguish the types in purely real terms
is to use real
lines passing through a pair of ovals of the sextic. In one case
such lines separate the cusps on the third oval as it is shown on
Figure \ref{6-pairs}, and in the other, they do not.

\midinsert\figure
{The six pairs of codes
that differ only by their type (the bottom pairs are the partners of the top pairs)}
 \label{6-pairs}
\endfigure
\hskip15mm\centerline{\epsfbox{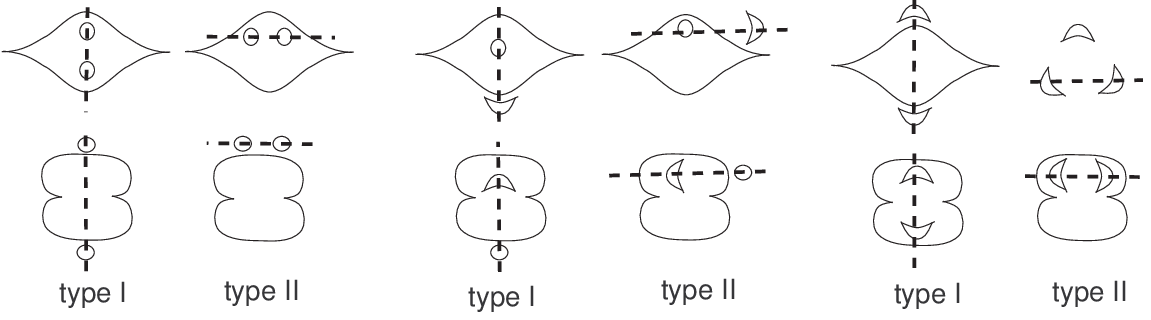} }
\endinsert
\midinsert\figure{Construction of the curves on Figure \ref{6-pairs} }
\label{I-IIconstructions}
\endfigure
\centerline{ \epsfbox{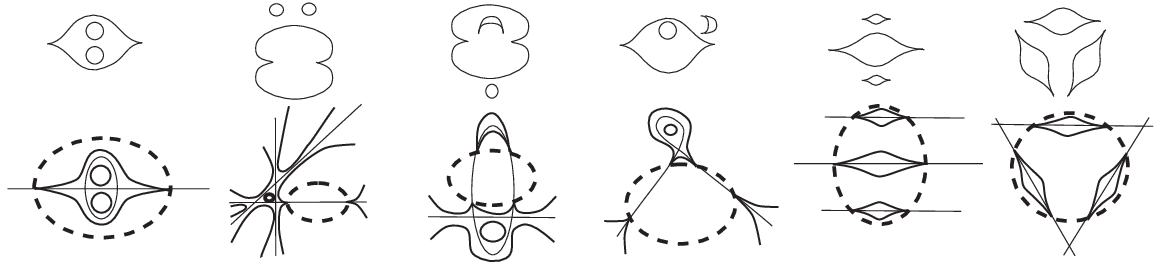} }
\endinsert

It would be interesting to find a conceptual explanation to this observation.
 Our current proof that the types are as indicated on Figure \ref{6-pairs}, is based on the construction of these curves shown on Figure \ref{I-IIconstructions}  (the construction in question consists
in a small perturbation, $p^2+\varepsilon q^3$,
where $p$ is a reducible cubic curve and $q$ is a conic).
Since sextics-partners have the same type, it is sufficient to
consider one representative from each pair of partners.

Recall that according to our definitions, a real Zariski sextic
is the apparent contour of a generic projection of a nonsingular real cubic surface, where
the latter is supposed to
have neither real nor complex singular points and a generic apparent contour is supposed to
have no singular points other than ordinary (real or complex) cusps.
Therefore, it is natural to ask how
the deformation classification of generic apparent contours may change if we allow for
cubic surfaces to have non-real singular points and for Zariski sextics to have
non-real singularities other than cusps.
 The answer is straightforward, it says that the deformation
classification does not change. (However, to
make the formulation and solution of the problem
depend only on the real locus of Zariski sextics in question, it
is better to exclude the case of sextics with empty real locus; in
all the other cases, the real locus is sufficient for extending the
real central projection correspondence stated in Proposition
\ref{real-central correspondence} to this new setting.).

\subsection{Transversal pairs of conic and cubic}
As is known, see Proposition \ref{torus-type}, a Zariski sextic is uniquely defined by a pair
of  homogeneous polynomials of degree $2$ and $3$
defining a conic and a cubic intersecting transversely.
Thus, it may be
worth to ask about deformation classification of such pairs of polynomials.
Understanding this question literally, that is as a classification of pairs $p,q$ where $p$, $q$ are a conic and a cubic (not necessarily nonsingular) intersecting transversally each other, one
 easily gets the
following answer. Over $\C$, there is only one deformation class.
And over $\R$, there are 4 classes; they
are  distinguished just by the number of real
intersections points: 6, 4, 2, or none.

Curiously
enough, the latter result being so far from the principal object of our investigation is however also related to classification of
cubic surfaces. Namely, if we impose the assumption of non-singularity on the conic,
then aforesaid transversal pairs describe cubic
surfaces with one node. Over $\C$ the pairs
of  homogeneous polynomials of degree $2$ and $3$
defining a non-singular conic and a (possibly singular) cubic intersecting its transversely
form a single deformation class,
while over $\R$ we get 10 deformation classes, more
than before, because of the well defined sign of the degree 2 polynomial on the nonorientable part of the complement of the
set of its zeros in the real projective plane (now, the sign can not be changed, because we have forbidden to the conic to become singular) and
a possibility to have 0 intersection points in two ways, with empty and with non empty conic.
As a consequence,
real cubic surfaces with one node form 5 deformation classes (twice less, since reversing of the sign of the polynomial
of the cubic surface leads to the change of the sign of the polynomial defining the conic); the
fact which was observed exactly in this way by F.~Klein, who used
it as a step to his classification of real nonsingular cubic surfaces up to deformation.

A different, but somehow related and typical
problem is the classification of pairs of transversally intersecting nonsingular curves. As was shown by G.~Polotovsky
\cite{Pol}, in the case of a real conic and a real
cubic one obtains 25 deformation classes: each deformation class is determined by the topology of the
arrangement of real points in the real projective plane; there are 7 extremal classes shown in Figure \ref{2+3}
and the other ones are obtained from them by
two moves: erasing an oval (of the cubic or of the conic)
containing no intersection points, and shifting a piece of curve containing
a pair of consecutive (both on the conic and the cubic) intersection points,
so that these two intersections disappear.

\midinsert\centerline{\epsfbox{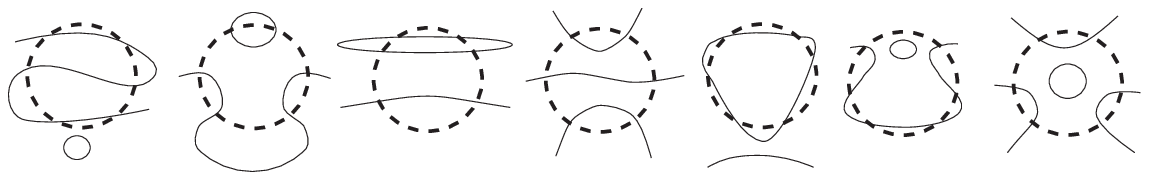} }\figure{The extremal mutual positions of a cubic and a conic} \label{2+3}
\endfigure
\endinsert

The both classification problems for the pairs formed by a conic and a cubic
are different from the problem of classification of Zariski sextics by an evident, although a deep, reason. In
the case of Zariski curves there is a more subtle non-singularity condition:
the sextic $p^2+q^3$ should not have singular points other than the intersection points of the conic $q$ with the cubic $p$.

\subsection{Ordering of IDs}
To characterize the set of complete codes of Zariski sextics one
can give a short list of the extremal ones, from which all the others
can be obtained by certain simplifying moves.
 The extremal codes include seven M-codes, namely, $1\la4\ra$,
 $\a_1\+1\la\b\ra$, and $\a\+1\la\b_1\ra$, $1\le \a,\b\le3$, $\a+\b=4$,
 the 3-nest code $1\la1\la1\ra\ra$, and the code $1\la1_1\+1\ra$.
 The simplifying moves are: \roster\item cancelation of an empty smooth oval (for instance,
$2\+1\la2_1\ra$ gives $1\+1\la 2_1\ra$),
 \item cancelation of an empty oval with a pair of
outward cusps (for instance, $2\+1\la2_1\ra$ gives $2\+1\la 1_1\ra$),
 \item
fusion of an empty oval that has a pair of outward cusps with a
principal oval  (the latter was defined in Section 2.6)
as is shown on Figure \ref{simplification}
(for instance, $2_1\+1\la2\ra$ gives $1_1\+1_1\la 2\ra$, while $2\+1\la2_1\ra$
gives $2\+1_{-1}\la 1_1\ra$).
\endroster
 \midinsert\centerline{\epsfbox{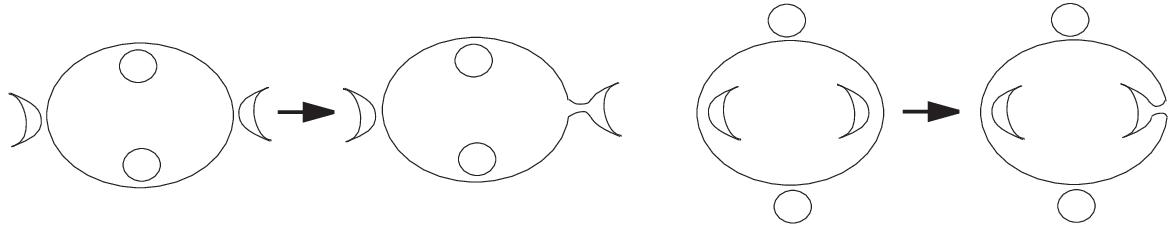} }\figure{Fusion moves}
\label{simplification}
\endfigure
\endinsert
All the complete codes of Zariski sextics can be obtained by such
moves from the extremal ones. However, one can obtain also a few extra ones,
which are not the codes of Zariski sextics. Namely,
$1\+1_{-3}$, $2\+1_{-2}$, and $3\+1_{-1}$, are obtained by
fusion-move applied to all the internal ovals of
$\a\+1\la(4-\a)_1\ra$, $\a=1,2,3$.
These are the only three
exceptions: if we apply simplification moves to any of the extremal codes
in any other way, then we necessarily obtain again a code of a Zariski sextic.

\subsection{Nonsingular partners} The partner duality that played an important role in our
classification of real Zariski sextics (in matching the homological types against the IDs)
exists as well in the case of nonsingular real sextics.
Implicitly, it appears already
in Hilbert's sixteenth problem statement.
Furthermore, it performed a dramatic and somehow decisive role in Gudkov's classification of real nonsingular sextics. It was the subject of his habilitation thesis, and when he had shown the preliminary version to one of his "thesis referees", V.V.~Morozov
 (Professor at the Kazan University), the latter have objected the resulting classification
exactly because of a small irregularity with respect to the "reversion symmetry" in it.
It is by repairing this asymmetry that Gudkov has come to his final result. Such a reversion symmetry revealed itself forcefully again in Rokhlin's and Nikulin's treatment
of deformation classes of  real nonsingular sextics. Up to the best of our knowledge, a conceptual explanation of this partner duality/reversion symmetry
was never explicitly presented in the literature. In fact, such an explanation for nonsingular sextics is literarily the same
as for Zariski sextics, both in the lattice-arithmetical and in geometric terms.
Namely, for each deformation class with one exception, the eigenlattice $L_-$
contains an $\U$-summand and the lattice-arithmetical form of the partner duality
consists in transferring the $\U$-summand to the opposite eigenlattice $L_+$
and then exchanging of the eigenlattices. In geometrical
terms it means that each partner in a partner pair can be deformed to a triple conic, near the triple conic the family looks as
$Q^3+tf_2Q^2+t^2f_4Q+f_6=0$ (cf., Introduction), and switching of the sign of $t$ (which corresponds to passing through the triple conic) replaces the curves from one deformation class by the curves from its partner class. The only exceptional deformation class having no partner class
is formed by real sextics of type I with the code $1\langle 4\rangle$.

\subsection{Promiscuity}
In the case of real sextics with arbitrary singularities,
an analogue of the partnership relation can be also defined, but it looks at first glance rather like
a sort of promiscuity, at least from one side.
 For instance,  in the case of sextics with a node (cf., \cite{Iten}),
there is a well-defined partnership map (neither injective nor surjective)
that transforms the deformation classes of real sextics with an {\it internal node} to the ones
 with an {\it external node} as is illustrated on Figure \ref{partnership-map}.
\midinsert
\centerline{ \epsfbox{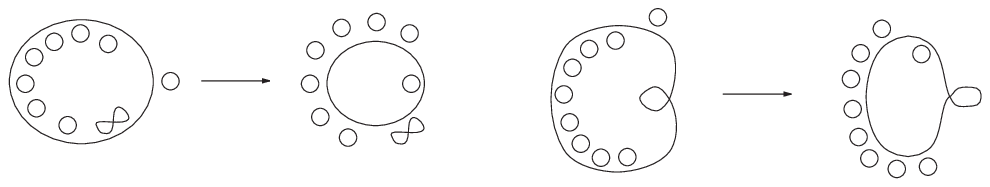} }
\figure{The partnership map for real nodal sextics}\endfigure\label{partnership-map}
\endinsert

{\sl
\hskip1.8in If we shadows have offended,

\hskip1.8in Think but this, and all is mended,

\hskip1.8in That you have but slumber'd here

\hskip1.8in While these visions did appear.

\hskip1.8in And this weak and idle theme,

\hskip1.8in No more yielding but a dream,

\hskip1.8in Gentles, do not reprehend:

\hskip1.8in if you pardon, we will mend ...
}
\vskip3mm
\hskip2in  ``Midsummer Night Dream'', W. Shakespeare

\enddocument